\documentclass[12pt,a4paper]{article}
\usepackage{amsmath,amssymb}
\usepackage{latexsym}
\usepackage{graphicx}
\usepackage{color}
\usepackage{multirow}
\usepackage{indentfirst}
\usepackage{subfigure}
\usepackage{mathrsfs}
\usepackage[pdftex]{hyperref} 

\newtheorem{exam}{\hspace{6mm}Example}[section]

\begin{document}
\baselineskip=1.2pc

\begin{center}
{\large \bf  A quasi-Lagrangian moving mesh discontinuous Galerkin method for hyperbolic conservation laws}
\end{center}

\centerline{Dongmi Luo%
\footnote{School of Mathematical Sciences and Fujian Provincial
Key Laboratory of Mathematical Modeling and High-Performance
Scientific Computing, Xiamen University, Xiamen, Fujian 361005, China. E-mail: dongmiluo@stu.xmu.edu.cn.},
Weizhang Huang%
\footnote{Department of Mathematics, University of Kansas, Lawrence, Kansas 66045, U.S.A.
E-mail: whuang@ku.edu.},
Jianxian Qiu%
\footnote{School of Mathematical Sciences and Fujian Provincial
Key Laboratory of Mathematical Modeling and High-Performance
Scientific Computing, Xiamen University, Xiamen, Fujian 361005, China. E-mail: jxqiu@xmu.edu.cn.}
}

\vspace{20pt}

\begin{abstract}
A moving mesh discontinuous Galerkin method is presented for the numerical solution of hyperbolic conservation laws. The method is a combination of the discontinuous Galerkin method and the mesh movement strategy which is based on the moving mesh partial differential equation approach and moves the mesh continuously in time and orderly in space. It discretizes hyperbolic conservation laws on moving meshes in the quasi-Lagrangian fashion with which the mesh movement is treated continuously and no interpolation is needed for physical variables from the old mesh to the new one. Two convection terms are induced by the mesh movement and their discretization is incorporated naturally in the DG formulation. Numerical results for a selection of one- and two-dimensional scalar and system conservation laws are presented. It is shown that the moving mesh DG method achieves the theoretically predicted order of convergence for problems with smooth solutions and is able to capture shocks and concentrate mesh points in non-smooth regions. Its advantage over uniform meshes and its insensitiveness to mesh smoothness are also demonstrated.
\end{abstract}


{\bf AMS(MOS) subject classification:} 65M60, 35L65
\pagenumbering{arabic}

\newpage

\section{Introduction}
\label{sec1}
\setcounter{equation}{0}
\setcounter{figure}{0}
\setcounter{table}{0}

We consider the numerical solution of hyperbolic conservation laws in the form
\begin{equation}
\label{eq1}
\left\{
\begin{array}{l}
u_t + \nabla\cdot F=0,\\
u(\textbf{x},0)=u_0(\textbf{x}),
\end{array}
\right.
\end{equation}
where $\textbf{x}\in \Omega \subset \mathbb{R}^d$  $(d = 1, 2)$, $\Omega$ is a bounded domain, $F=(f_1(u),\cdots,f_d(u))$, and $u$, $f_1(u),\; \cdots,\; f_d(u)$ are either scalars or vectors.
The major difficulty in solving nonlinear  hyperbolic conservation laws in (\ref{eq1}) is that the solution can develop discontinuities even if the initial condition is smooth. The discontinuous Galerkin (DG) method is an increasingly popular approach to solve the equations.
The DG method was first introduced by Reed and Hill \cite{EH00} for solving linear hyperbolic problems associated
with neutron transfer.  Then a major development of the DG method for time-dependent nonlinear  hyperbolic conservation laws was carried out by Cockburn et al. in a series of papers \cite{EH03, EH02, EH01, EH04}. The DG method can capture the weak discontinuity without any modification. But the nonlinear limiters must be applied to control the spurious oscillations in the numerical solution for strong shocks. One type of these limiters is based on the slope methodology such as the minmod type limiters \cite{EH03, EH02, EH01, EH04}. These limiters are effective in controlling oscillations. However, the accuracy of the DG method may decrease if they are mistakenly used in smooth regions. Another type of limiter is based on the weighted essentially non-oscillatory (WENO) methodology \cite{EH05, EH06}, which can achieve both high-order accuracy and non-oscillatory properties. The WENO limiters in \cite{EH13,EH14} and the Hermite WENO (HWENO) limiters in \cite{EH09, EH10, EH11, EH12} belong to this type limiter and require a wide stencil especially for the high-order accuracy. Recently, a simple WENO limiter \cite{EH15, EH16} and compact HWENO limiter \cite{EH17, EH18} for Runge-Kutta DG (RKDG) were developed. The key idea of these limiters is to reconstruct the whole polynomial in the target element instead of point values or moments.

The physical phenomena in a variety of fields may develop dynamically singular solutions, such as shock waves and detonation waves. If we use the globally uniform mesh method, the computation can become prohibitively expensive when dealing with two or higher dimensional systems. Then adaptive mesh methods that can increase the accuracy of the numerical approximations and decrease the computational cost are in critical need. In general, there are three types of mesh adaptive methods. The first one is $h$-methods which are widely used and generate a new mesh by adding or removing points to an old mesh. The second one is $p$-methods with which the order of the polynomial is increased or decreased from place to place according to the solution error. The last one is $r$-methods which are also called moving mesh methods and relocate the mesh points while keeping the total number of mesh points and the mesh connectivity unchanged. A number of moving mesh methods have been developed in the past; e.g., see
the books or review articles \cite{Bai94a,Baines-2011,BHR09,EH22,Tan05} and references therein
and some recent applications \cite{wise2017,Fei-Zhang,zhang2017}.


One of the advantages of DG method is that its numerical solution is discontinuous at edges of mesh elements
and flexible for mesh adaptation strategies. 
There exist a few research works in this aspect. A combination of the $hp$-method and DG method was developed in \cite{EH19} for the hyperbolic conservation laws. 
Li and Tang \cite{EH07} solved two-dimensional conservation laws using rezoning moving mesh methods
where the physical variables are interpolated from the old mesh to the new one using conservative interpolation
schemes.  The methods are shown to perform well although it is unclear that the method can be high order. Machenzie and Nicola \cite{EH20} solved the Halmiton-Jacobi equations by the DG method using a moving mesh method based on the moving mesh partial differential equation (MMPDE) strategy. Uzunca, Karas\"{o}zen, and K\"{u}\c{c}\"{u}kseyhan \cite{EH29} employed the moving mesh symmetric interior penalty Galerkin (SIPG) method to solve PDEs with traveling waves.

In this paper a moving mesh DG method is proposed for the numerical solution of hyperbolic conservation laws. The method is different from those \cite{EH07, EH33} in mesh movement strategy and discretization of physical equations. We use here the MMPDE moving mesh strategy \cite{EH39, EH28, EH22, ZCHQ} which employs a meshing functional based on the equidistribution and alignment conditions and a matrix-valued function that provides the information needed
to control the size, shape, and orientation of mesh elements from place to place. 
Moreover, the newly developed discretization of the MMPDE \cite{EH23} is used here, which makes the implementation more easier and much more reliable since there is a theoretical guarantee for mesh nonsingularity \cite{EH37}.
In contrast to the moving mesh finite difference WENO method \cite{EH08} where the smoothness of the mesh is extremely important for accuracy, we find that the moving mesh DG method presented in this work
is not sensitive to the smoothness of the mesh.
Furthermore, the method discretizes hyperbolic conservation laws on moving meshes in the quasi-Lagrangian fashion with which the mesh movement is treated continuously and no interpolation is needed for physical variables
from the old mesh to the new one. Numerical results for a selection of one- and two-dimensional scalar and system conservation laws will be presented to demonstrate that the moving mesh DG method achieves the theoretically predicted order of convergence for problems with smooth solutions and is able to capture shocks and concentrate mesh points in non-smooth regions. The sensitivity of the accuracy of the method to mesh smooth will also be discussed.


The organization of the paper is as follows. In Section \ref{sec2},
the DG method, identification of troubled cells, and limiters on moving meshes are described in detail.
The MMPDE moving mesh strategy is described in Section~\ref{sec:mmpde}.
In Section \ref{sec:numerics}, one- and two-dimensional numerical examples are presented to demonstrate
the accuracy and the mesh adaptation capability of the scheme. Conclusions are drawn in Section~\ref{sec4}.

\section{DG method on the moving mesh}
\label{sec2}
\setcounter{equation}{0}
\setcounter{figure}{0}
\setcounter{table}{0}

In this section we describe a Runge-Kutta DG (RKDG) method for the numerical
solution of conservation laws in the form of (\ref{eq1}) on a moving triangular mesh
and pay special attention to limiters and the identification of troubled cells.
We will discuss the generation of adaptive moving meshes in the next section.

\subsection{The RKDG method on moving meshes}

For the moment we assume that 
a moving mesh $\mathscr{T}_h(t)$ for the domain $\Omega$ is given at time instants
\[
t_0=0 < t_1 < \cdots t_n < t_{n+1} < \cdots \le T.
\]
The appearances are denoted by $\mathscr{T}_h^n$, $n=0, 1, ...$.
Since they belong to the same mesh, they have the same number of elements $(N)$ and vertices $(N_v)$
and the same connectivity, and differ only in the location of the vertices.
Denote the coordinates of the vertex of $\mathscr{T}_h^n$ by $\textbf{x}_j^n,\; j=1,2,\cdots,N_v$.
Between any time interval $t_n$ and $t_{n+1}$, the coordinates and velocities of the vertices of
the mesh are defined as
 \begin{align}
 &\textbf{x}_j(t) =\frac{t-t_n}{\Delta t_n}\textbf{x}_j^{n+1}+\frac{t_{n+1}-t}{\Delta t_n}\textbf{x}_j^n,\quad j=1,2,\cdots,N_v, \notag\\
 &\dot{\textbf{x}}_j(t)=\frac{\textbf{x}_j^{n+1}-\textbf{x}_j^n}{\Delta t_n},\quad j=1,2,\cdots,N_v.\notag
 \end {align}
 where $\Delta t_n = t_{n+1}-t_{n}$.
The DG finite element space is defined as
\[
V_h^k=\{p(\textbf{x},t):p|_K\in P^k(K),\forall K\in \mathscr{T}_h(t)\},
\]
where $P^k(K)$ is the space of polynomials of degree $\leq k$ defined on $K$. 
Notice that $P^k(K)$ can be expressed as 
\[
P^k(K) = \text{span}\{\phi_1(\textbf{x},t),\cdots,\phi_L(\textbf{x},t)\},
\]
where $L=k+1$ for one dimensional case, $L=\frac{(k+1)(k+2)}{2}$ for two dimensional case, 
and $\{\phi_1(\textbf{x},t),\; \cdots,\; \phi_L(\textbf{x},t)\}$ is a basis of $P^k(K)$.
Notice that the time dependence of the basis functions comes from the time dependence
of the location of the vertices. 

The semi-discrete DG approximation for (\ref{eq1}) is to
 find $u_h(\cdot,t)\in V_h^k$, $t \in (0, T]$ such that
 \begin{align}
  \label{v1}
  \int_{K} \frac{\partial u_h}{\partial t} vdxdy+\int_{\partial K} F\cdot \vec{n}vds-\int_{K} F\cdot \nabla vdxdy=0,
\qquad  \forall v\in P^k(K), \quad \forall K \in \mathscr{T}_h(t) 
\end {align}
where $\vec{n}=(n_x,n_y)$ is the outward unit normal vector of the triangular boundary $\partial K$. 
Expressing $u_h$ as
 \[
 u_h(\textbf{x},t)=\sum_{j=1}^L u_j^K(t)\phi_j(\textbf{x},t),\quad \forall \textbf{x}\in K
 \]
 we can find its time derivative as
\[
\frac{\partial u_h}{\partial t}=\sum_{j=1}^L \frac{du^K_j(t)}{dt}\phi_j(\textbf{x},t)+\sum_{j=1}^L u^K_j(t)\frac{\partial\phi_j(\textbf{x},t)}{\partial t}.
\]
It is not difficult to show \cite{EH40} that
\[
\frac{\partial \phi_j(\textbf{x},t)}{\partial t}=-\nabla{\phi_j(\textbf{x},t)}\cdot\dot X(\textbf{x},t),
\]
where $\dot X(\textbf{x},t)$ is the linear interpolation of the nodal mesh speed in $K$, viz.,
\[
\dot X(\textbf{x},t)=\sum_{l=1}^3\dot{\textbf{x}}_l\Phi_l(\textbf{x},t)
\]
where $\dot{\textbf{x}}_l$, $l = 1, 2, 3$ are the nodal mesh speed and  $\Phi_l(\textbf{x},t),\; l=1,2,3$
are the linear basis functions.
Then, we get
\[
\frac{\partial u_h}{\partial t}=\sum_{j=1}^N\frac{du^K_j(t)}{dt}\phi_j(\textbf{x},t)-\nabla{u_h(\textbf{x},t)}\cdot\dot X(\textbf{x},t).
\]
Inserting this to (\ref{v1}) yields
 \begin{align*}
  \int_{K} \left (\sum_{j=1}^N\frac{du^K_j(t)}{dt}\phi_j(\textbf{x},t)-\nabla{u_h(\textbf{x},t)}\cdot\dot X(\textbf{x},t)\right )vdxdy+\int_{\partial K} F\cdot \vec{n}vds-\int_{K} F\cdot\nabla vdxdy=0. 
\end{align*}
Integrating the second term by parts, we obtain
\begin{align}
  \label{v3}
  \int_{K} \sum_{j=1}^N\frac{du^K_j(t)}{dt}\phi_j vdxdy+\int_{\partial K} (F-u_h\dot X)\cdot\vec{n}vds
  -\int_{K} \left ( F\cdot \nabla v-u_h\nabla\cdot(\dot Xv)\right )dxdy=0. 
\end{align}
From the above equation we can see that two advection terms are induced by the mesh movement. Denote 
\[
H(u_h) \equiv (F(u_h)-u_h\dot X(\textbf{x},t))\cdot \vec{n}, \qquad 
H_1(u_h)\equiv F\cdot\nabla v-u_h\nabla\cdot(\dot Xv).
\]
{ Applying} the Gauss quadrature rule to the second and third terms in (\ref{v3}), we get
\begin{align}
 & \int_{\partial K} H(u_h)vds \approx \sum\limits_{e}
 \sum\limits_{G_e} H(u_h(\textbf{x}_{G_e}))v(\textbf{x}_{G_e})w_{G_e}|e|, \notag \\
& \int_{K} H_1(u_h)dxdy \approx \sum\limits_{G}H_1(u_h(\textbf{x}_{G}))w_G|K|, \notag
 \end {align}
where $e$ represents the edges of the element $K$,  $|K|$ is the volume of the element $K$,
and $\textbf{x}_{G}$ and $\textbf{x}_{G_e}$  represent the Gaussian points on $K$ and $e$, respectively.
The summations $\sum_e$, $\sum_{G}$, and $\sum_{G_e}$ are taken over the edges of $\partial K$,
Gauss points on $K$, and Gauss points on $e$, respectively.
Replacing the flux $H$ by the numerical flux $\Hat H$, we obtain
\begin{equation}
\label{vv}
\begin{cases}
\int_{K} \sum\limits_{j=1}^N\frac{du^K_j(t)}{dt}\phi_j(\textbf{x},t) vdxdy+\sum\limits_e\sum\limits_{G_e} \Hat Hv(\textbf{x}_{G_e}^{int})w_{G_e}|e| \\
\qquad \qquad \qquad \qquad \qquad
-\sum\limits_{G}H_1(u_h(\textbf{x}_{G}))w_G|K|=0,\quad \forall K\in \mathscr{T}_h,\quad v\in V_h^k\\
\int_K (u_h(\textbf{x},0)-u_0(\textbf{x}))vdxdy=0,\quad \forall K\in \mathscr{T}_h,\quad v\in V_h^k
\end{cases}
\end{equation}
where the numerical flux has the form $\Hat H=\Hat H(u(\textbf{x}_{G_e}^{int}),u(\textbf{x}_{G_e}^{ext}))$
and $u(\textbf{x}_{G_e}^{int})$ and $u(\textbf{x}_{G_e}^{ext})$ are defined as the values from the interior
and exterior of $K$, respectively, i.e.,
\[
u(\textbf{x}_{G_e}^{int},t)=\lim\limits_{\textbf{x}\to\textbf{x}_{G_e},\textbf{x}\in K}u_h(\textbf{x},t),
\quad
u(\textbf{x}_{G_e}^{ext},t)=\lim\limits_{\textbf{x}\to\textbf{x}_{G_e},\textbf{x}\notin K}u_h(\textbf{x},t).
\]
The numerical flux $\Hat H(a,b)$ is required to satisfy the following conditions.
\begin{enumerate}
\item[(i)] $\Hat H(a,b)$ is Lipschitz continuous in both arguments $a$ and $b$;
\item[(ii)] $\Hat H$ is consistent with $H(u)$, namely, $\hat H(u,u)=H(u)$;
\end{enumerate}
In this work, we use the local Lax-Friedrichs flux,
\[
\Hat H(a,b)=\frac{1}{2}\left (\frac{}{} H(a)+H(b)-\alpha_{e,K} (b-a)\frac{}{}\right ),
\]
where $\alpha_{e,K}$ is the numerical viscosity constant taken as the largest eigenvalues in magnitude of
\[
\frac{\partial}{\partial u}(F(\bar u_K)-\bar u_K\dot X(\textbf{x},t))\cdot \vec{n}, \qquad
\frac{\partial}{\partial u}(F(\bar u_{K'})-\bar u_{K'}\dot X(\textbf{x},t))\cdot \vec{n},
\]
where $K$ and $K'$ are the elements sharing the common edge $e$ and
\begin{equation}
\label{average-1}
\bar u_K=\frac{1}{|K|}\int_K u_h dxdy,\quad \bar u_{K'}=\frac{1}{|K'|}\int_{K'} u_h dxdy .
\end{equation}


Finally the semi-discrete scheme (\ref{vv}) is discretized in time. Here, we use
an explicit, the third order TVD Runge-Kutta scheme \cite{EH21}. Casting (\ref{vv}) in the form
\[
\frac{\partial u_h}{\partial t} =L_h(u_h,t),
\]
the scheme reads as
\begin{align}
 \label{tt}
&u_h^*=u_h^n+\Delta t_n L_h(u_h^n,t_n), \notag\\
&u_h^{**}=\frac{3}{4}u_h^n+\frac{1}{4}(u_h^*+\Delta t_n L_h(u_h^*,t_n+\Delta t_n)),  \\
&u_h^{n+1}=\frac{1}{3}u_h^n+\frac{2}{3}(u_h^{**}+\Delta t_n L_h(u_h^{**},t_n+\frac{1}{2}\Delta t_n )). \notag
\end{align}

The time step $\Delta t_n$ is chosen to ensure the stability of the method. For a fixed mesh,
the time step is commonly taken as
\[
\Delta t'=\hbox{CFL}\cdot \frac{\min\limits_{j}R_j^n}{\max\limits_u |F'(u_h^n)\cdot \vec{n}|},
\]
where $R_j^n$ is the radius of the inscribed circle of the $j$th element at $t_n$.
For a moving mesh, we need to consider
the effects of mesh movement and thus take the time step as
\[
\Delta t''=\hbox{CFL}\cdot \frac{\min\limits_{j}R_j^{n+1}}{\max\limits_u |(F'(u_h^n)-\dot X^n)\cdot \vec{n} |}.
\]
Finally, the time step is taken as  $\Delta t_n=\min\{\Delta t',\Delta t''\}$.

%
%
 
\subsection{Identification of troubled cells }

We now discuss the use of the TVB limiter to detect troubled cells \cite{EH10, EH13, EH15}.
Let $\textbf{x}_{m_l},\, l=1,2,3$ be the midpoints of the edges of the target element $K$,
$\textbf{x}_{b_i},\; i=1,2,3$ be the barycenters of the neighboring triangles $K_i,\; i=1,2,3$
and $\textbf{x}_{b_0}$ the barycenter of $K$. For the first edge we have
\begin{align}
\label{tvb}
\textbf{x}_{m_1}-\textbf{x}_{b_0}=\alpha_1(\textbf{x}_{b_1}-\textbf{x}_{b_0})+\alpha_2(\textbf{x}_{b_l}-\textbf{x}_{b_0}),
\quad l = 2, 3.
\end{align}
We choose $l$ such that the parameters $\alpha_1, \alpha_2$ in (\ref{tvb}) are nonnegative, which depend only on $\textbf{x}_{m_1}$ and the geometry of the elements. Then we define
\begin{align}
&\widetilde{u}(\textbf{x}_{m_1},t)\equiv u_h(\textbf{x}_{m_1},t)-\bar u_K,\notag \\
&\Delta u(\textbf{x}_{m_1},t) \equiv \alpha_1(\bar u_{K_1}-\bar u_K)+\alpha_2(\bar u_{K_l}-\bar u_K),\notag
\end {align}
where $\bar u$ is the average of $u$ (cf. (\ref{average-1})).
The value is modified by the standard minmod limiter
\[
\widetilde{u}^{(mod)}=\widetilde{m}({\widetilde{u}(\textbf{x}_{m_1},t)},\gamma\Delta u(\textbf{x}_{m_1},t)),
\]
where $\gamma >1$ is a parameter (which is taken as $\gamma=1.5$ in our computation)
and $\widetilde m$ is defined by the TVB modified minmod function as
\begin{equation}
\widetilde{m}(a_1,a_2)=
\left
\{
\begin{array}{ll}
a_1,&\quad \hbox{if}\quad |a_1|\leq \left ( 2\max\limits_j(R_j^n)\right )^2\\
m(a_1,a_2),&\quad \hbox{otherwise}
\end{array}
\right.
\end{equation}
\begin{equation}
m(a_1,a_2)=
\left
\{
\begin{array}{ll}
\hbox{sign}(a_1) \min(|a_1|,|a_2|),&\quad \hbox{if}\quad \hbox{sign}(a_1)=\hbox{sign}(a_2)\\
0,&\quad \hbox{otherwise} .
\end{array}
\right.
\end{equation}
The procedure is repeated for the other two edges. 
Finally, $K$ is marked as a troubled cell for reconstructions
if there is  $\widetilde{u}^{(mod)}\neq\widetilde u(\textbf{x}_{m_1},t)$ at least on one edge.

Note that the nonnegative parameters $\alpha_1, \alpha_2$ exist for a uniform mesh in general.
However, this is not always true for a moving mesh. To fix the problem, we set the element $K$
as a troubled cell if the nonnegative parameters do not exist on the element $K$.
We also remark that there exist other methods such as KXRCF \cite{EH42}
to identify troubled cells.

\subsection{The limiting procedure}
In this subsection we present a limiting procedure from \cite{EH17,EH18}. The key idea of the procedure is to reconstruct the entire polynomial on a troubled cell as a convex combination of the DG solution polynomial on this cell and the "modified" DG solution polynomials on its immediately neighboring cells. The modification procedure is in a least square manner \cite{EH41}. 

Assume that $K$ is a troubled cell and $K_i\; (i=1,2,3)$ are the neighboring cells of $K$.
The DG solution $u_h$ on $K,\; K_1,\; K_2$, and $K_3$ are denoted by
$p_0(\textbf{x}),\; p_1(\textbf{x}),\; p_2(\textbf{x})$, and $p_3(\textbf{x})$, respectively,
and their averages over the corresponding cells by $\bar p_0$, $\bar p_1$, $\bar p_2$, and $\bar p_3$.

\subsubsection{The scalar case}

The limiting procedure for two-dimensional scalar conservation laws is given in the following.

{\bf Step 1.} We first modify the polynomial $p_1(\textbf{x})$ into a polynomial $\widetilde{p}_1(\textbf{x})$ on $K_1$ in the least square sense, i.e.,  $\widetilde{p}_1(\textbf{x})$ is defined as the solution of the minimization problem 
\begin{align}
\label{min1}
\min\limits_{\forall \varphi(\textbf{x})\in P^k(K_1)}
\int_{K_1}(\varphi(\textbf{x})-p_1(\textbf{x}))^2 dxdy
+\sum\limits_{\ell\in {L_1}}(\int_{K_{\ell}}\varphi(\textbf{x})-p_{\ell}(\textbf{x})dxdy)^2,
\end {align}
subject to
\[
\frac{1}{|K|}\int_{K} \varphi(\textbf{x}) dxdy=\bar p_0,
\]
where 
$
L_1=\{2,3\}\cap \{ \ell:|\bar p_{\ell}-\bar p_0|<\hbox{max}(|\bar p_2-\bar p_0|,|\bar p_3-\bar p_0|) \}.
$
Polynomials $\widetilde{p}_2(\textbf{x})$ and $\widetilde{p}_3(\textbf{x})$ can be obtained similarly.

{\bf Step 2.} The linear weights $\gamma_0, \gamma_1, \gamma_2$, and $\gamma_3$ are chosen.
In principle, they can be taken as any set of positive numbers with a unitary sum since there is no
constraint on the linear weights in order to maintain the accuracy of the method. In our computation,
we take $\gamma_0=0.997, \gamma_1=0.001, \gamma_2=0.001$, and $\gamma_3=0.001$
attempting to keep the reconstructed polynomial as close as possible to the original one.

{\bf Step 3.} The smoothness indicators, denoted by $\beta_l,\; l=0,1,2,3$ are computed (cf.  \cite{EH18}) as 
\[
\beta_l=\sum\limits^k_{|s|=1}\int_{K}|K|^{|s|-1}\left (\frac{1}{|s|!}\frac{\partial^{|s|}\widetilde{p}_l
(\textbf{x})}{\partial x^{s_1}\partial y^{s_2}}\right )^2dxdy,\quad l=0,1,2,3
\]
where $s=(s_1,s_2)$ and $|s|=s_1+s_2.$
They measure how smooth the function $\widetilde{p}_l(\textbf{x})$ is in the target cell $K$.

{\bf Step 4.} The nonlinear weights based on smoothness indicators and linear weights are defined as
\[
\omega_l=\frac{\Bar \omega_j}{\sum\limits_m\Bar \omega_m},
\qquad \Bar \omega_m=\frac{\gamma_m}{(\lambda+\beta_m)^2},\quad l=0,1,2,3
\]
where $\lambda$ is a small number to avoid the denominator to become zero, in this paper, we take $\lambda=10^{-6}$.

{\bf Step 5.} Finally, we obtain the reconstructed polynomial as 
\[
p_0^{new}=w_0{p}_0(\textbf{x})+w_1\widetilde{p}_1(\textbf{x})+w_2\widetilde{p}_2(\textbf{x})
+w_3\widetilde{p}_3(\textbf{x})
\]
and define $u_h^{new}|_{K}=p_0^{new}(\textbf{x})|_{K}$.

\subsubsection{The system case}
We now describe the limiting procedure for systems which uses the local characteristic field decomposition
for better nonoscillatory properties. To be specific,  we consider the Euler system in two dimensions as
\begin{equation}
\label{2d}
u_t+f_1(u)_x+f_2(u)_y\equiv \frac{\partial }{\partial t}
\left(
  \begin{array}{c}
    \rho  \\
     \rho\mu   \\
      \rho\nu   \\
      E \\
   \end{array}
 \right)
 +\frac{\partial }{\partial x}
 \left(
  \begin{array}{c}
    \rho\mu  \\
     \rho\mu^2+P   \\
      \rho\mu\nu   \\
      \mu(E+P) \\
   \end{array}
 \right)
 +\frac{\partial }{\partial y}
 \left(
  \begin{array}{c}
    \rho\nu  \\
     \rho\mu\nu   \\
      \rho\nu^2+P   \\
      \nu(E+P) \\
   \end{array}
 \right)
 =0,
\end{equation}
with $u(\textbf{x},0)=u_0(\textbf{x}),$ where $\rho$ is the density, $\mu$ and $\nu$ are the velocity components
in the $x$- and $y$-direction, respectively, $E$ is the energy density, and $P$ is the pressure. The equation of state is $E=\frac{P}{\gamma-1}+\frac{1}{2}\rho(\mu^2+\nu^2)$ with $\gamma=1.4$.

Let $F^{\prime}(\bar u_K)\cdot \vec{n}_i =(f'_1(\bar u_K),f'_2(\bar u_K)) \cdot \vec{n}_i$ be the Jacobian matrices and $\vec{n}_i=(n_{ix},n_{iy}),\; i=1,2,3$ be the outward unit normals to the edges of the element $K$. Then the left and right eigenvectors of the Jacobian matrices are given by
\begin{equation}
L_i(\bar u_K)=
 \left(
  \begin{array}{cccc}
    \frac{B_2+(\mu n_{ix}+\nu n_{iy})/c}{2}  & -\frac{B_1\mu+n_{ix}/c}{2} & -\frac{B_1\nu+n_{iy}/c}{2} & \frac{B_1}{2} \\
     n_{iy}\mu-n_{ix}\nu & -n_{iy} & n_{ix} & 0  \\
      1-B_2 & B_1\mu & B_1\nu & -B_1  \\
    \frac{B_2-(\mu n_{ix}+\nu n_{iy})/c}{2}  & -\frac{B_1\mu-n_{ix}/c}{2} & -\frac{B_1\nu-n_{iy}/c}{2} & \frac{B_1}{2} \\
   \end{array}
 \right) ,
\end{equation}
\begin{equation}
R_i(\bar u_K)=
 \left(
  \begin{array}{cccc}
    1 & 0 & 1 & 1  \\
     \mu-cn_{ix} & -n_{iy} & \mu & \mu+cn_{ix}  \\
      \nu-cn_{iy} & n_{ix}& \nu & \nu+cn_{iy}  \\
    H-c(\mu n_{ix}+\nu n_{iy}) & -n_{iy}\mu+n_{ix}\nu & \frac{\mu^2+\nu^2}{2} & H+c(\mu n_{ix}+\nu n_{iy}) \\
   \end{array}
 \right) ,
\end{equation}
where $c=\sqrt{\frac{\gamma P}{\rho}}$, $B_1=\frac{\gamma-1}{c^2}$, $B_2=\frac{B1(\mu^2+\nu^2)}{2}$, and $H=\frac{E+P}{\rho}$.
The limiting procedure for the Euler system is as follows.

{\bf Step 1.} For each direction $\vec{n}_i,\; i=1,2,3,$ project the polynomials $p_0$, $p_1$, $p_2$, and $p_3$ into the characteristic directions,  i.e., $ \widetilde p_{i_l}=L_i\cdot p_l,\; l=0,1,2,3.$ Then the limiting procedure in the previous
subsection for the scalar case is applied and the new modified polynomial is denoted by $\widetilde{p}_{i,0}^{new}$.
Next,  $\widetilde{p}_{i,0}^{new}$ is projected back into the physical space, that is,
$p_{i,0}^{new}=R_i\cdot\widetilde{p}_{i,0}^{new},\; i=1,2,3.$

{\bf Step 2.} The final new polynomial on the troubled cell $K$ is computed as
\[
p_0^{new}=\frac{\sum\limits_{i=1}^{3}p_{i,0}^{new}|K_i|}{\sum\limits_{i=1}^{3}|K_i|}
\]
and we define $u_h^{new}|_{K}=p_0^{new}(\textbf{x})|_{K}$.

\section{The MMPDE moving mesh strategy}
\label{sec:mmpde}

In this section we discuss the generation of $\mathscr{T}_h^{n+1}$ using the MMPDE moving mesh method.
To this end, we assume that the mesh $\mathscr{T}_h^{n}$ at $t = t_n$ and the numerical approximation $u_h^n$
of a physical variable $u$
are given. We also assume that a reference computational mesh $\Hat {\mathscr{T}}_c = \{\Hat {\pmb{\xi}}_j\}^{N_v}_{j=1}$, a deformation of the physical mesh, has been chosen. This mesh is fixed for the whole computation
and can be taken as the initial physical mesh. For the purpose of mesh generation, we need to use another mesh,
called the computational mesh ${\mathscr{T}}_c =  \{{\pmb{\xi}}_j\}^{N_v}_{j=1}$, which is also a deformation
of the physical mesh and will be used as an intermediate variable.

The MMPDE method views any nonuniform mesh as a uniform one in some metric specified by a metric tensor
$\mathbb{M} = \mathbb{M}(\textbf{x})$. The metric tensor $\mathbb{M}(\textbf{x})$ is a symmetric and positive definite matrix for each $\textbf{x}$ and uniformly positive definite on $\Omega$. It provides the information
needed to control the size, shape and orientation of the mesh elements throughout the domain. 
In our computation we use
\begin{align}
\label{mer}
\mathbb{M}=\hbox{det}(\mathbb{I}+|H(u_h^n)|)^{-\frac{1}{d+4}}(\mathbb{I}+|H(u^n_h)|),
\end {align}
where $d$ is the dimension of the domain ($d = 1$ for one dimension and $d=2$ for two dimensions),
$\mathbb{I}$ is the $d\times d$ identity matrix, $H(u_h^n)$ is a recovered Hessian from the numerical solution $u^n_h$, $|H(u^n_h)|=Q\hbox{diag}(|\lambda_1|,\cdots,|\lambda_d|)Q^T$ with $Q\hbox{diag}(\lambda_1,\cdots,\lambda_d)Q^T$ being the eigen-decomposition of $H(u^n_h)$, and $\hbox{det}(\cdot)$ is the determinant of a matrix. 
The metric tensor (\ref{mer}) is known to be optimal for the $L^2$ norm of linear interpolation error \cite{EH30}. In our computation the Hessian is recovered using the least square fitting \cite{EH31}.
It is common practice in moving mesh computation to smooth the metric tensor/monitor function
for smoother meshes. To this end, we apply a low-pass filter \cite{EH22, EH08} to the smoothing of the metric tensor
several sweeps every time it is computed.

For scalar equations, we use the solution to the equation as $u$ in computing $\mathbb{M}$. 
For the Euler system, motivated by the choice in \cite{EH32}, we take $u$ to be the quantity
$S=0.5 \sqrt{1+\beta (\frac{\rho}{\max(\rho)})^2}+ 0.5 \sqrt{1+\beta(\frac{E}{\max(E)})^2}$.
Its nodal value at $\textbf{x}_j$ ($j=1,2,\cdots,N_v$) is computed as
\begin{align}
\label{ent}
&S_j=0.5 \sqrt{1+\beta (\frac{\rho_j}{\max\limits_{1\leqslant m\leqslant N_v}(\rho_m)})^2}+0.5\sqrt{1+\beta(\frac{E_j}{\max\limits_{1\leqslant m\leqslant N_v}(E_m)})^2}, \\
&{\rho}_j=\frac{\sum\limits_{K\in\omega_j}|K|{\rho}_K}{\sum\limits_{K\in\omega_j}|K|},\qquad 
 {E}_j=\frac{\sum\limits_{K\in\omega_j}|K|{E}_K}{\sum\limits_{K\in\omega_j}|K|},\notag
\end{align}
where $\omega_j$ is the element patch associated with $\textbf{x}_j$ and $\beta$ is a positive parameter. 
The choice of $\beta$ is given in Section~\ref{sec:numerics}.

We now describe the MMPDE moving mesh strategy. A mesh $\mathscr{T}_h$ is uniform in the metric $\mathbb{M}$
with reference to a computational mesh $\mathscr{T}_c$ will be referred to as an $\mathbb{M}$-uniform mesh
with respect to $\mathscr{T}_c$. It is known \cite{EH38} that such a mesh satisfies
the equidistribution and alignment conditions
\begin{align}
|K|\sqrt{\hbox{det}(\mathbb{M}_K)}=\frac{\sigma_h|K_c|}{|\Omega_c|},\quad \forall K\in \mathscr{T}_h \label{ec}\\
\frac{	1}{d}\text{tr}((F'_K)^{-1}\mathbb{M}_K^{-1}(F'_K)^{-T})=\hbox{det}((F'_K)^{-1}\mathbb{M}_K^{-1}(F'_K)^{-T})^{\frac{1}{d}},\quad \forall K\in\mathscr{T}_h \label{al}
\end {align}
where $K_c$ is the element in $ \mathscr{T}_c$ corresponding
to $K$, $F_K$ is the affine mapping from $K_c$ to $K$ and $F_K'$ is its Jacobian matrix,
$\mathbb{M}_K$ is the average of $\mathbb{M}$ over $K$, $\hbox{tr}(\cdot)$ is the trace of a matrix, and
\[
|\Omega_c|=\sum\limits_{K_c\in\mathscr{T}_c}|K_c|,\quad
\sigma_h =\sum\limits_{K\in\mathscr{T}_h}|K|\hbox{det}(\mathbb{M}_K)^{\frac{1}{2}} .
\]
The equidistribution condition (\ref{ec}) determines the size of the element $K$ through
the metric tensor $\mathbb{M}$. The volume $|K|$ is smaller in regions where
$\text{det}(\mathbb{M}_K)^{\frac{1}{2}}$ is larger. On the other hand, the alignment condition (\ref{al})
determines the shape and orientation of $K$ through $\mathbb{M}_K$ and $K_c$.

The objective of the MMPDE moving mesh method is to generate a mesh satisfying the above two conditions
as closely as possible. This is done by minimizing the energy function
\begin{align}
\label{fl}
I_h(\mathscr{T}_h,\mathscr{T}_c) &= \sum\limits_{K\in\mathscr{T}_h}|K|\sqrt{\hbox{det}(\mathbb{M}_K)}(\hbox{tr}((F'_K)^{-1}\mathbb{M}^{-1}_K(F'_K)^{-T}))^{\frac{3 d}{4}} \notag\\
 &\qquad + d^{\frac{3 d}{4}}\sum\limits_{K\in\mathscr{T}_h} |K|\sqrt{\hbox{det}(\mathbb{M}_K)}\left (\frac{|K_c|}{|K|\sqrt{\hbox{det}(\mathbb{M}_K)}}\right )^{\frac{3}{2}},
\end{align}
which is a Riemann sum of a continuous functional developed in \cite{EH27} based on equidistribution and alignment for variational mesh adaptation. 
Notice that $I_h(\mathscr{T}_h,\mathscr{T}_c)$ is a function of the vertices $\{ \pmb{\xi}_j\}$ of the computational mesh $\mathscr{T}_c$ and the vertices $\{\textbf{x}_j\}$ of the physical mesh $\mathscr{T}_h$.
Here we use the $\xi$-formulation where we take $\mathscr{T}_h = \mathscr{T}_h^n$ and minimize
$I_h(\mathscr{T}_h^n,\mathscr{T}_c)$ by solving its gradient system with respect to $\{ \pmb{\xi}_j\}$.
Thus, the mesh equation reads as
\begin{align}
\label{MM}
\frac{d\pmb{\xi}_j}{dt}=-\frac{P_j}{\tau}\left (\frac{\partial I_h}{\partial \pmb{\xi}_j}\right )^T,\quad j=1,2,\cdots,N_v
\end {align}
where $\frac{\partial I_h}{\partial \pmb{\xi}_j}$ is considered as a row vector, $\tau>0$ is a parameter used to adjust the time scale of the mesh movement to respond the changes in $\mathbb{M}$, and $P_j$ is a positive function used to make the MMPDE to have desired invariant properties. Here, we take $P_j=\hbox{det}(\mathbb{M}(\textbf{x}_j))^{\frac{p-1}{2}}$ so that (\ref{MM}) is invariant under scaling transformations of $\mathbb{M}$.
Using the notion of scalar-by-matrix differentiation, we can find the analytical formulations
of the derivatives in (\ref{MM}) \cite{EH23} and rewrite the mesh equation as
\begin{align}
\label{xim}
\frac{d\pmb{\xi}_j}{dt}=\frac{P_j}{\tau}\sum\limits_{K\in\omega_j}|K|\pmb{v}^K_{j_K},\quad j=1,2,\cdots,N_v
\end{align}
where $j_K$ is the local index of $\textbf{x}_j$ in $K$ and $\pmb{v}^K_{j_K}$ is the local velocity for $\textbf{x}_j$
contributed by $K$. The local velocities contributed by $K$ to its vertices are given by
\begin{equation}
\left[
  \begin{array}{c}
    ( \pmb{v}_1^K  )^T  \\
      \vdots   \\
     ( \pmb{v}_d^K  )^T     \\
   \end{array}
 \right]
 =-E_K^{-1}\frac{\partial G}{\partial \mathbb{J}}-\frac{\partial G}{\partial \hbox{det}(\mathbb{J})}\frac{\partial \hbox{det}(\Hat E_K)}{\partial \hbox{det}(E_K)}\Hat E_K^{-1},\quad \pmb{v}^K_0=-\sum\limits_{i=1}^d\pmb{v}^K_d,
\end {equation}
where the vertices of $K$ and $K_c$ are denoted by $\textbf{x}_j^K,\; j = 0, 1, ..., d$ and
$\pmb{\xi}_j^K,\; j = 0, 1, ..., d$, respectively,
$E_K=[\textbf{x}_1^K-\textbf{x}_0^K,\cdots,\textbf{x}_d^K-\textbf{x}_0^K]$ and $E_{K_c}=[\pmb{\xi}_1^K-\pmb{\xi}_0^K,\cdots,\pmb{\xi}_d^K-\pmb{\xi}_0^K]$ are the edge matrices of $K$ and $K_c$,
the function $G=G(\mathbb{J},\hbox{det}(\mathbb{J}), \mathbb{M}_K)$ with $\mathbb{J}=(F'_K)^{-1}=E_{K_c}E_K^{-1}$ is associated with the energy function (\ref{fl}), and its definition and derivatives are given by
\begin{align}
&G(\mathbb{J},\hbox{det}(\mathbb{J}),\mathbb{M})=\sqrt{\hbox{det}(\mathbb{M})}(\hbox{tr}(\mathbb{J}\mathbb{M}^{-1}\mathbb{J}^T))^{\frac{3 d}{4}}+d^{\frac{3 d}{4}}\sqrt{\hbox{det}(\mathbb{M})}\left (\frac{{\hbox{det}(\mathbb{J})}}{\sqrt{\hbox{det}(\mathbb{M})}}\right )^{\frac{3}{2}},\notag\\
&\frac{\partial G}{\partial \mathbb{J}}=\frac{3d}{2}\sqrt{\hbox{det}(\mathbb{M})}(\hbox{tr}(\mathbb{J}\mathbb{M}^{-1}\mathbb{J}^T))^{\frac{3 d}{4}-1}\mathbb{M}^{-1}\mathbb{J}^T,\notag\\
&\frac{\partial G}{\partial \hbox{det}(\mathbb{J})}=\frac{3}{2} d^{\frac{3d}{4}}\hbox{det}(\mathbb{M})^{-\frac{1}{4}}\hbox{det}(\mathbb{J})^{\frac{1}{2}}.  \notag
\end {align}

In practical computation, we first compute the edge matrices and the local velocities for all elements and then use (\ref{xim}) to obtain the nodal mesh velocities. The mesh equation is modified for the boundary mesh points. For fixed points, the mesh velocity can be set to be zero. For those on a boundary edge, the mesh velocities should be modified to ensure they stay on the boundary. 

The mesh equation (\ref{xim}) (with the proper  modifications for boundary vertices) can be integrated from $t^n$ to $t^{n+1}$, starting  with the reference computational mesh $\Hat {\mathscr{T}}_c$ as an initial mesh. Then the new computational mesh $\mathscr{T}_c^{n+1}$ is obtained and forms a correspondence $\psi_h$
with the physical mesh $\mathscr{T}_h^{n}$ with the property $\textbf{x}_j^n=\psi_h(\pmb{\xi}_j^{n+1}),\; j=1,2,\cdots,N_v$. Finally, the new physical mesh $\mathscr{T}_h^{n+1}$ is defined as $\textbf{x}_j^{n+1}=\psi_h(\Hat {\pmb{\xi}}_j),\; j=1,2,\cdots,N_v$, which can readily be computed using 
linear interpolation.

The equation (\ref{xim}) is called the ${\xi}$-formulation of the MMPDE moving mesh method.
Alternatively, we can use the $x$-formulation where we take $\mathscr{T}_c = \hat{\mathscr{T}}_c$ and minimize $I_h$
by integrating its gradient system with respect to $\{ \textbf{x}_j\}$; see \cite{EH38}.
Although its implementation is more complex and costly than the $\xi$-formulation, the $x$-formulation
has the advantage that it can be shown analytically \cite{EH37} that the moving mesh
governed by the $x$-formulation will stay free of tangling and cross-over for both convex or concave domains.
Such a theoretical result is not available for the $\xi$-formulation although the numerical experiment shows that
it also produces nonsingular moving meshes.

\section{Numerical examples}
\label{sec:numerics}
\setcounter{equation}{0}
\setcounter{figure}{0}
\setcounter{table}{0}

In this section we present numerical results obtained with the moving mesh DG method described
in the previous sections for a selection of one- and two-dimensional examples.
Recall that the method has been described in two dimensions. Its implementation in one dimension
is similar. The CFL number in time step selection is set to be 0.3 for $P^1$ elements and 0.15 for
$P^2$ elements. The parameter $\tau$ in \eqref{MM} is taken as $0.1$ for accuracy test problems and
$10^{-3}$ and $10^{-4}$ for one- and two-dimensional { systems} with discontinuities, respectively.
The parameter $\beta$ in \eqref{ent} is taken as 10 for one-dimensional examples and 1 for two-dimensional problems,
unless otherwise stated.  Moving and uniform meshes will be denoted by "MM'' and "UM'', respectively. 
Unless otherwise stated, three sweeps of a low-pass filter \cite{EH22, EH08} are applied to the smoothing of
the metric tensor every time it is computed.
In addition, for examples having an exact solution, the error of the computed solution is measured
in the global norm, i.e.,
\[
\|e_h\|_{L^q}=\left (\int_{0}^T\int_{\Omega}|e_h(\textbf{x},t)|^qd\textbf{x}dt\right )^{\frac{1}{q}}, \quad q=1,2,\infty.
\]

\subsection{One-dimensional examples}

\begin{exam}{\em
\label{exam4.1}
We first consider Burgers' equation
\[
u_t+\left (\frac{u^2}{2}\right )_x=0,\quad x\in (0,2)
\]
subject to the initial condition $u(x,0)=0.5+\hbox{sin}(\pi x)$ and a periodic boundary condition. We compute the solution up to $T = \frac{0.5}{\pi}$ when the solution is still smooth and the exact solution can be computed using Newton's iteration. The error is listed in Table \ref{ex4.1}, which shows the convergence of the second order for $P^1$ elements
and the third order for $P^2$ elements for the moving mesh DG method.

To test the convergence of the method for discontinuous solutions, we compute Burgers' equation
up to $T = \frac{1.5}{\pi}$ when the shock appears. The error of $L^1$ is listed in Table \ref{ex4nc30}. The results indicate that the computed solution is convergent although the convergence order decreases to one for both $P^1$ and $P^2$ elements.

To see how the smoothness of the mesh affects the accuracy of the method, we list the $L^1$ error
in Tables~\ref{ex4.1nc30} and \ref{ex4nc30} for the solutions computed with different numbers
of sweeps of the low-pass filter applied to the metric tensor. One can see that the results are almost
the same. In Table~\ref{ex4nc30} where the solution is discontinuous, the error is slightly worse
for the case with more sweeps. This is because more sweeps lead to smoother meshes with
less concentration near the shock. Overall, the results show that the accuracy of the method
is not sensitive to the smoothness of the mesh, which is in contrast with the situation
for the moving mesh finite difference WENO method \cite{EH08}.

We also compute Burgers' equation with a discontinuous initial condition 
\begin{equation}
u(x,0)=
\left
\{
\begin{array}{ll}
1,      \quad&\text{for}\quad -1\leqslant x \leqslant 0\\
0, \quad&\text{for}\quad    0<x\leqslant 1.\notag
\end{array}
\right.
\end{equation}
The error for both $P^1$ and $P^2$ elements at $T=1$ is listed in Table \ref{ex4.1.1}.
Once again, the results show that the computed solution is convergent at a rate of the first order
for both $P^1$ and $P^2$ elements.

}\end{exam}

\begin{table}
\caption{Example~\ref{exam4.1}: Solution error with periodic boundary conditions and $T=\frac{0.5}{\pi}$.}
\renewcommand{\multirowsetup}{\centering}
\begin{center}
\begin{tabular}{|c|c|c|c|c|c|c|c|c|c|c|c|c|}
\hline
$k$ & $N$  &20      & 40 & 80 & 160 &320 &640\\
\hline
\multirow{6}{1cm}{1}
 & $L^1$ &1.400e-3      & 3.609e-4 & 7.379e-5 & 1.707e-5  & 4.341e-6 & 1.102e-6\\
 &  Order      & \quad  & 1.96       & 2.29      &  2.11        & 1.98 &1.98 \\
 &$L^2$   & 4.196e-3      & 1.103e-3 & 2.100e-4 &  4.762e-5  & 1.216e-5 & 3.100e-6\\
 &Order   & \quad   & 1.93     & 2.39      &   2.14      & 1.97 & 1.97 \\
 &$L_{\infty}$ & 6.423e-3 & 1.754e-3 & 3.265e-4  &  5.917e-5 & 1.441e-5 & 3.638e-6\\
 &Order    & \quad  & 1.87      & 2.43     &   2.46      & 2.04 & 1.99   \\
 \hline
 \multirow{6}{1cm}{2}

 & $L^1$  & 4.678e-5      &6.596e-6 & 8.693e-7 &  1.094e-7 & 1.315e-8 & 1.533e-9 \\
 &  Order      & \quad & 2.83    & 2.92     & 2.99     & 3.06 & 3.10  \\
 &$L^2$  & 1.702e-4       & 3.031e-5 & 4.962e-6 &  7.312e-7 & 9.402e-8 & 1.056e-8\\
 &Order       & \quad &2.49     & 2.61       & 2.76    & 2.96 & 3.15 \\
 &$L_{\infty}$ & 3.419e-4 & 6.959e-5 & 1.350e-5  & 2.323e-6 & 3.542e-7 & 4.707e-8 \\
 &Order    & \quad       &2.30    & 2.37      & 2.54   & 2.71 &2.91   \\
 \hline

\end{tabular}
\end{center}
\label{ex4.1}
\end{table}


\begin{table}
\caption{Example~\ref{exam4.1}: Solution error in $L^1$ norm with periodic boundary conditions and $T=\frac{0.5}{\pi}$. Various numbers of sweeps of a low-pass filter have been applied to the smoothing of the metric tensor.}
\renewcommand{\multirowsetup}{\centering}
\begin{center}
\begin{tabular}{|c|c|c|c|c|c|c|c|c|c|c|c|c|}
\hline
$k$ & Sweeps$\backslash$ $N$  &20      & 40 & 80 & 160 &320 &640\\
\hline
\multirow{6}{1cm}{1}
  & 3 &1.400e-3      & 3.609e-4 & 7.379e-5 & 1.707e-5  & 4.341e-6 & 1.102e-6\\
 &  Order      & \quad  & 1.96       & 2.29      &  2.11        & 1.98 &1.98 \\
 & 30 &1.657e-3      & 3.797e-4 & 7.425e-5 & 1.716e-5  & 4.353e-6 & 1.103e-6\\
 &  Order      & \quad  & 2.13       & 2.35      &  2.11        & 1.98 &1.98 \\
& 100 &1.737e-3      & 3.993e-4 & 7.483e-5 & 1.726e-5  & 4.368e-6 & 1.104e-6\\
 &  Order      & \quad  & 2.12       & 2.42      &  2.12        & 1.98 &1.98 \\
 \hline
 \multirow{6}{1cm}{2}

& 3  & 4.678e-5      &6.596e-6 & 8.693e-7 &  1.094e-7 & 1.315e-8 & 1.533e-9 \\
 &  Order      & \quad & 2.83    & 2.92     & 2.99     & 3.06 & 3.10  \\
 & 30  & 4.228e-5      &5.494e-6 & 7.290e-7 &  9.522e-8 & 1.200e-8 & 1.462e-9 \\
 &  Order      & \quad & 2.94    & 2.91     & 2.94     & 2.99 & 3.03  \\
 & 100  & 4.343e-5      &5.336e-6 & 6.860e-7 &  8.920e-8 & 1.142e-8 & 1.421e-9 \\
 &  Order      & \quad & 3.02    & 2.96     & 2.94     & 2.97 & 3.01  \\
 \hline

\end{tabular}
\end{center}
\label{ex4.1nc30}
\end{table}

\begin{table}
\caption{Example~\ref{exam4.1}: Solution error in $L^1$ norm with periodic boundary conditions and $T=\frac{1.5}{\pi}$. Various numbers of sweeps of a low-pass filter have been applied to the smoothing of the metric tensor.
}
\renewcommand{\multirowsetup}{\centering}
\begin{center}
\begin{tabular}{|c|c|c|c|c|c|c|c|c|c|c|c|}
\hline
$k$ & Sweeps$\backslash$ $N$     & 80 & 160 &320 &640  & 1280 & 2560 & 5120\\
\hline
  & 3       & 7.301e-4 & 2.211e-4  & 7.985e-5 & 3.231e-5 &1.388e-5 & 6.258e-6 & 2.897e-6\\
 &  Order      & \quad      &  1.72        & 1.47 &1.31 & 1.22 & 1.15 & 1.11\\
\multirow{2}{1cm}{1}
  & 30       & 1.405e-3 & 4.787e-4  & 1.725e-4 & 7.038e-5 &3.230e-5 & 1.567e-5 & 7.727e-6\\
 &  Order      & \quad      &  1.55     & 1.47 &1.29 & 1.12 & 1.04 & 1.02\\
   & 100       & 1.929e-3 & 7.550e-4  & 2.945e-4 & 1.212e-4 &5.649e-5 & 2.862e-5 & 1.479e-5\\
 &  Order      & \quad      &  1.35     & 1.36 &1.28 & 1.10 & 0.98 & 0.95\\
 \hline
   & 3       & 3.718e-4 & 1.305e-4 & 5.374e-5 & 2.468e-5 & 1.189e-5 & 5.822e-6 & 2.856e-6\\
 &  Order      & \quad     & 1.51    & 1.28 & 1.12 & 1.05 & 1.03 & 1.03 \\
 \multirow{2}{1cm}{2}
  & 30       & 8.943e-4 & 3.317e-4 & 1.266e-4 & 5.521e-5 & 2.691e-5 & 1.374e-5 & 7.065e-6\\
 &  Order      & \quad    & 1.43    & 1.39 & 1.20 & 1.04 & 0.97 & 0.96 \\
  & 100       & 1.283e-3 & 5.466e-4 & 2.226e-4 & 9.583e-5 & 4.634e-5 & 2.404e-5 & 1.261e-5\\
 &  Order      & \quad     & 1.23    & 1.30 & 1.22 & 1.05 & 0.95 & 0.93 \\
 \hline
\end{tabular}
\end{center}
\label{ex4nc30}
\end{table}

\begin{table}
\caption{Example~\ref{exam4.1}: Solution error at $T=1$.}
\renewcommand{\multirowsetup}{\centering}
\begin{center}
\begin{tabular}{|c|c|c|c|c|c|c|c|c|c|c|c|c|}
\hline
$k$ & $N$  &20      & 40 & 80 & 160 &320 &640\\
\hline
\multirow{2}{1cm}{1}
  & $L^1$ & 7.742e-3      & 2.179e-3 & 8.449e-4 & 3.414e-4  & 1.529e-4 & 6.974e-5\\
 &  Order      & \quad     & 1.83     &  1.37  & 1.31&1.16  &1.13 \\
 \hline
 \multirow{2}{1cm}{2}
  & $L^1$  & 6.976e-3      &2.046e-3 & 8.446e-4 & 3.556e-4 & 1.761e-4 & 8.018e-5 \\
 &  Order      & \quad & 1.77    & 1.28     & 1.25    & 1.01 & 1.14  \\
 \hline
\end{tabular}
\end{center}
\label{ex4.1.1}
\end{table}


\begin{exam}{\em
\label{exam4.2}
To see the accuracy of the method for system problems, we compute the Euler equations,
\begin{equation}
\label{Euler}
\begin{pmatrix}\rho \\ \rho u \\ E \end{pmatrix}_t+\begin{pmatrix} \rho u \\ \rho u^2+P \\ u(E+P) \end{pmatrix}_x=0,
\end{equation}
where
$\rho$ is the density, $u$ is the velocity, $E$ is the energy density, and $P$ is the pressure. The equation of state is $E=\frac{P}{\gamma-1}+\frac{1}{2}\rho u^2$ with $\gamma=1.4$. The initial condition is
\[
\rho(x,0)=1+0.2\text{sin}(\pi x), \quad u(x,0)=1,\quad P(x,0)=1,
\]
and a periodic boundary condition is used. The exact solution for this problem is
$$\rho(x,t)=1+0.2\hbox{sin}(\pi(x-t)),\quad u(x,t)=1,\quad P(x,t)=1.$$
The final time is $T=1.0$. The parameter $\beta$ in \eqref{ent} is set to be 100. The error in computed density
is listed in Table~\ref{ex4.2}. From the table one can see that the $(k+1)^{\text{th}}$ order of accuracy of the scheme
is achieved for this nonlinear system.

}\end{exam}

\begin{table}
\caption{Example~\ref{exam4.2}: Error in computed density for periodic boundary conditions and $T=1.0$, $\beta=100$.}
\renewcommand{\multirowsetup}{\centering}
\begin{center}
\begin{tabular}{|c|c|c|c|c|c|c|c|c|c|c|c|c|}
\hline
$k$ & $N$       &10 & 20 & 40 & 80 & 160 &320 \\
\hline
\multirow{6}{1cm}{1}
 & $L^1$       & 5.576e-3 & 1.350e-3 & 3.310e-4 &  8.226e-5 & 2.054e-5  & 5.145e-6 \\
 &  Order      & \quad    & 2.05     & 2.03     & 2.01      &  2.00        & 2.00 \\
 &$L^2$       & 5.187e-3 & 1.265e-3 & 3.149e-4 & 8.041e-5 &  2.062e-5  & 5.246e-6 \\
 &Order       & \quad    & 2.04     & 2.01     & 1.97      &   1.96      & 1.97   \\
 &$L_{\infty}$&1.298e-2  &3.167e-3  & 7.791e-4 & 2.099e-4  &  5.576e-5 & 1.431e-5 \\
 &Order    & \quad       & 2.03     & 2.02     & 1.89    &   1.91      & 1.96    \\
 \hline
 \multirow{6}{1cm}{2}

 & $L^1$       & 3.201e-4 & 4.424e-5 &5.784e-6 & 6.934e-7 &  7.947e-8 & 9.223e-9 \\
 &  Order      & \quad    & 2.86     & 2.94     & 3.06     & 3.13     & 3.11  \\
 &$L^2$        & 3.513e-4 & 5.249e-5 & 7.216e-6 & 8.745e-7 &  9.626e-8 & 1.055e-8\\
 &Order       & \quad    & 2.74     & 2.86     & 3.04       & 3.18    & 3.19  \\
 &$L_{\infty}$&1.082e-3  &1.875e-4  & 2.926e-5 & 3.904e-6  & 4.509e-7 & 4.682e-8 \\
 &Order    & \quad       & 2.53    & 2.68    & 2.91      & 3.11   & 3.27   \\
 \hline

\end{tabular}
\end{center}
\label{ex4.2}
\end{table}

\begin{figure}[hbtp]
 \begin{center}
 \mbox{\subfigure[MM: $N=100$, UM: $N=100$]
 {\includegraphics[width=8cm]{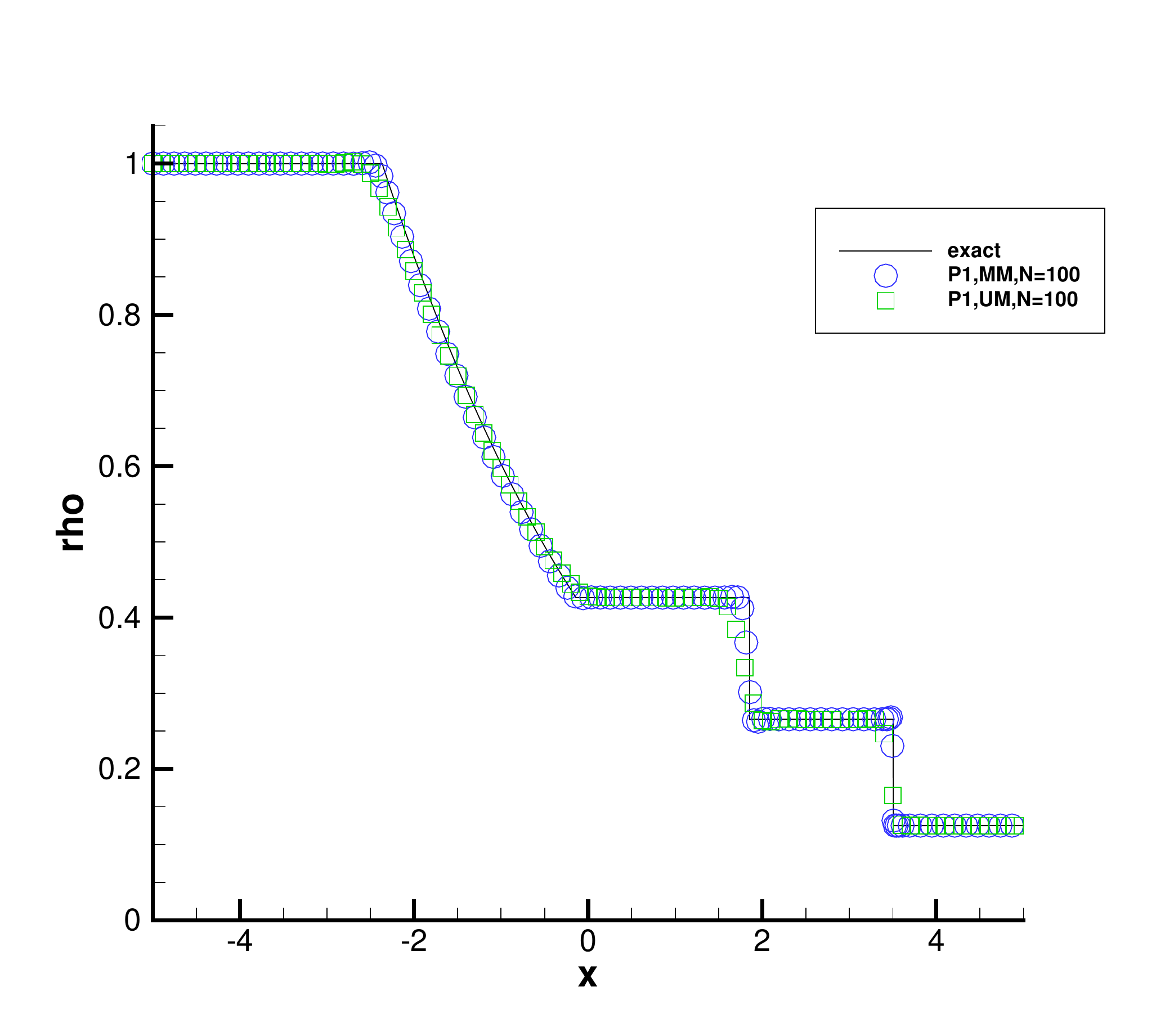}}\quad
   \subfigure[close view of (a) near shock]
   {\includegraphics[width=8cm]{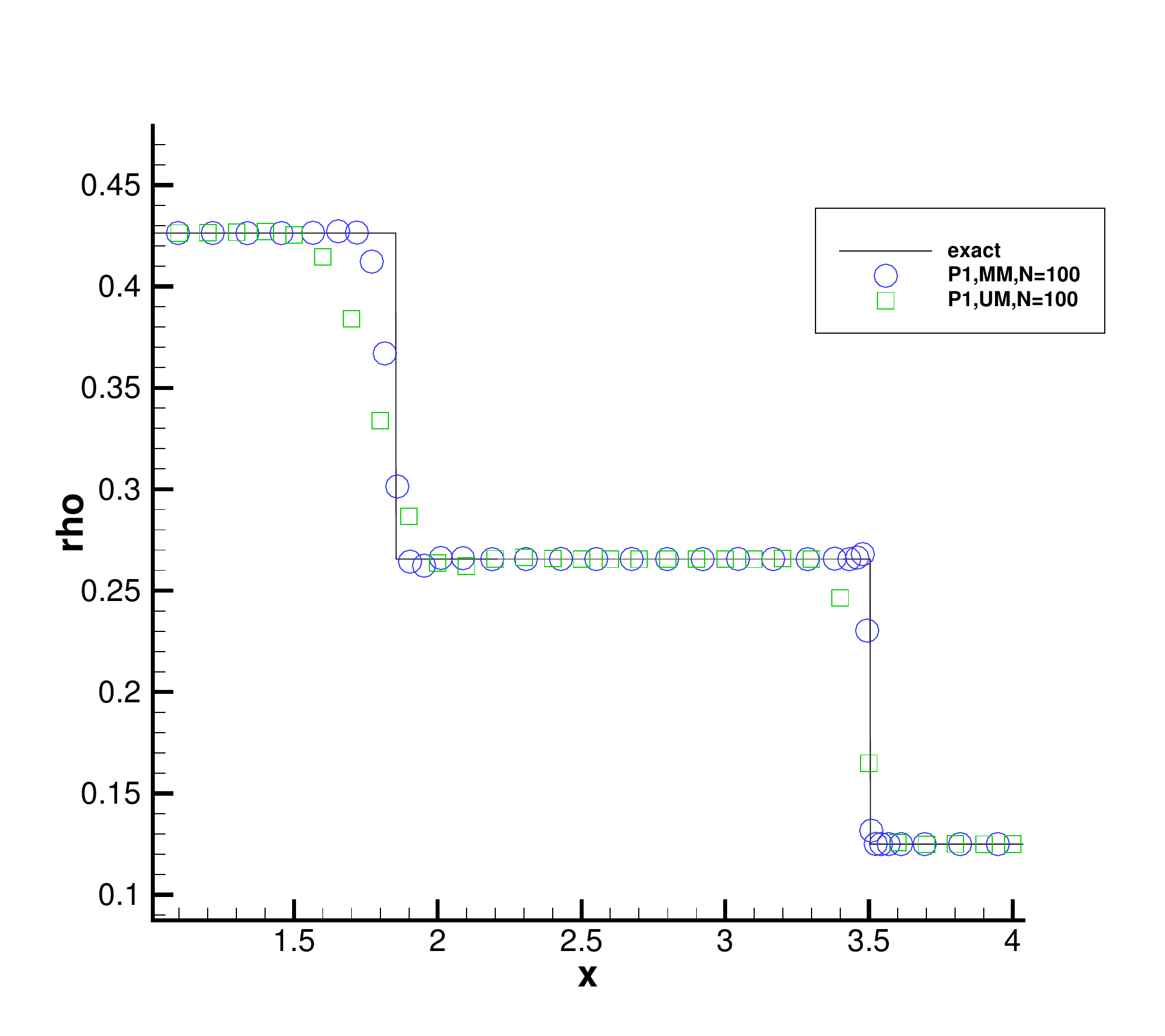}}
   }
 \mbox{\subfigure[MM: $N=100$, UM: $N=200$]
 {\includegraphics[width=8cm]{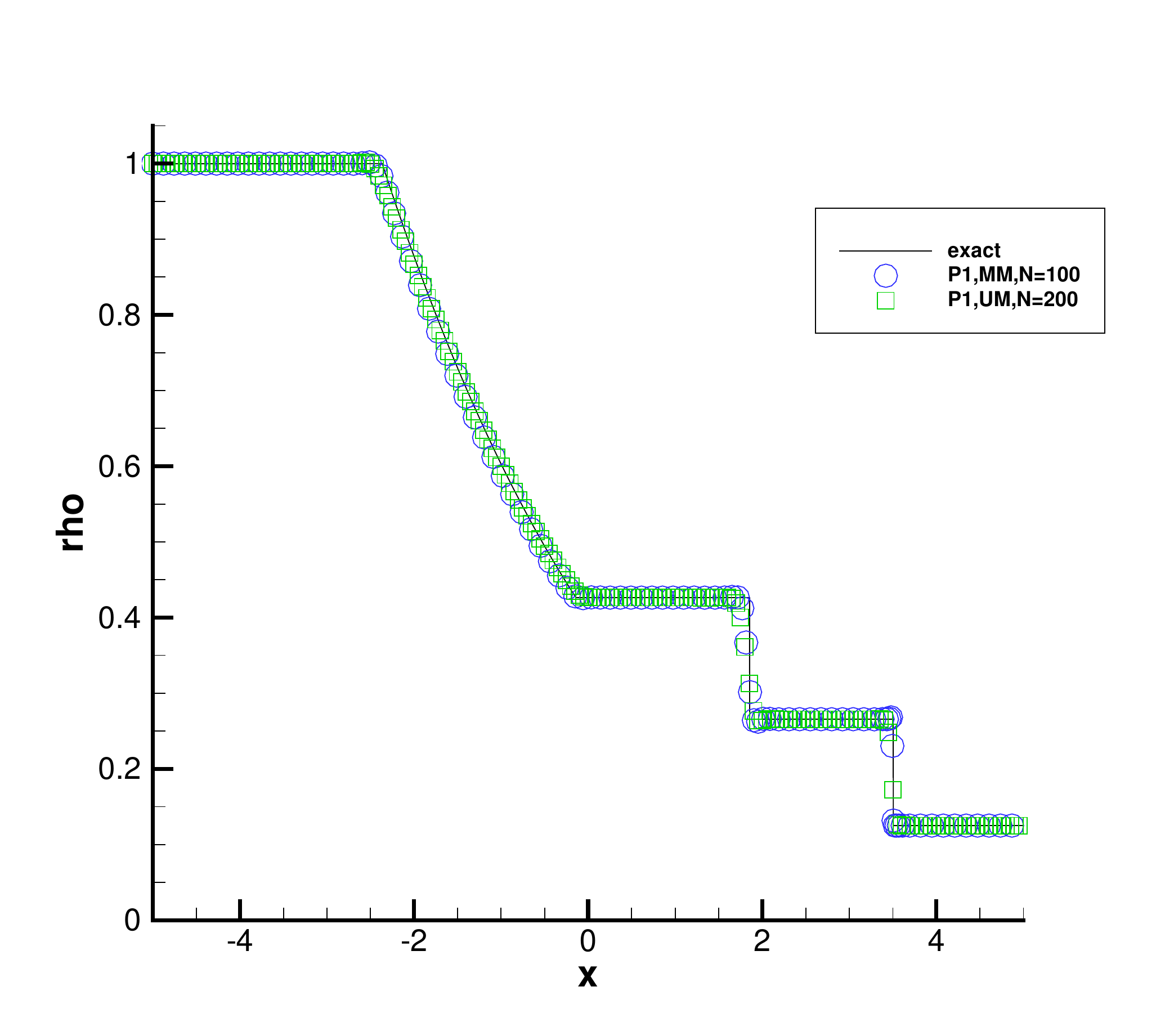}}\quad
   \subfigure[close view of (c) near shock]
   {\includegraphics[width=8cm]{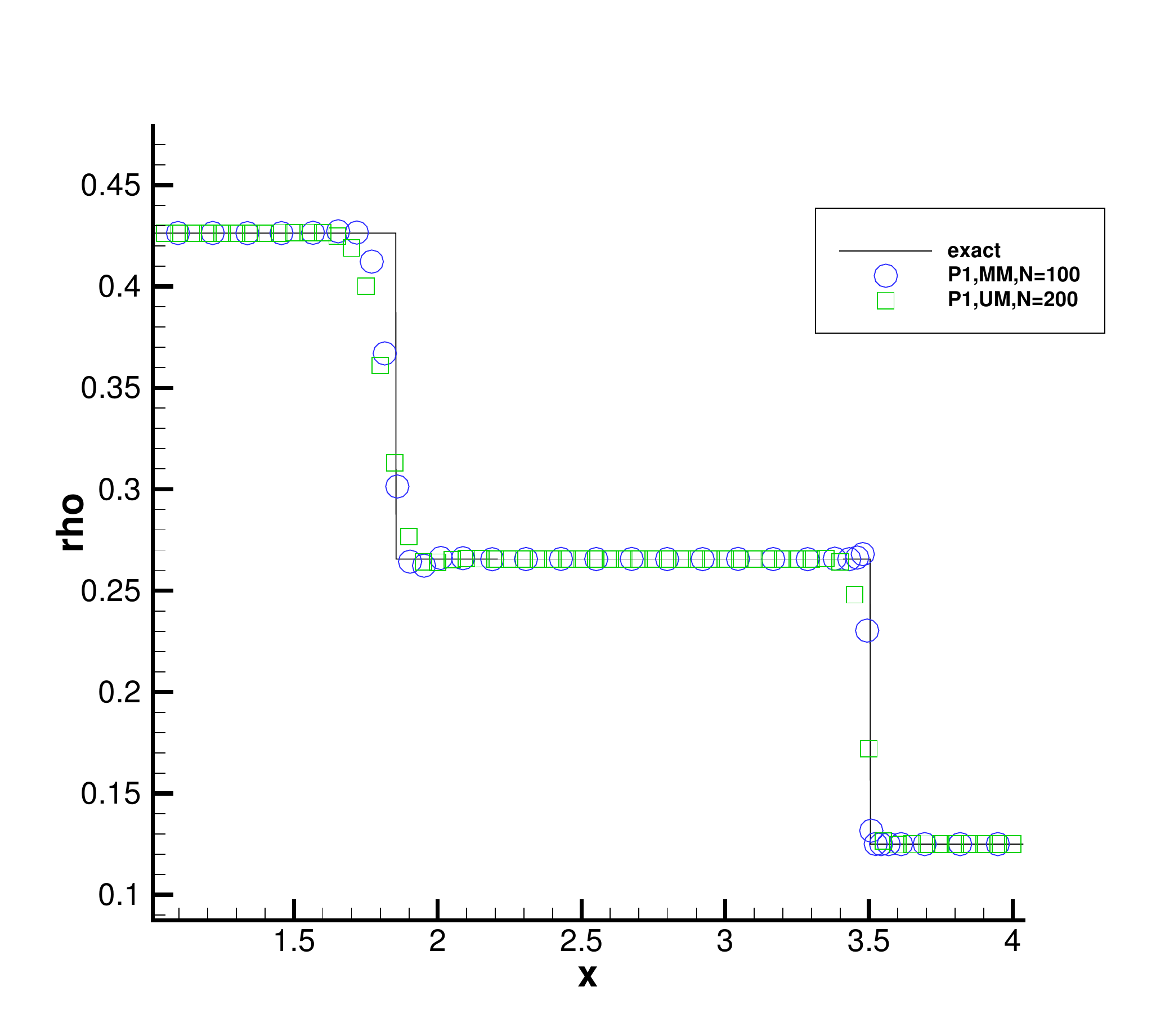}}
   }
   \mbox{\subfigure[MM: $N=100$, UM: $N=400$]
 {\includegraphics[width=8cm]{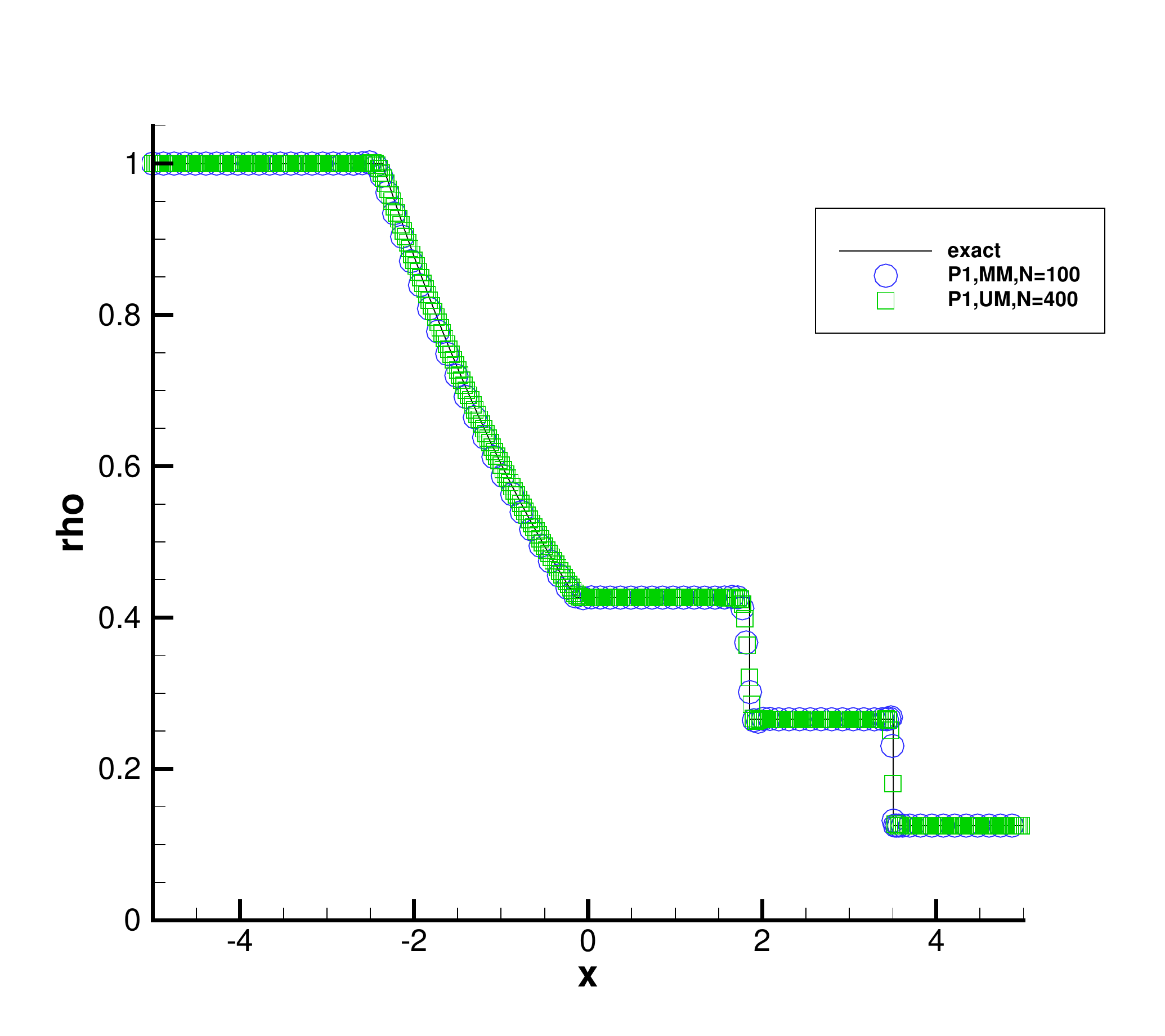}}\quad
   \subfigure[close view of (e) near shock]
   {\includegraphics[width=8cm]{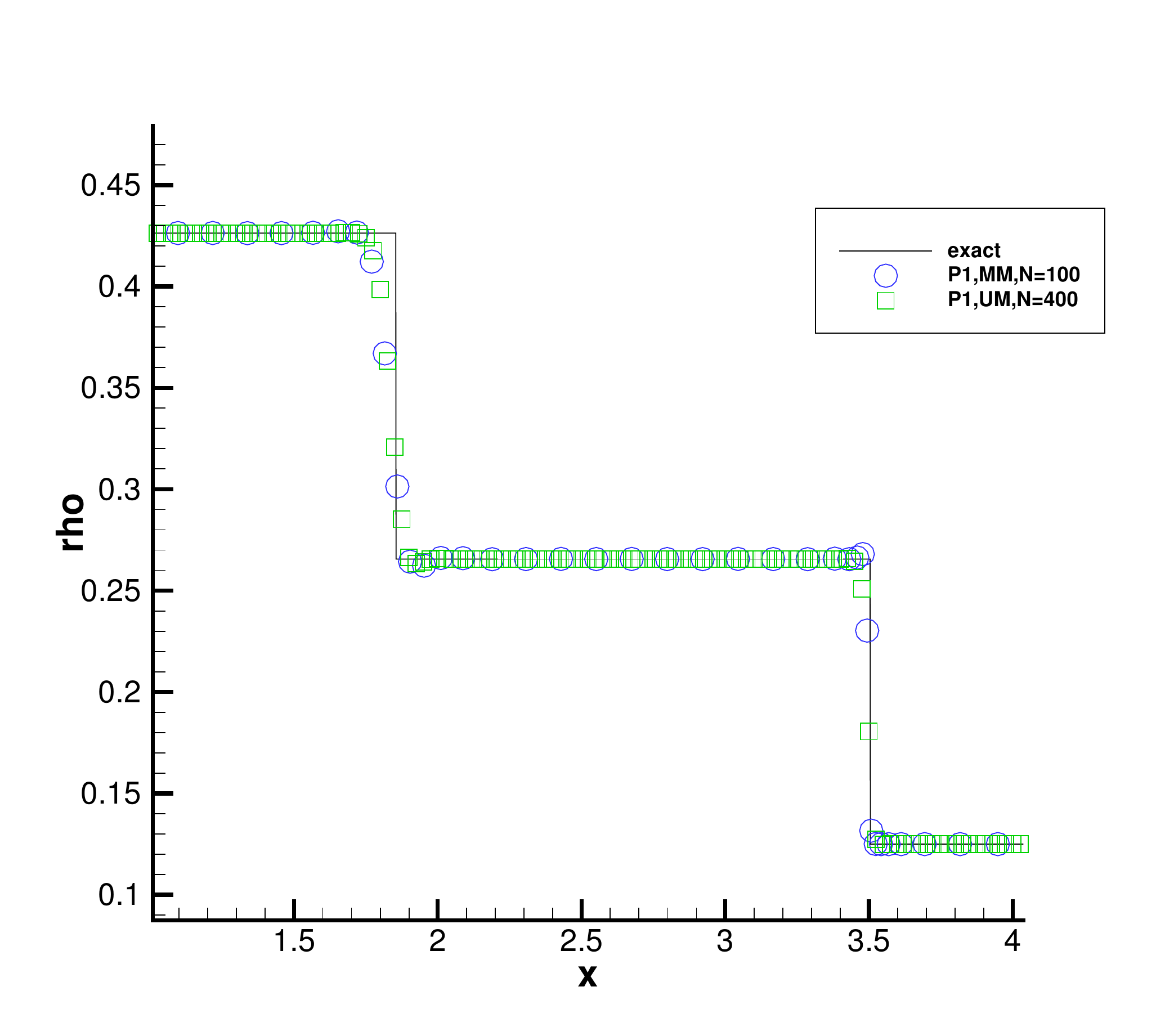}}
   }

   \caption{Example~\ref{exam4.3} (Sod Problem). The moving mesh solution (density) with $N=100$ is compared with the uniform mesh solutions with $N=100$, $200$, and $400$. $P^1$ elements are used.}
   \label{fig:edge1}
   \end{center}
   \end{figure}

\begin{figure}[hbtp]
 \begin{center}
 \mbox{\subfigure[MM: $N=100$, UM: $N=100$]
 {\includegraphics[width=8cm]{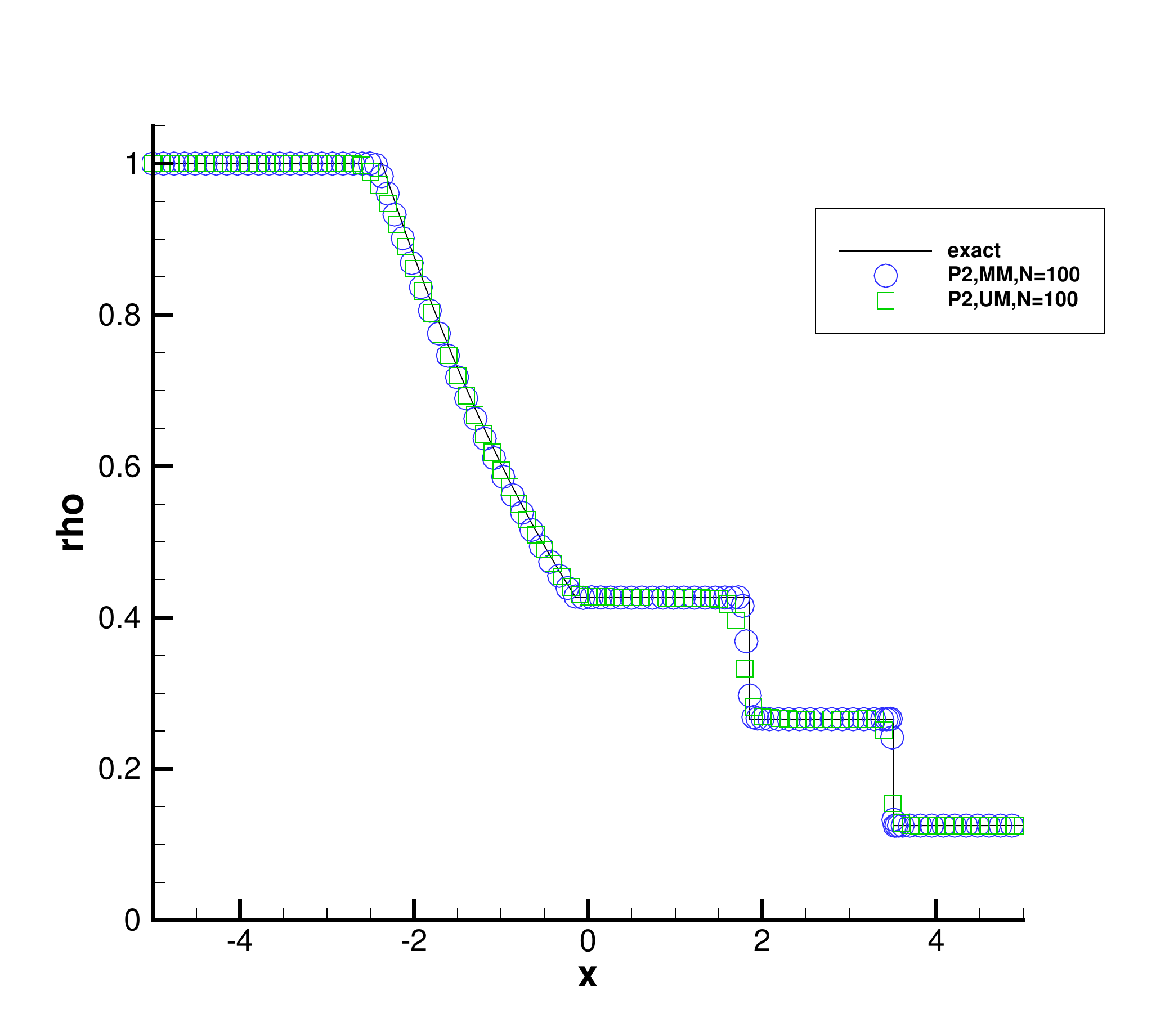}}\quad
   \subfigure[close view of (a) near shock]
   {\includegraphics[width=8cm]{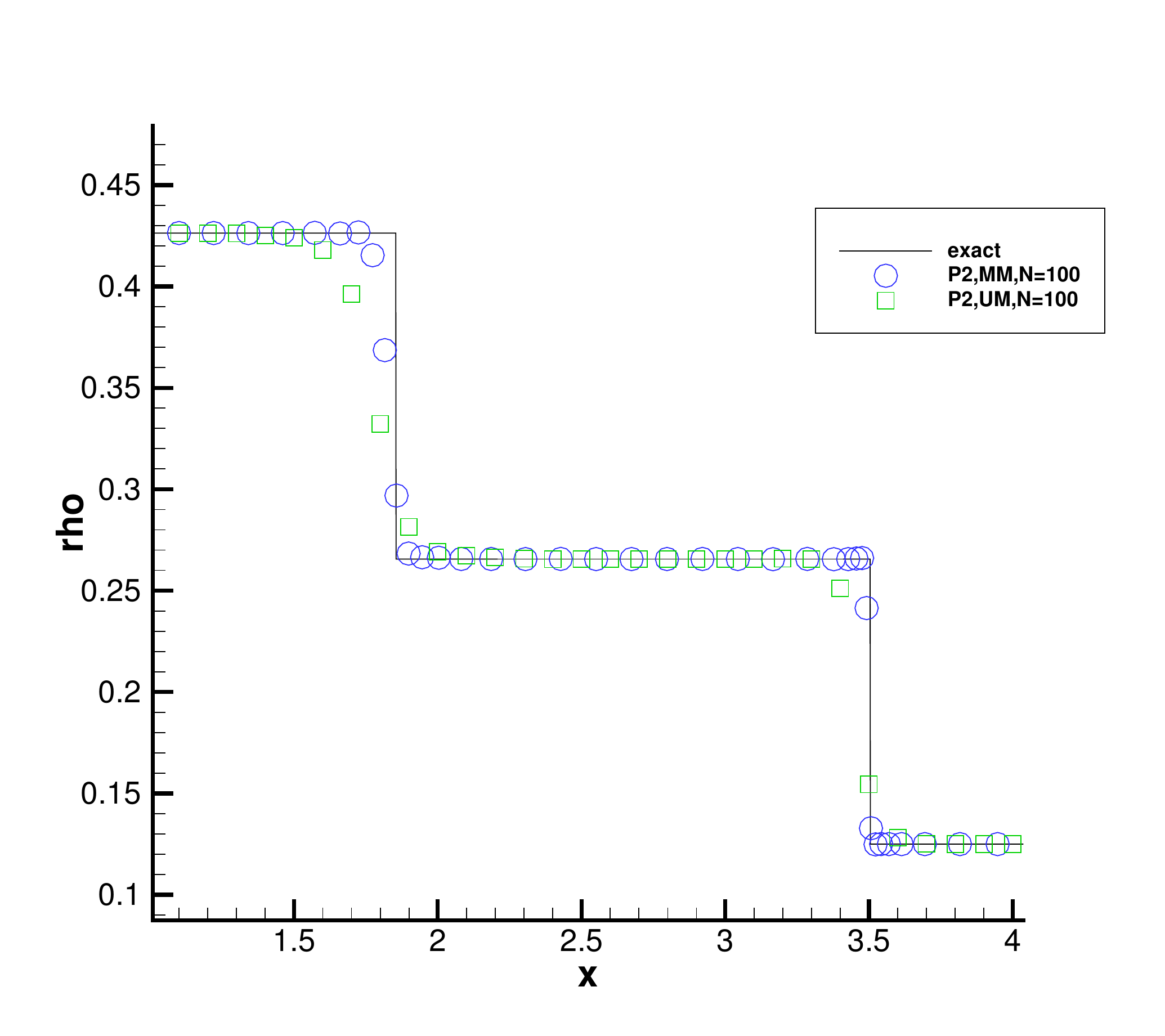}}
   }
 \mbox{\subfigure[MM: $N=100$, UM: $N=200$]
 {\includegraphics[width=8cm]{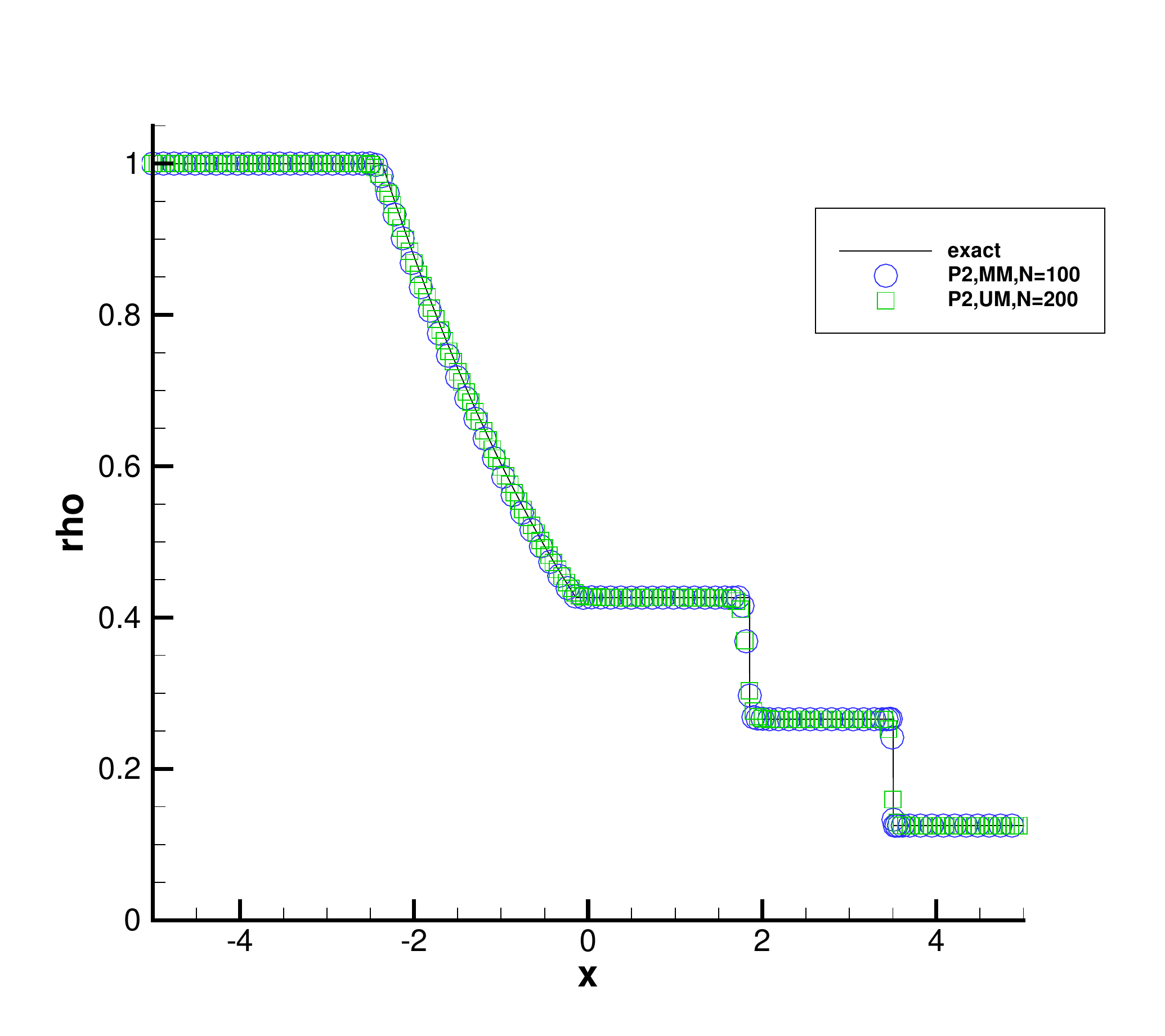}}\quad
   \subfigure[close view of (c) near shock]
   {\includegraphics[width=8cm]{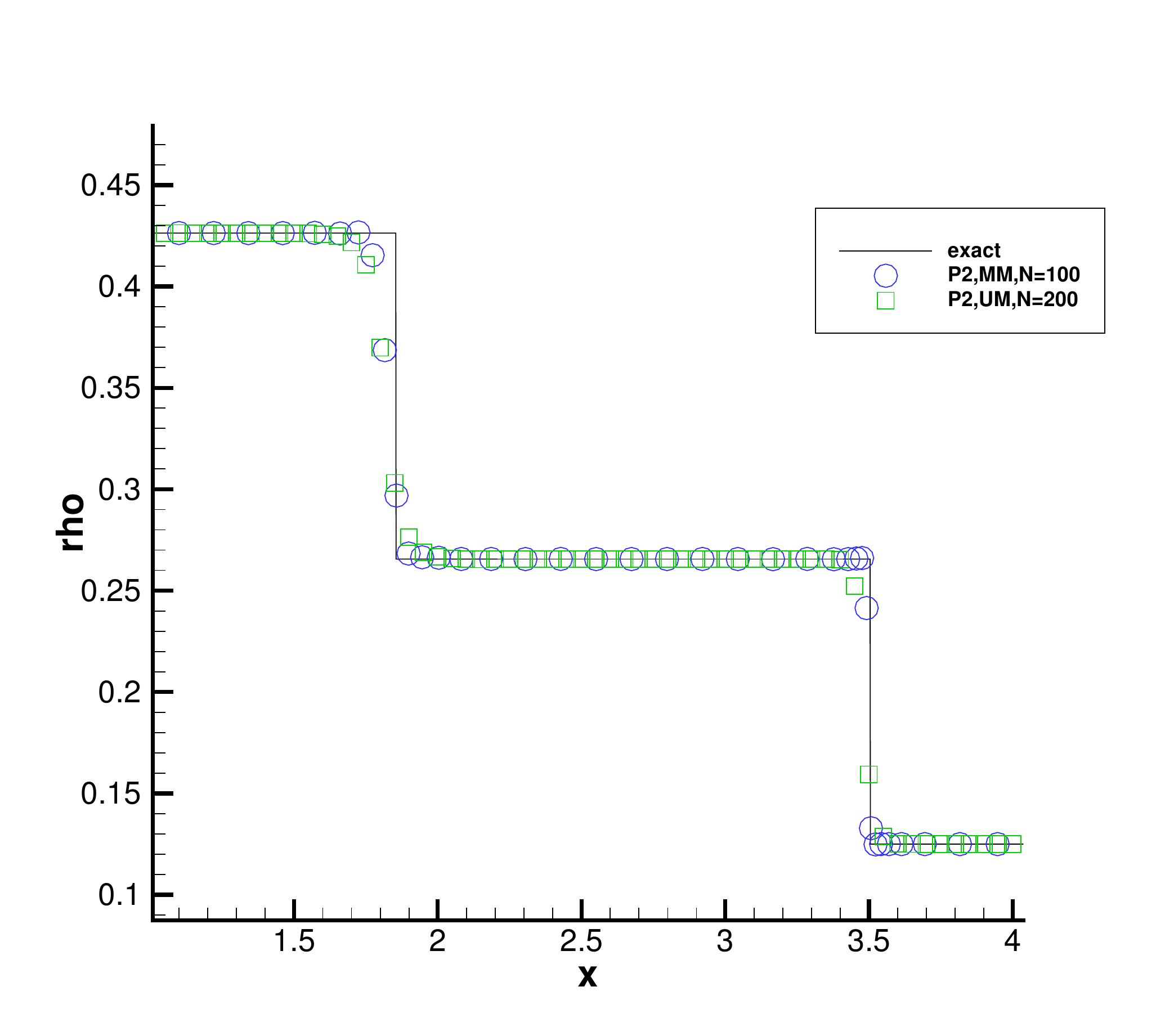}}
   }
   \mbox{\subfigure[MM: $N=100$, UM: $N=400$]
 {\includegraphics[width=8cm]{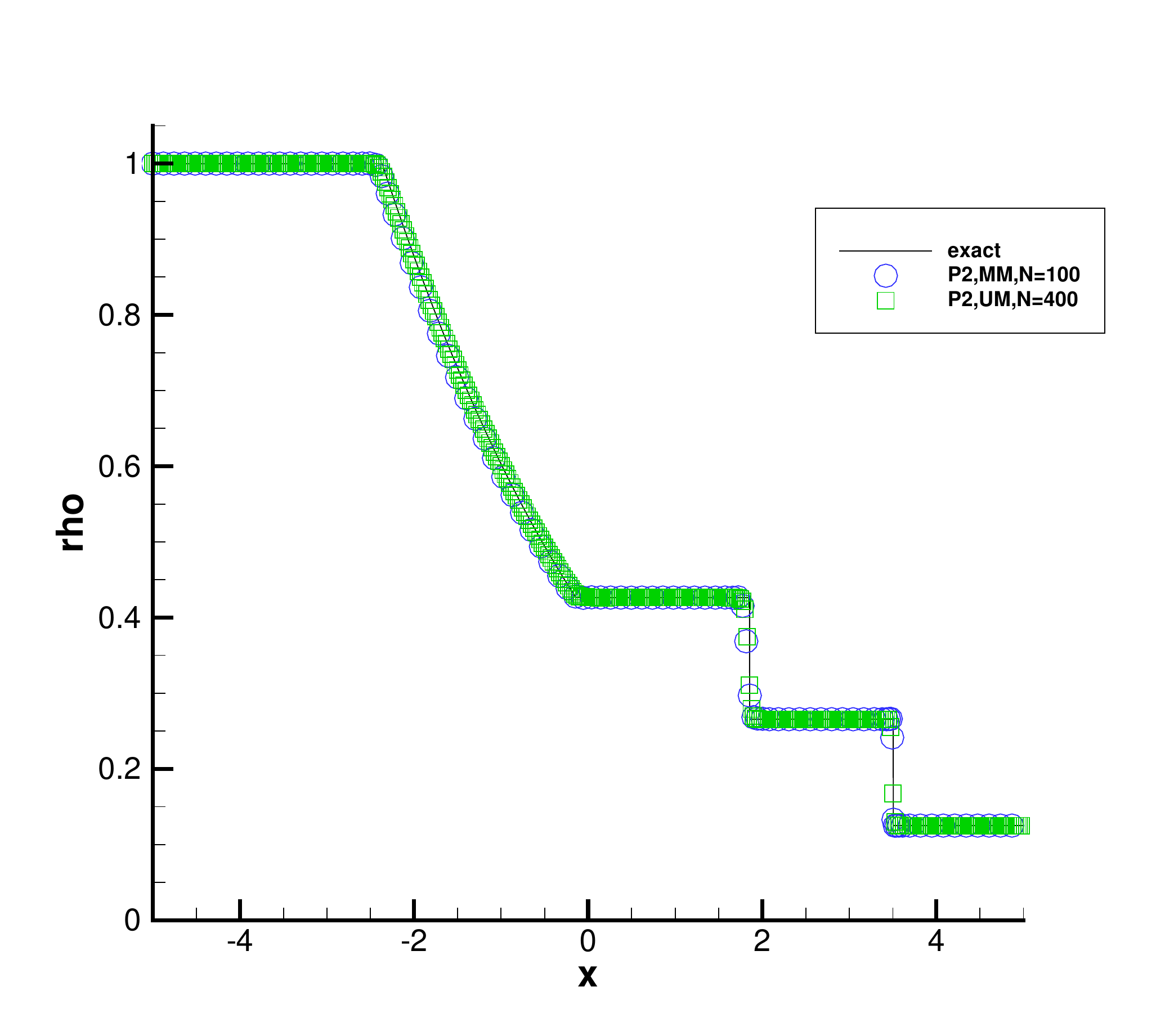}}\quad
   \subfigure[close view of (e) near shock]
   {\includegraphics[width=8cm]{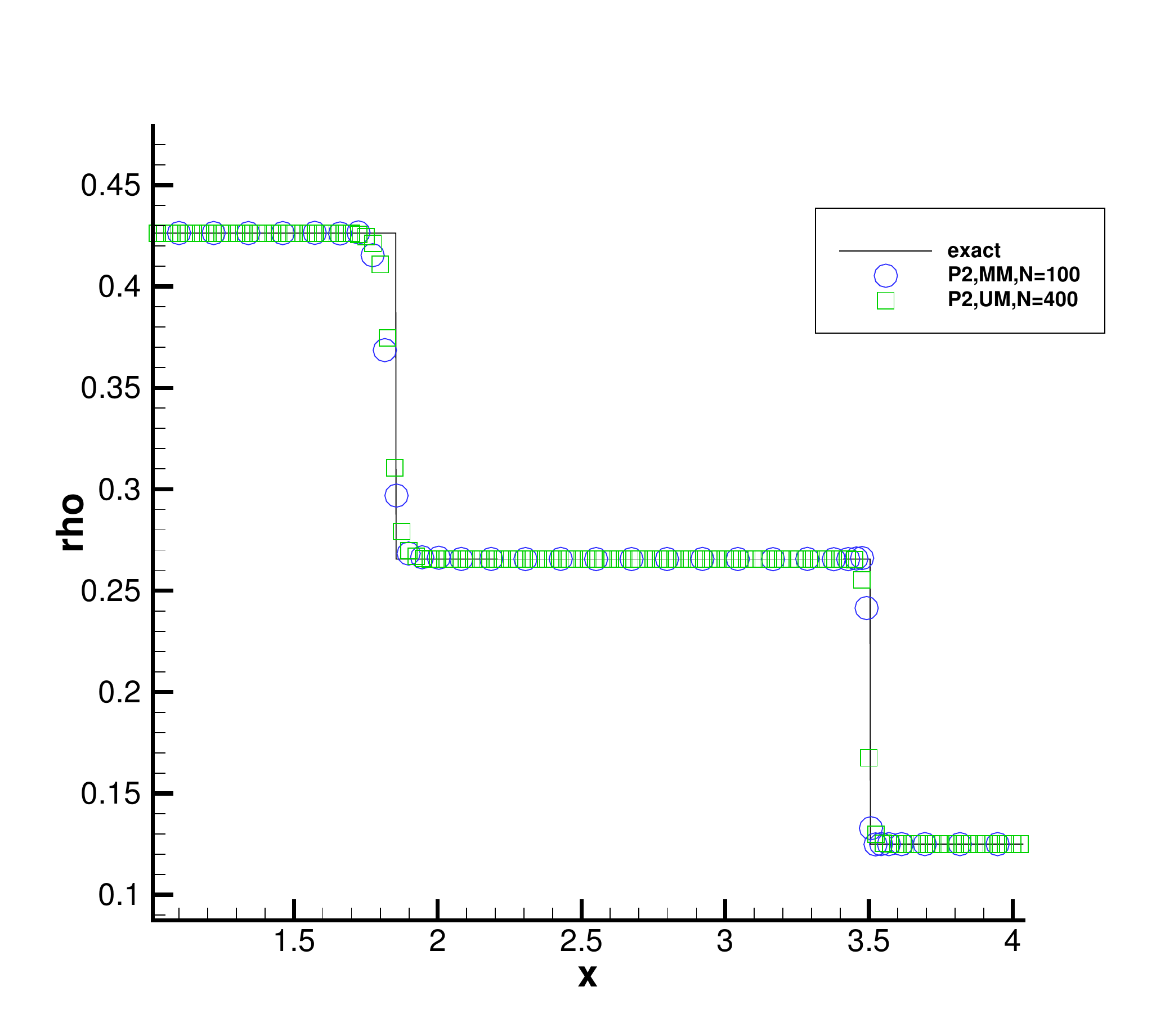}}
   }

   \caption{Example~\ref{exam4.3} (Sod Problem). The moving mesh solution (density) with $N=100$ is compared with the uniform mesh solutions with $N=100$, $200$, and $400$. $P^2$ elements are used.}
   \label{fig:edge2}
   \end{center}
   \end{figure}

\begin{figure}[hbtp]
 \begin{center}
 \mbox{\subfigure[$P^1$ elements]
 {\includegraphics[width=8cm]{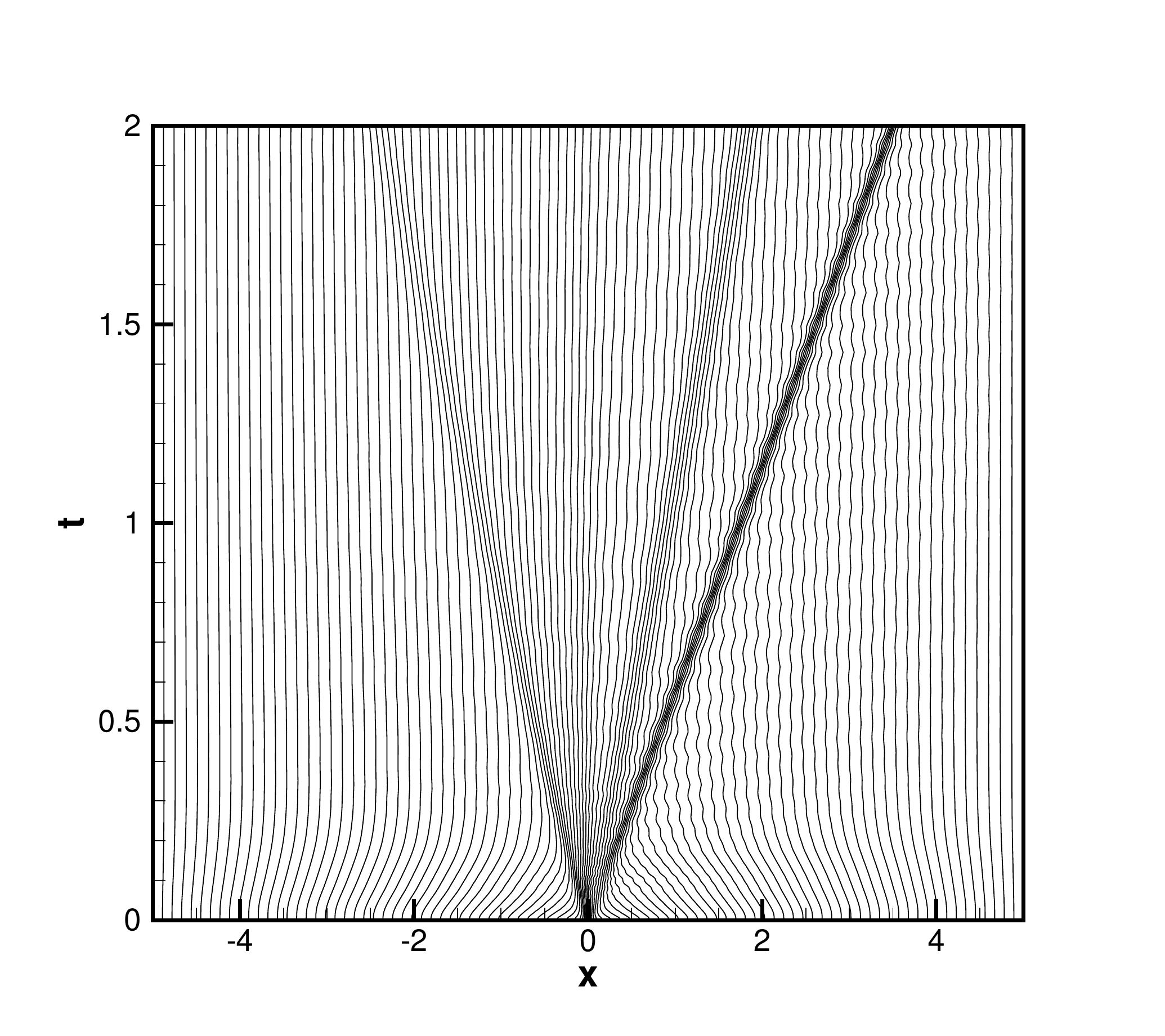}}
  \subfigure[$P^2$ elements]
  {\includegraphics[width=8cm]{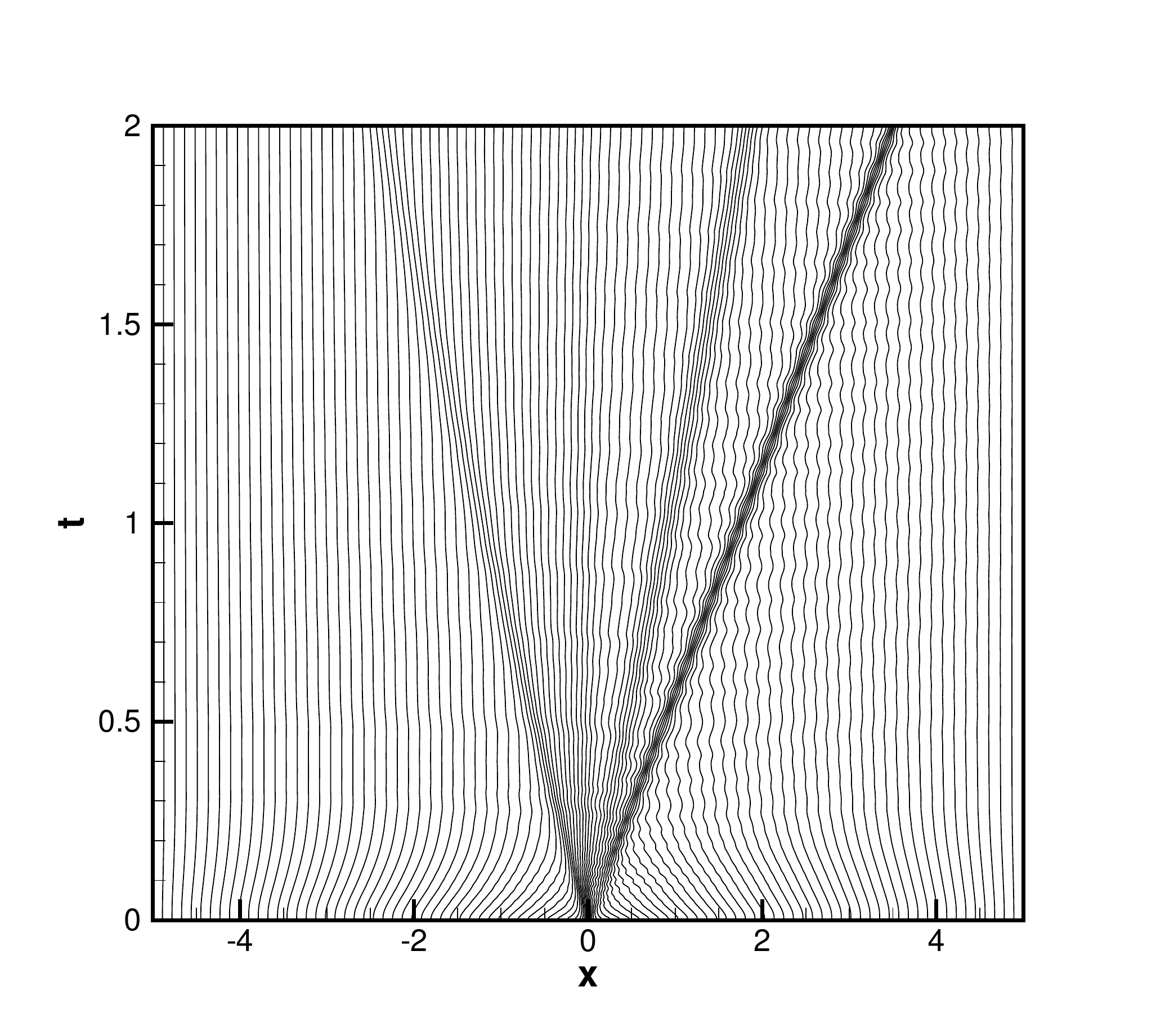}}
  }
 \caption{Example~\ref{exam4.3} (Sod Problem). The trajectories of a moving mesh with $N=100$ are plotted.} \label{trfig1}
  \end{center}
   \end{figure}

\begin{exam}{\em
\label{exam4.3}
This example is the Sod problem of the Euler equation (\ref{Euler}) subject to the inflow/outflow boundary condition and a Riemann initial condition
\begin{equation}
(\rho,u,p)=
\left
\{
\begin{array}{ll}
(1,0,1),      \quad&\text{for}\quad x<0\\
(0.125,0,0.1), \quad&\text{for}\quad x>0.\notag
\end{array}
\right.
\end{equation}
The computational domain is $(-5,5)$ and the final time is $T=2.0$. 

The moving mesh DG solution (density) obtained with $N=100$ is compared with the uniform mesh solutions obtained with $N=100,\; 200$, and $400$ in Figs. \ref{fig:edge1} and \ref{fig:edge2}. One can see that the moving mesh solutions are more accurate than the uniform mesh solutions for the same number of points and comparable with those with $N=400$ for both $P^1$ and $P^2$ elements.

The trajectories of a moving mesh are plotted in Fig. \ref{trfig1}. From the figure, one can observe that the points are concentrated at $x=0$ initially where the initial condition is discontinuous. As time evolves, the moving mesh can capture not only the shock but also the contact discontinuity well. In addition, the points are also clustered at the front and the tail of the rarefaction since the Hessian is used in the computation of the metric tensor $\mathbb{M}$.


}\end{exam}

%
%
%

%
%
%

\begin{figure}[hbtp]
 \begin{center}
 \mbox{\subfigure[MM: $N=100$, UM: $N=100$]
 {\includegraphics[width=8cm]{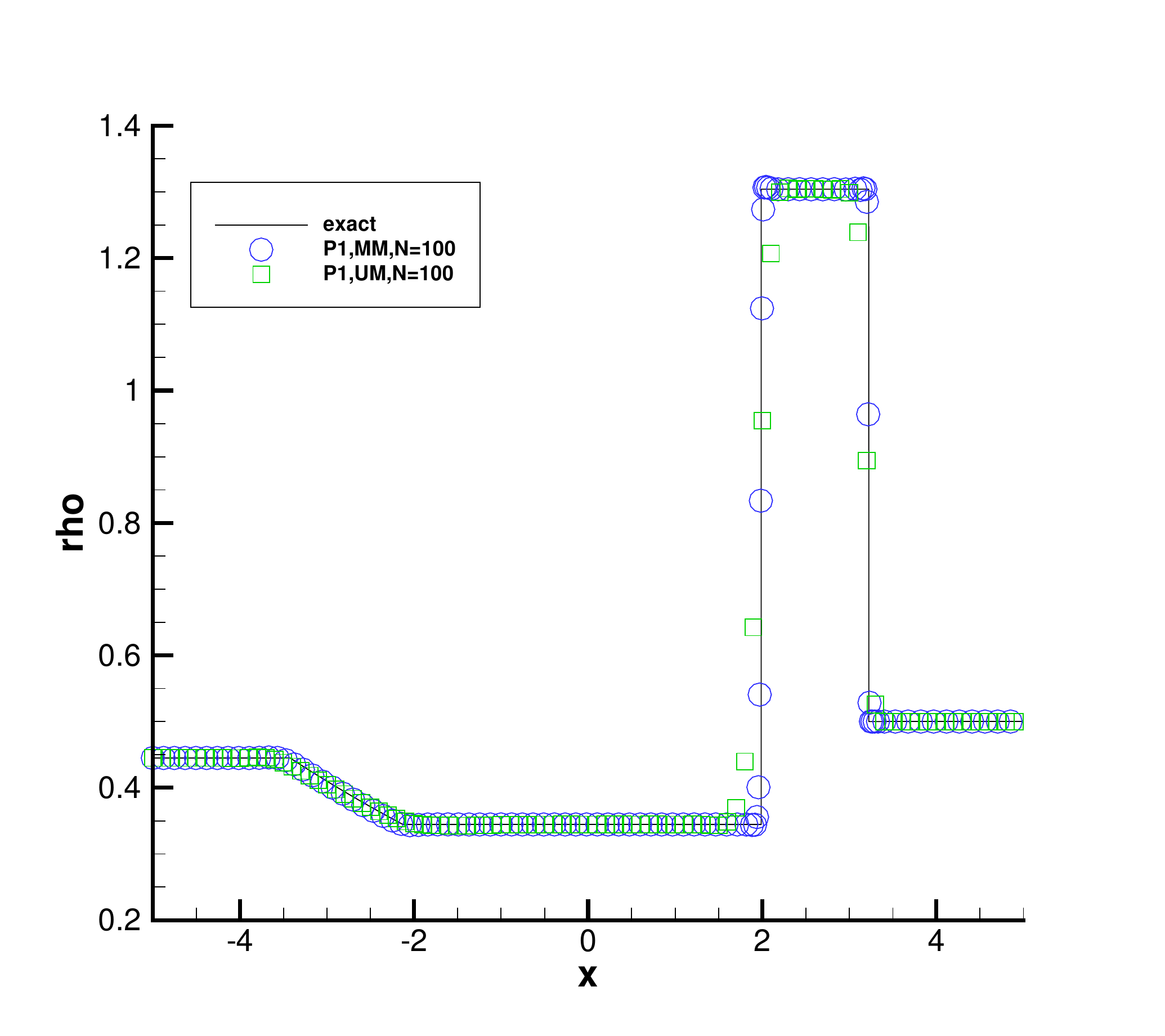}}\quad
   \subfigure[close view of (a) near shock]
   {\includegraphics[width=8cm]{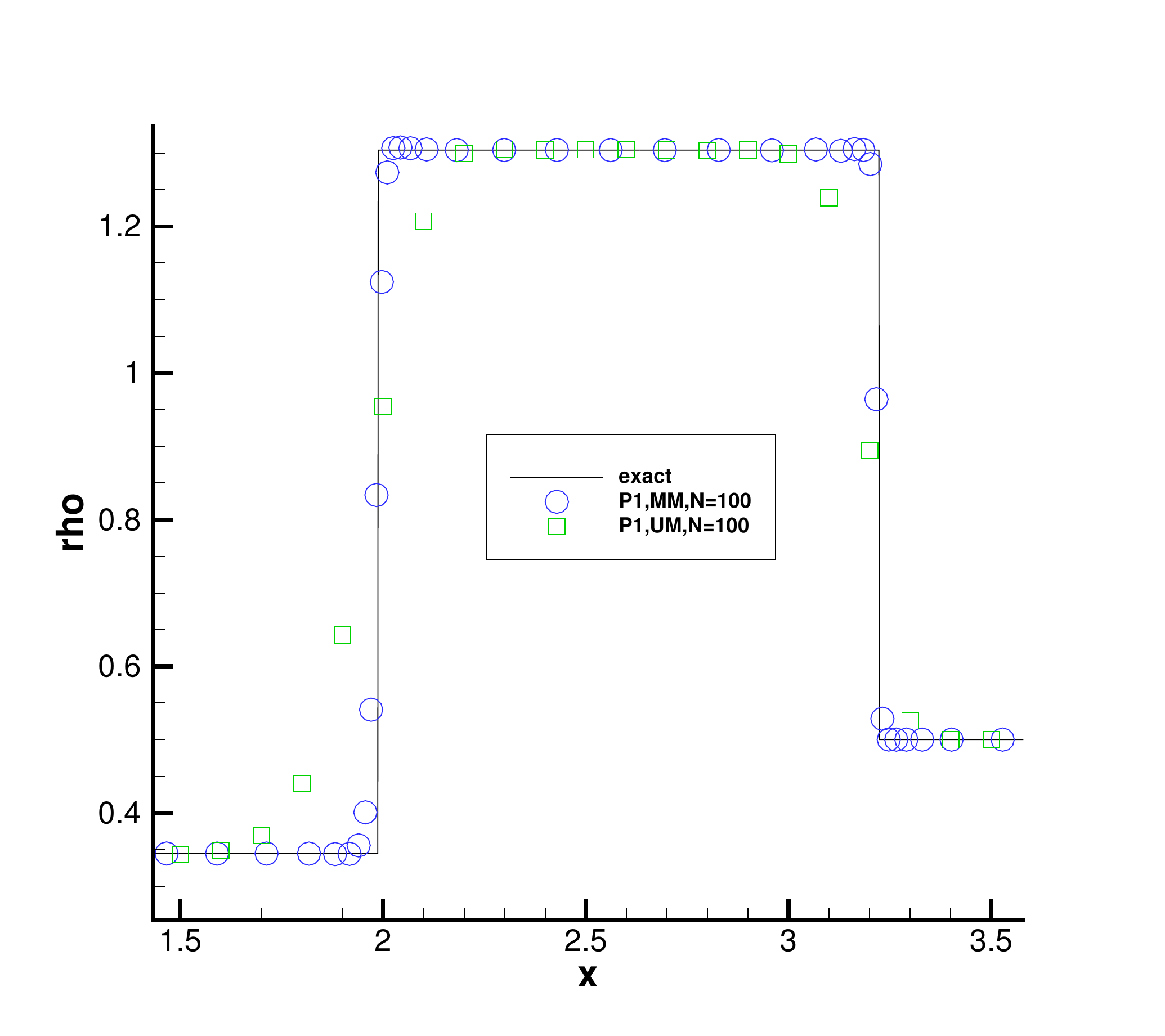}}
   }
 \mbox{\subfigure[MM: $N=100$, UM: $N=200$]
 {\includegraphics[width=8cm]{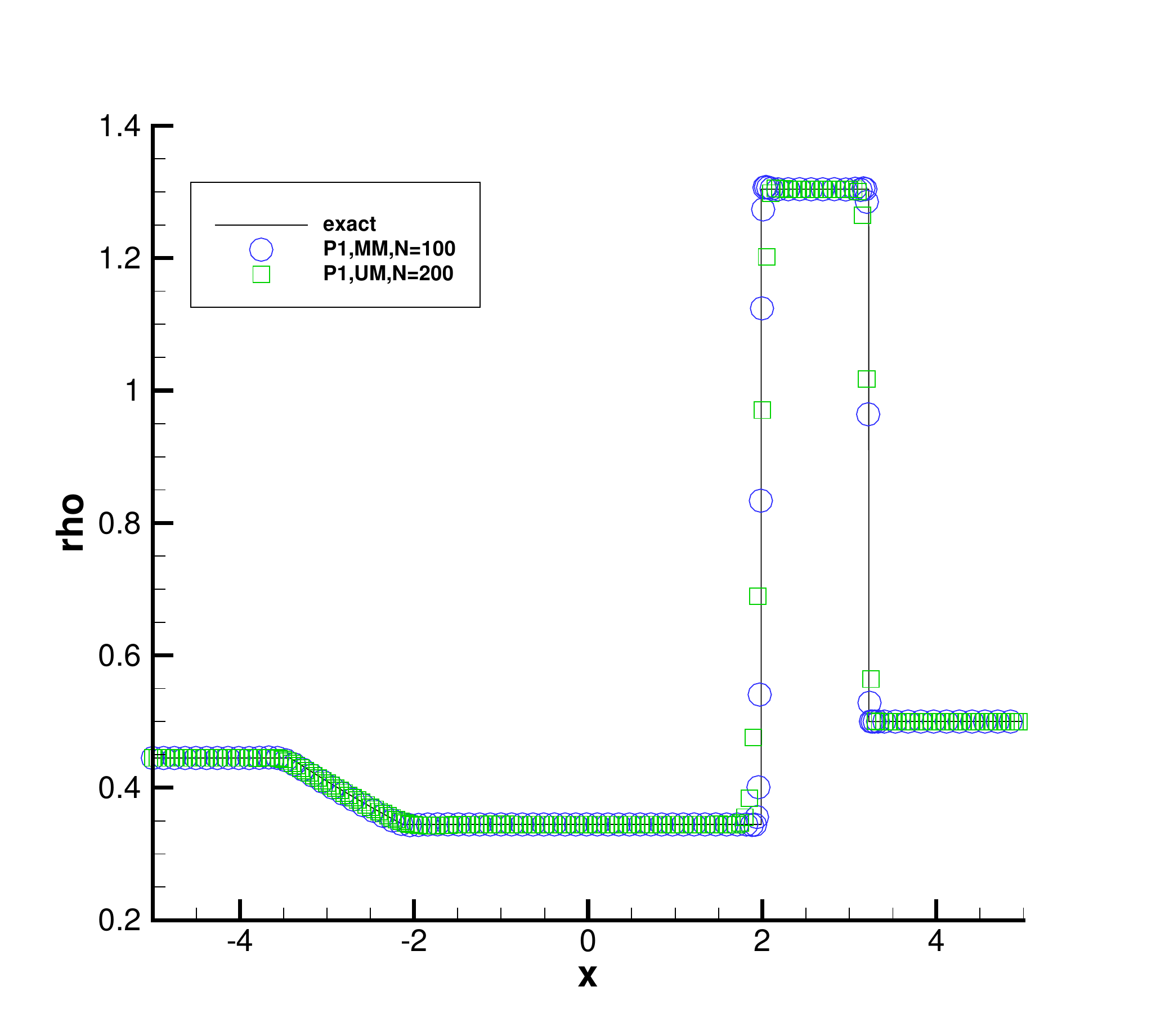}}\quad
   \subfigure[close view of (c) near shock]
    {\includegraphics[width=8cm]{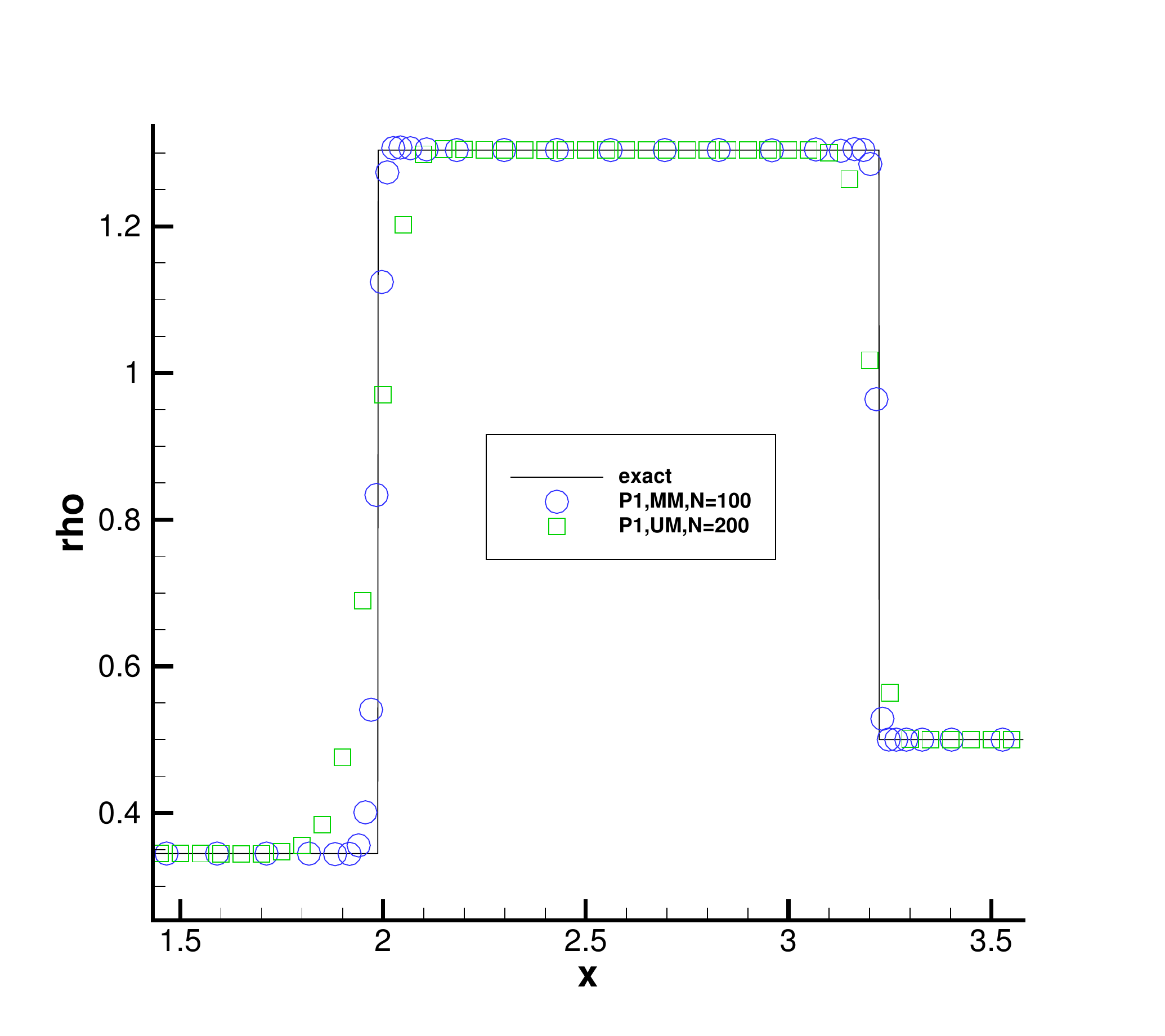}}
   }
   \mbox{\subfigure[MM: $N=100$, UM: $N=400$]
 {\includegraphics[width=8cm]{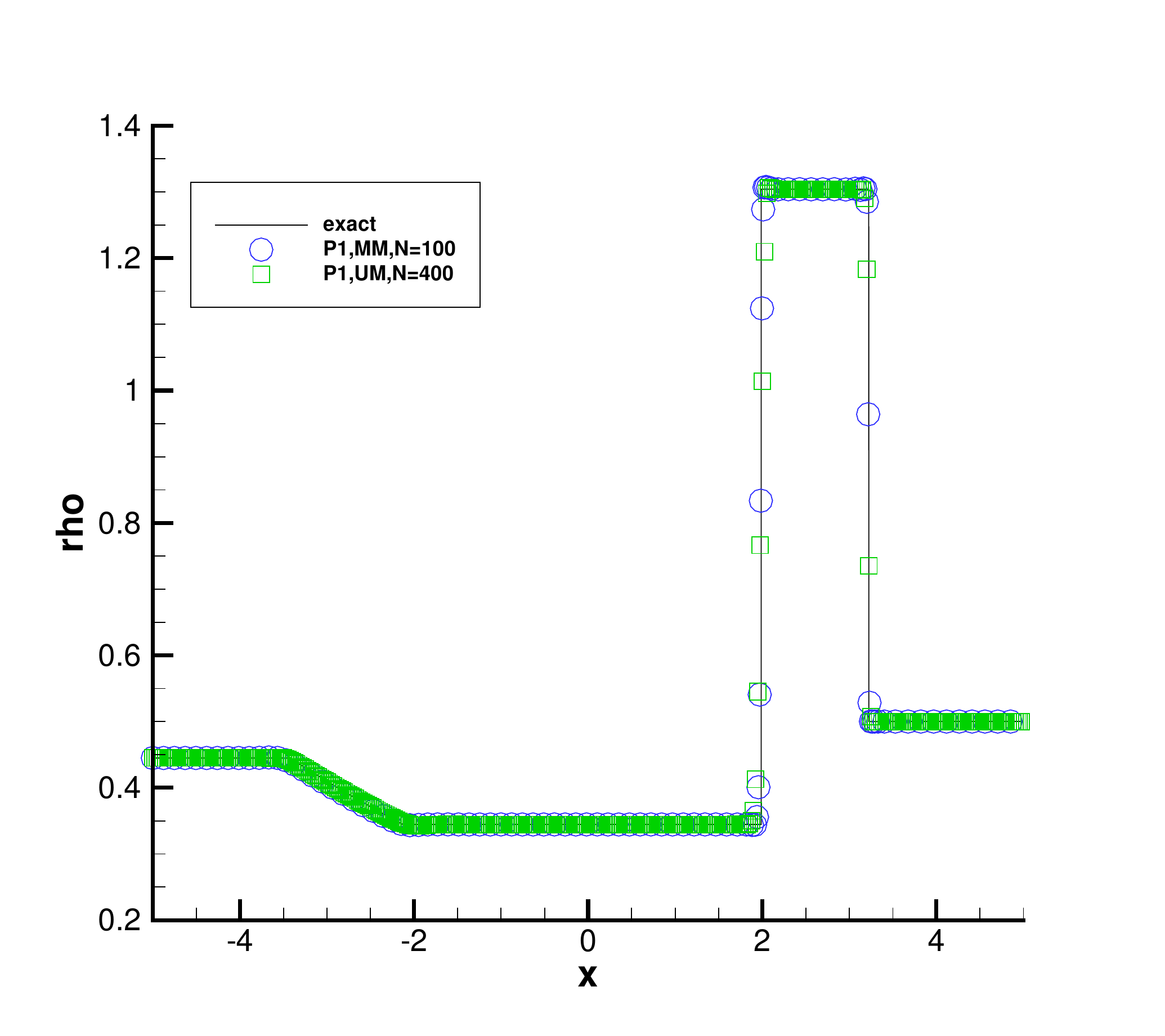}}\quad
   \subfigure[close view of (e) near shock]
   {\includegraphics[width=8cm]{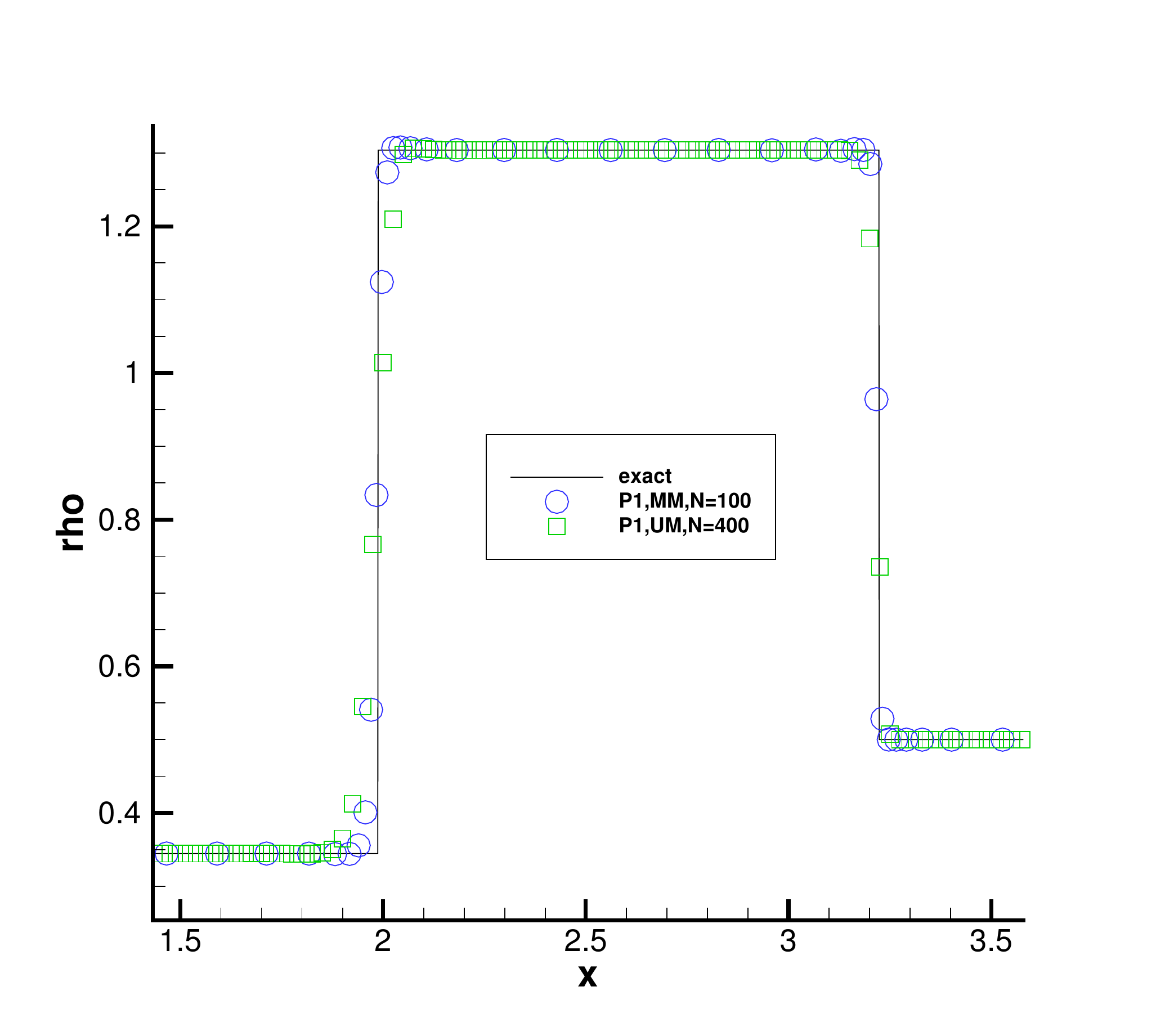}}
   }
   \caption{Example~\ref{exam4.4} (Lax Problem). The moving mesh solution (density) with $N=100$ is compared with the uniform mesh solutions with $N=100$, $200$, and $400$. $P^1$ elements are used.}
   \label{fig:edge3}
   \end{center}
   \end{figure}

\begin{figure}[hbtp]
 \begin{center}
 \mbox{\subfigure[MM: $N=100$, UM: $N=100$]
 {\includegraphics[width=8cm]{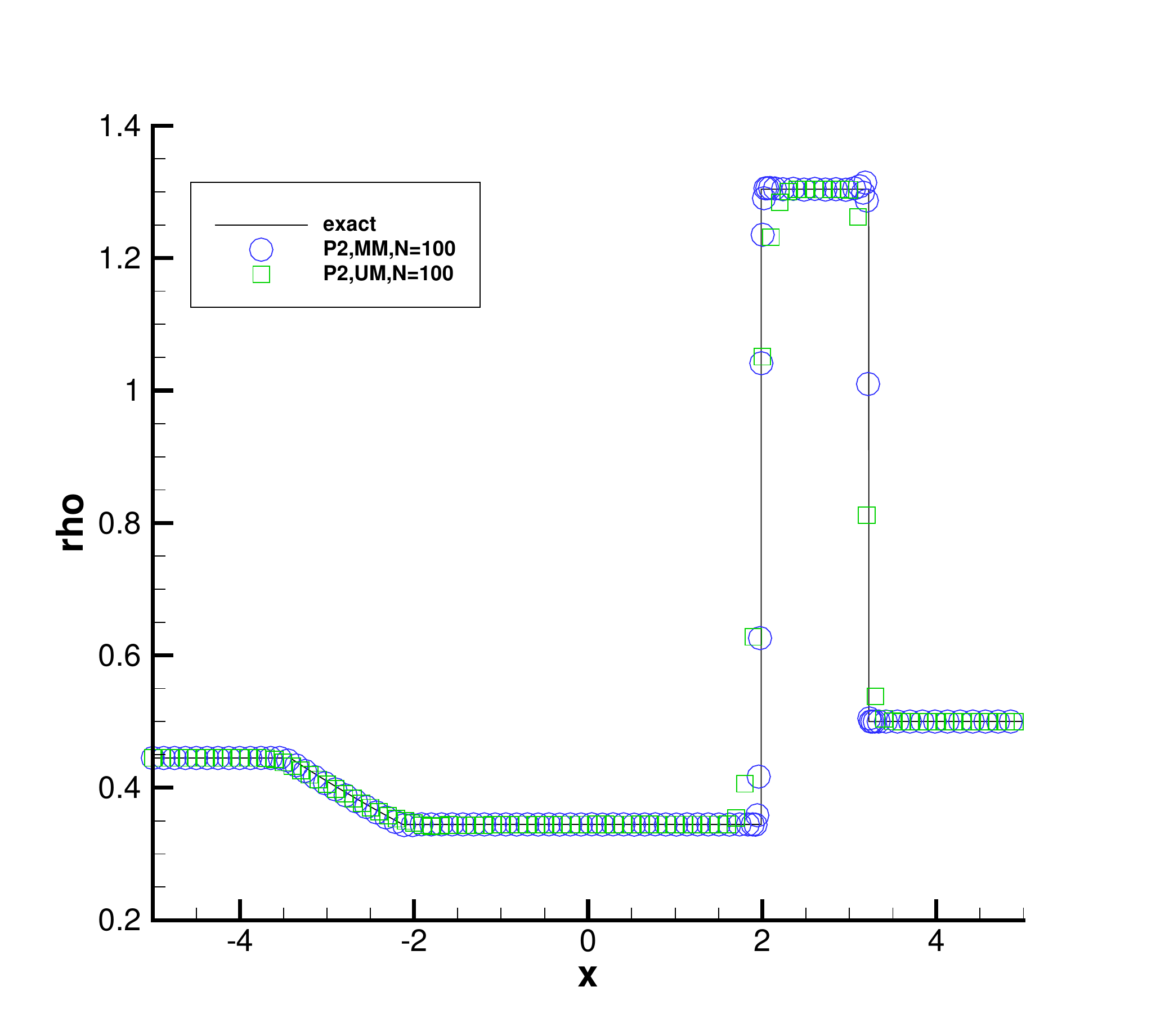}}\quad
   \subfigure[close view of (a) near shock]
   {\includegraphics[width=8cm]{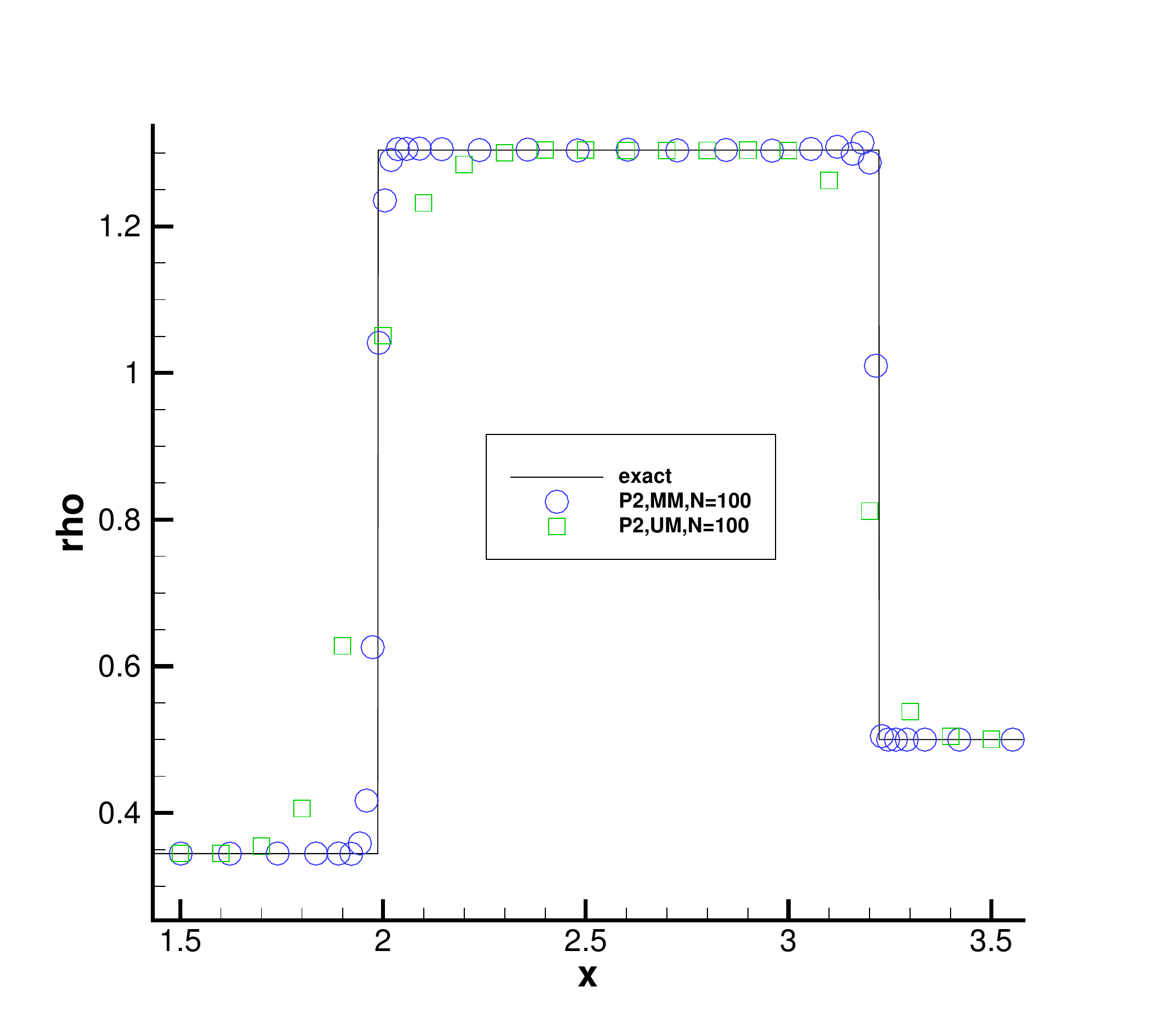}}
   }
 \mbox{\subfigure[MM: $N=100$, UM: $N=200$]
 {\includegraphics[width=8cm]{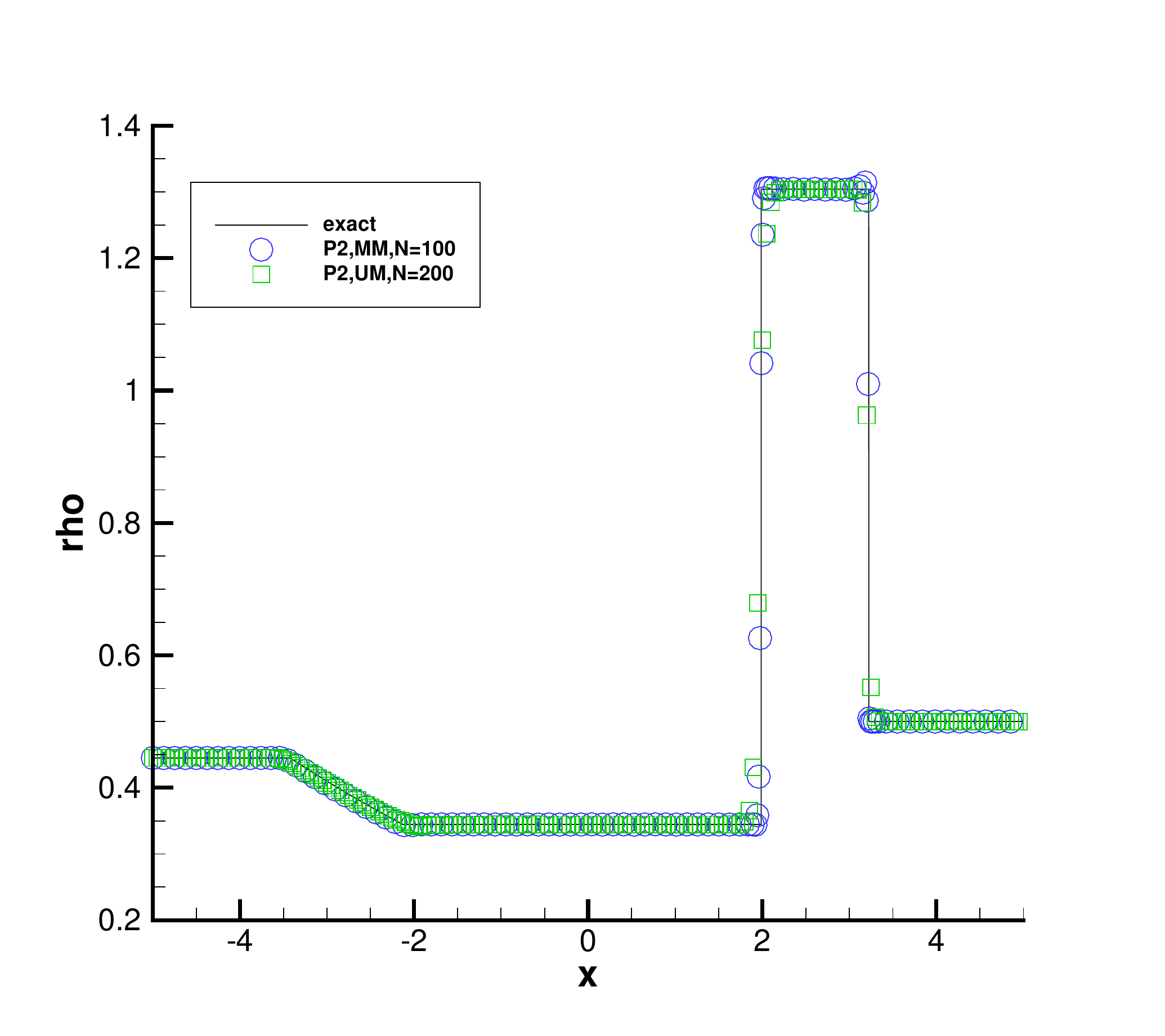}}\quad
   \subfigure[close view of (c) near shock]
   {\includegraphics[width=8cm]{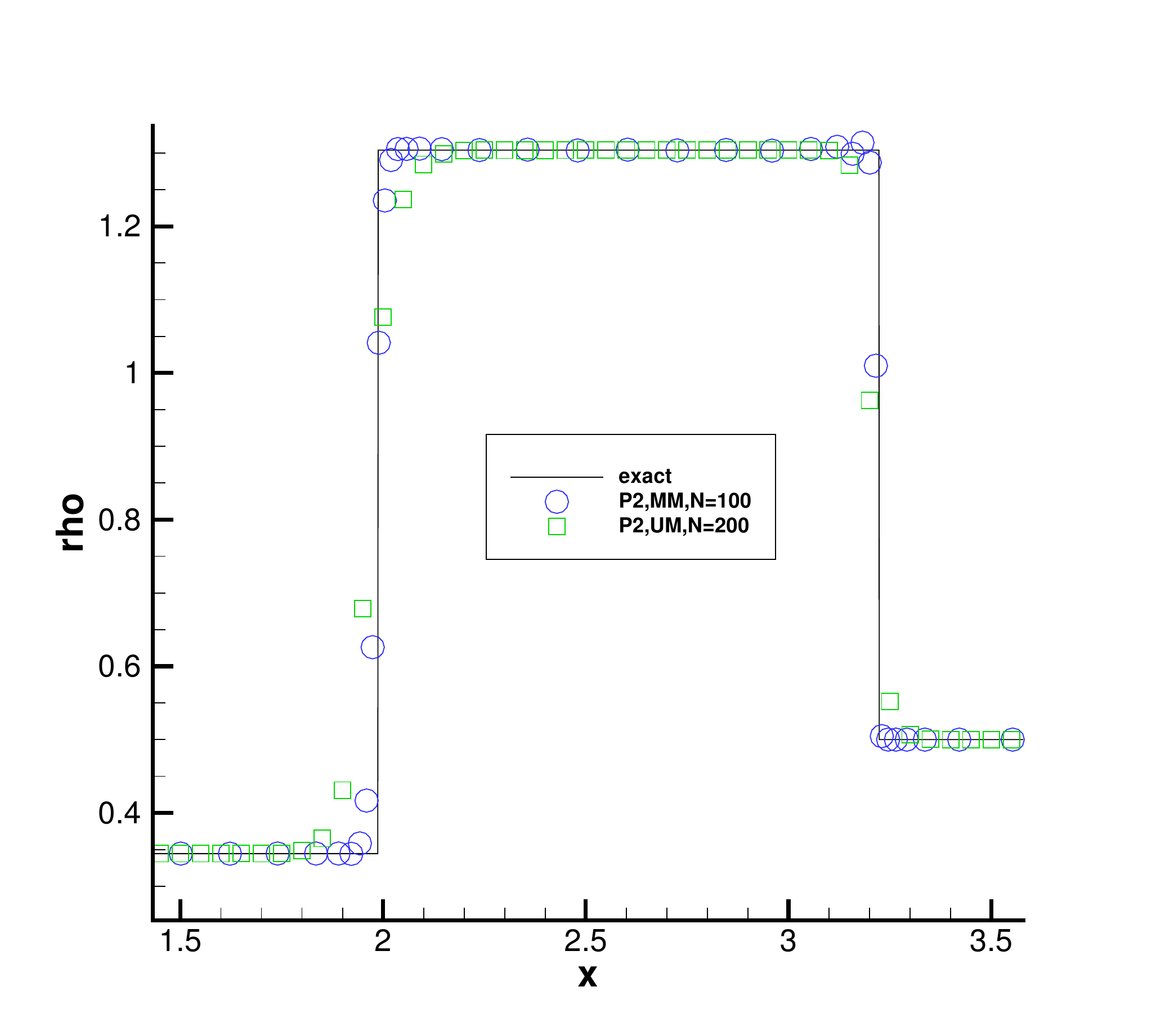}}
   }
   \mbox{\subfigure[MM: $N=100$, UM: $N=400$]
 {\includegraphics[width=8cm]{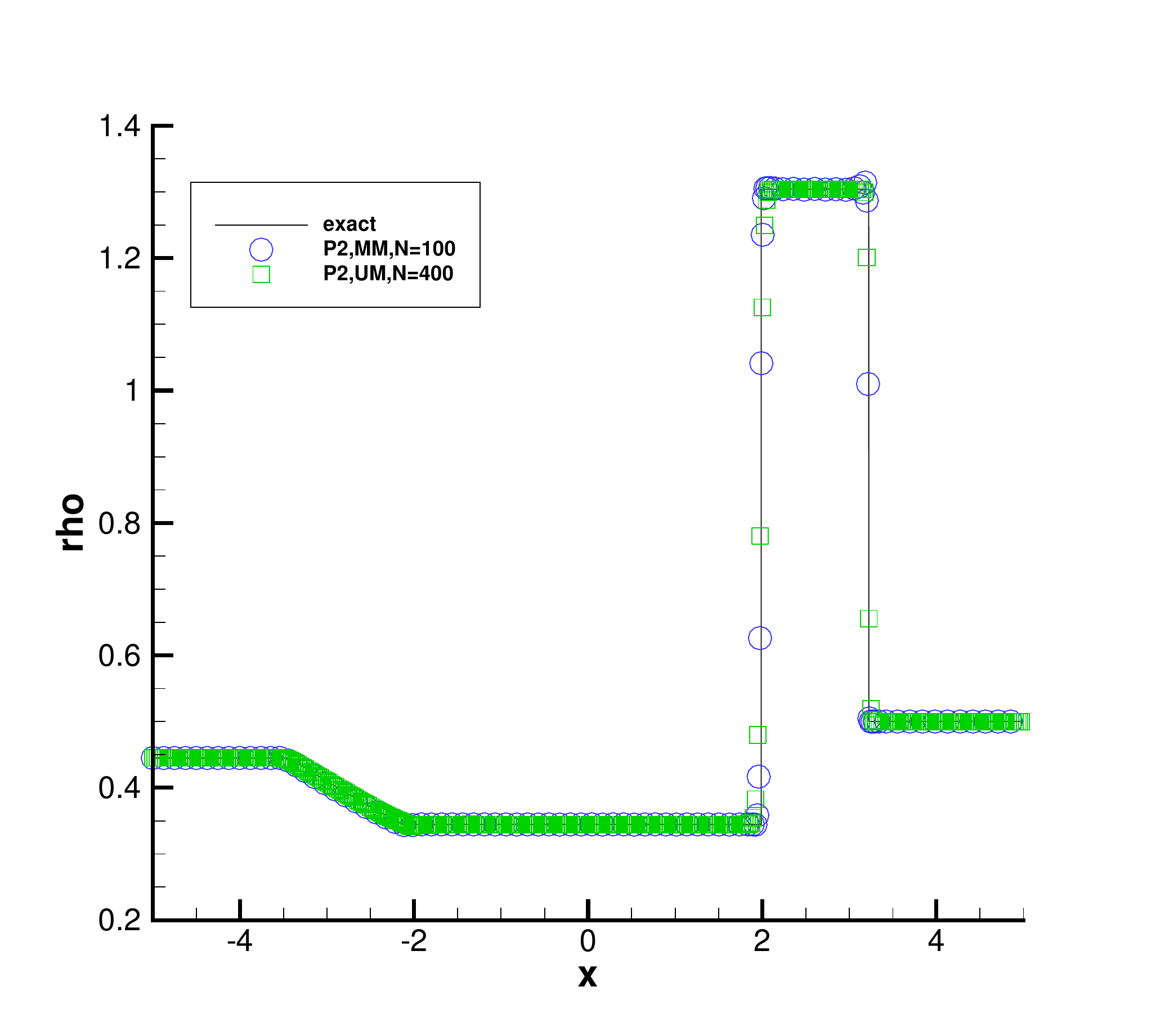}}\quad
   \subfigure[close view of (e) near shock]
   {\includegraphics[width=8cm]{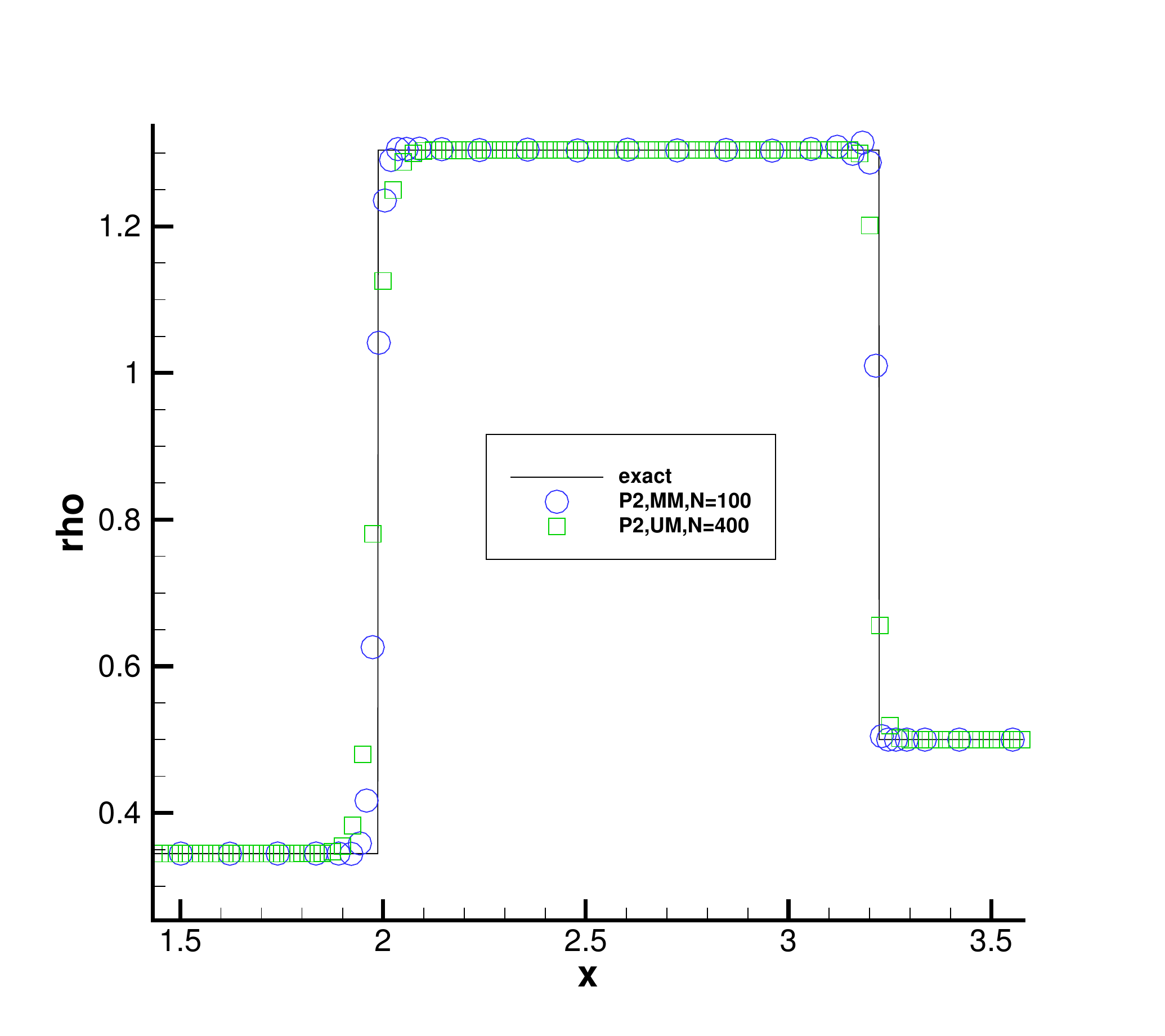}}
   }

   \caption{Example~\ref{exam4.4} (Lax Problem). The moving mesh solution (density) with $N=100$ is compared with the uniform mesh solutions with $N=100$, $200$, and $400$. $P^2$ elements are used.}
   \label{fig:edge4}
   \end{center}
   \end{figure}

\begin{figure}[hbtp]
 \begin{center}
 \mbox{\subfigure[$P^1$]
 {\includegraphics[width=8cm]{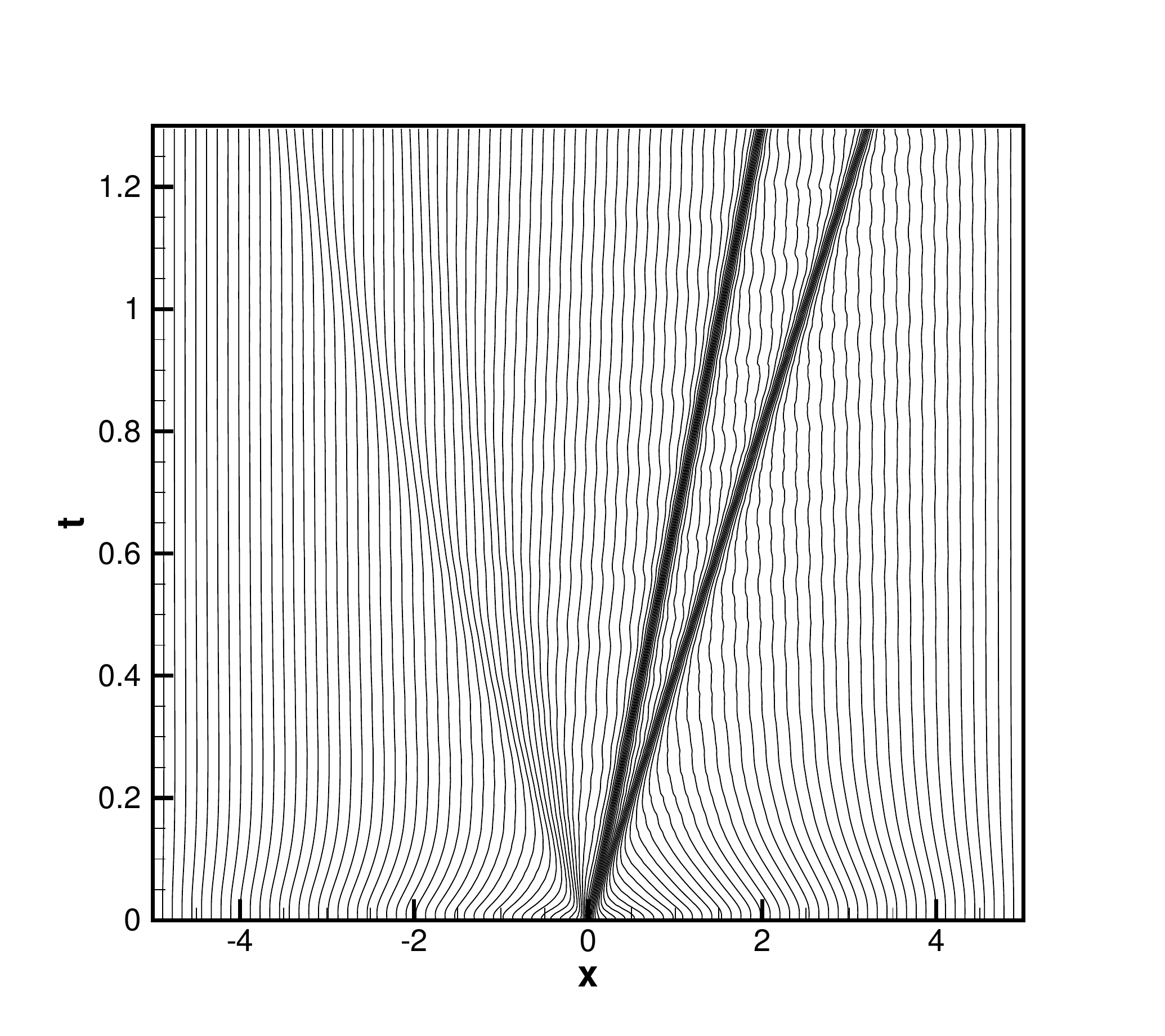}}
 \subfigure[$P^2$]
 {\includegraphics[width=8cm]{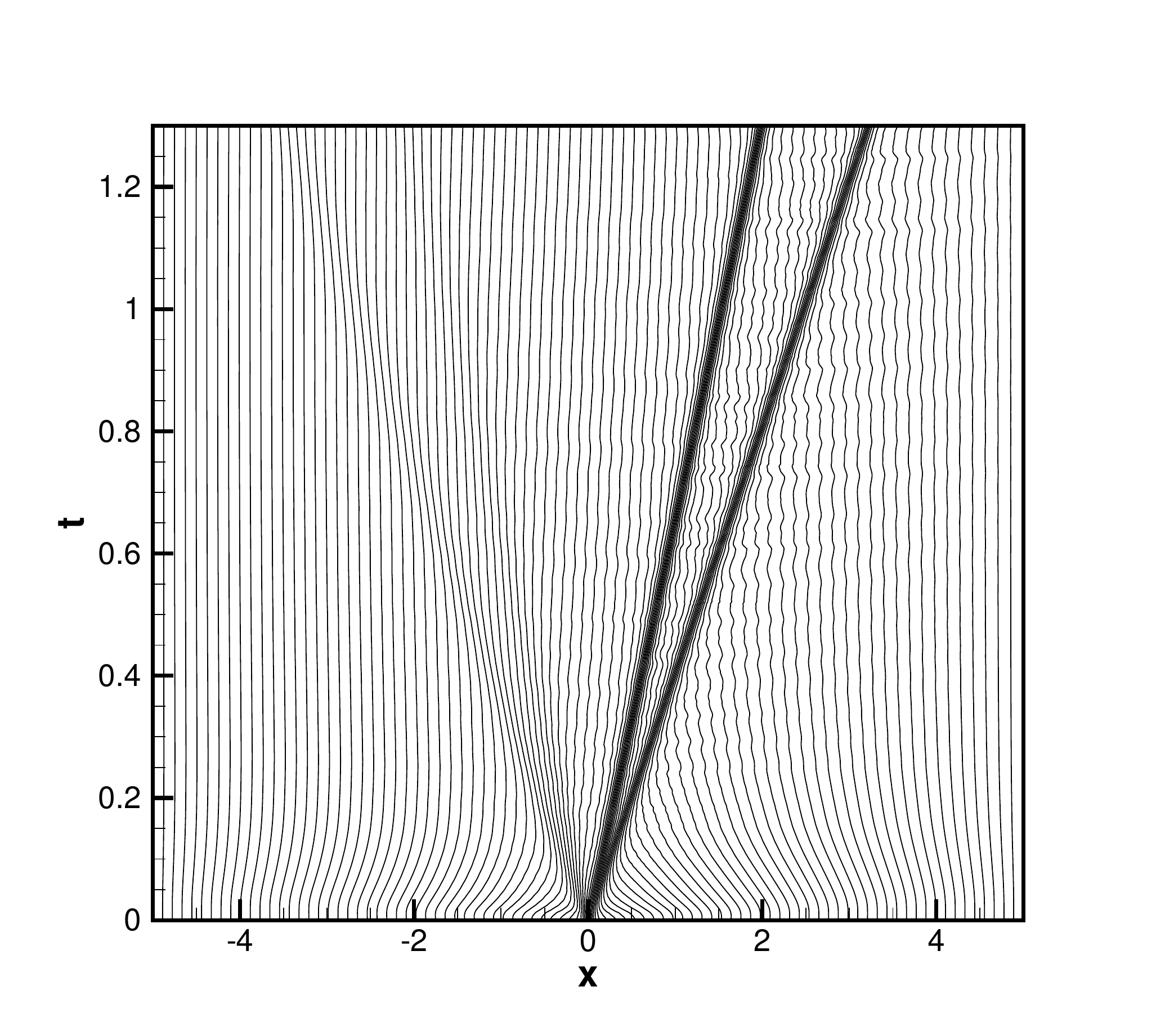}}
 }
 \caption{Example~\ref{exam4.4} (Lax Problem). The trajectories of a moving mesh with $N=100$ are plotted.}
 \label{trfig2}
  \end{center}
   \end{figure}

\begin{exam}{\em
\label{exam4.4}
In this example we consider the Lax problem of the Euler equations (\ref{Euler}) with the initial condition
\begin{equation}
(\rho,u,p)=
\left
\{
\begin{array}{ll}
(0.445,0.698,3.528), \quad &\text{for}\quad x<0\\
(0.5,0,0.571), \quad &\text{for}\quad x>0\notag
\end{array}
\right.
\end{equation}
and the inflow/outflow boundary condition. The computational domain is $(-5,5)$ and the integration is stopped at $T=1.3$. 
The computed density is plotted in Figs. \ref{fig:edge3} and \ref{fig:edge4} and the trajectories of the moving mesh are shown in Fig. \ref{trfig2}.

From the figures, one can see that the moving mesh solution (density) obtained with $N=100$ is more accurate than that obtained with a uniform mesh with $N=400$ for both $P^1$ and $P^2$ elements in this example. From the trajectories, one can see that the points are concentrated at the shock, contact discontinuity and the rarefaction the same as for the Sod problem.

   
}\end{exam}

\begin{figure}[hbtp]
 \begin{center}
 \mbox{\subfigure[MM: $N=150$, UM: $N=150$]
 {\includegraphics[width=8cm]{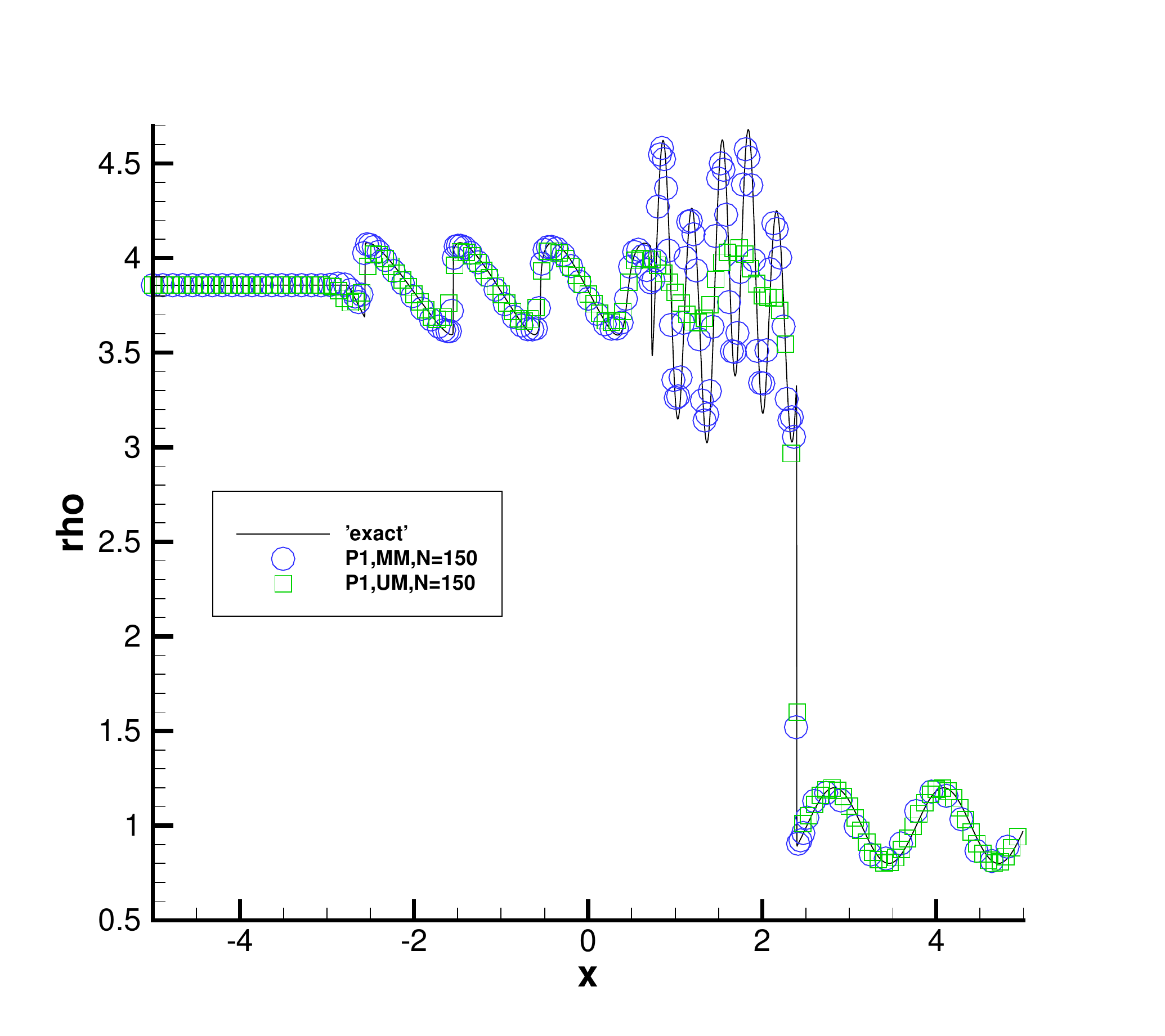}}\quad
   \subfigure[close view of (a) ]
   {\includegraphics[width=8cm]{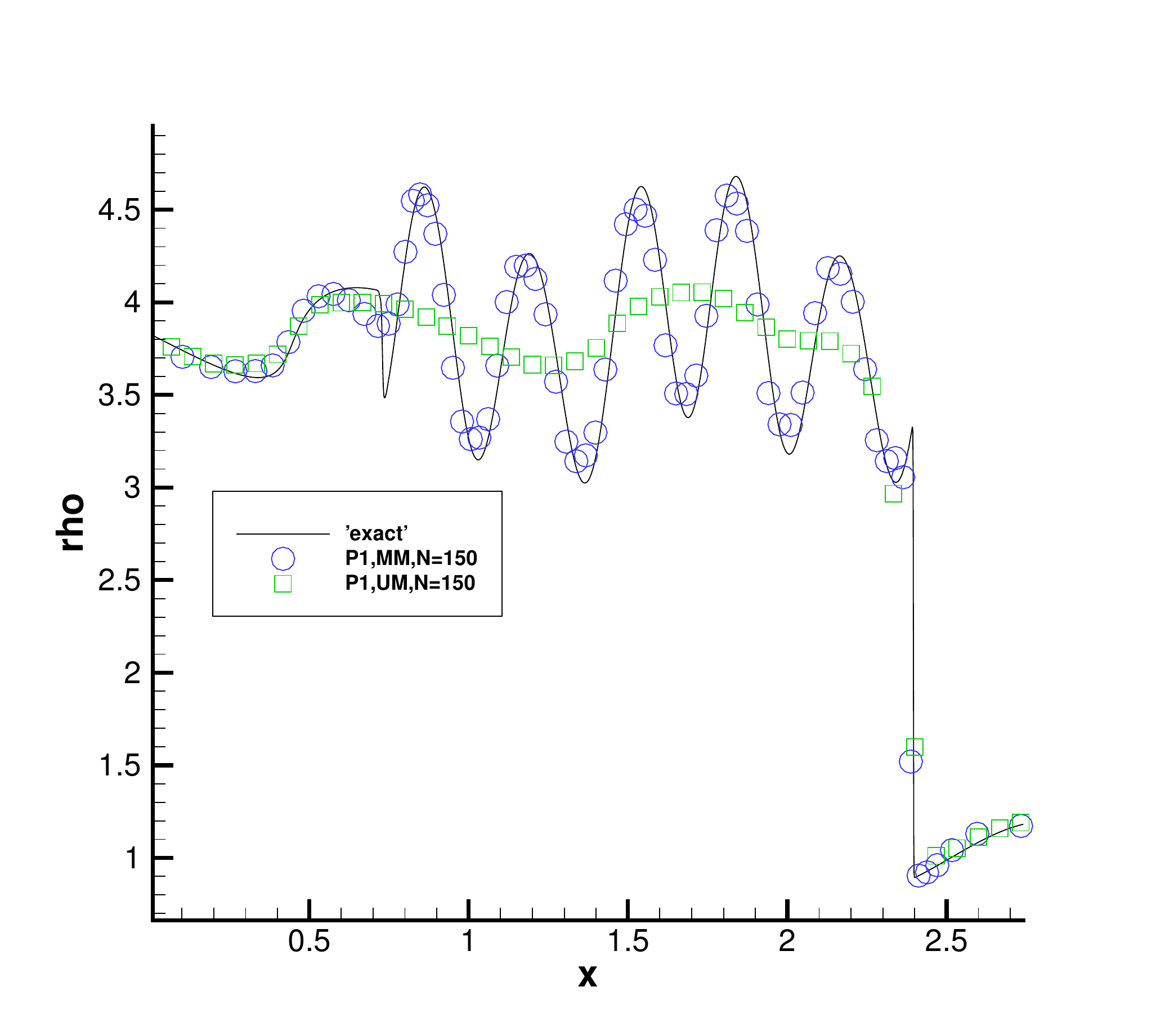}}
   }
 \mbox{\subfigure[MM: $N=150$, UM: $N=400$]
 {\includegraphics[width=8cm]{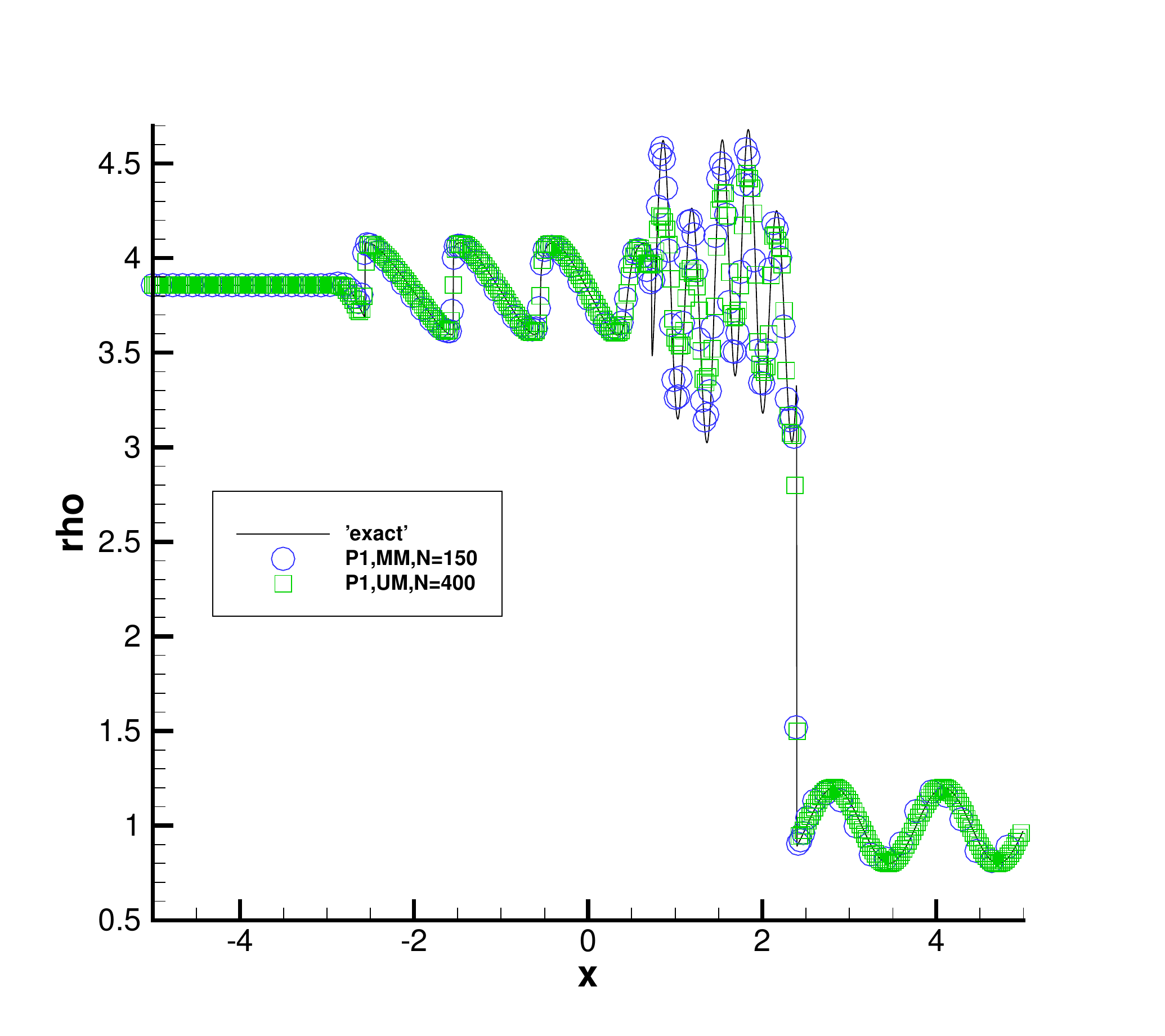}}\quad
   \subfigure[close view of (c) ]
   {\includegraphics[width=8cm]{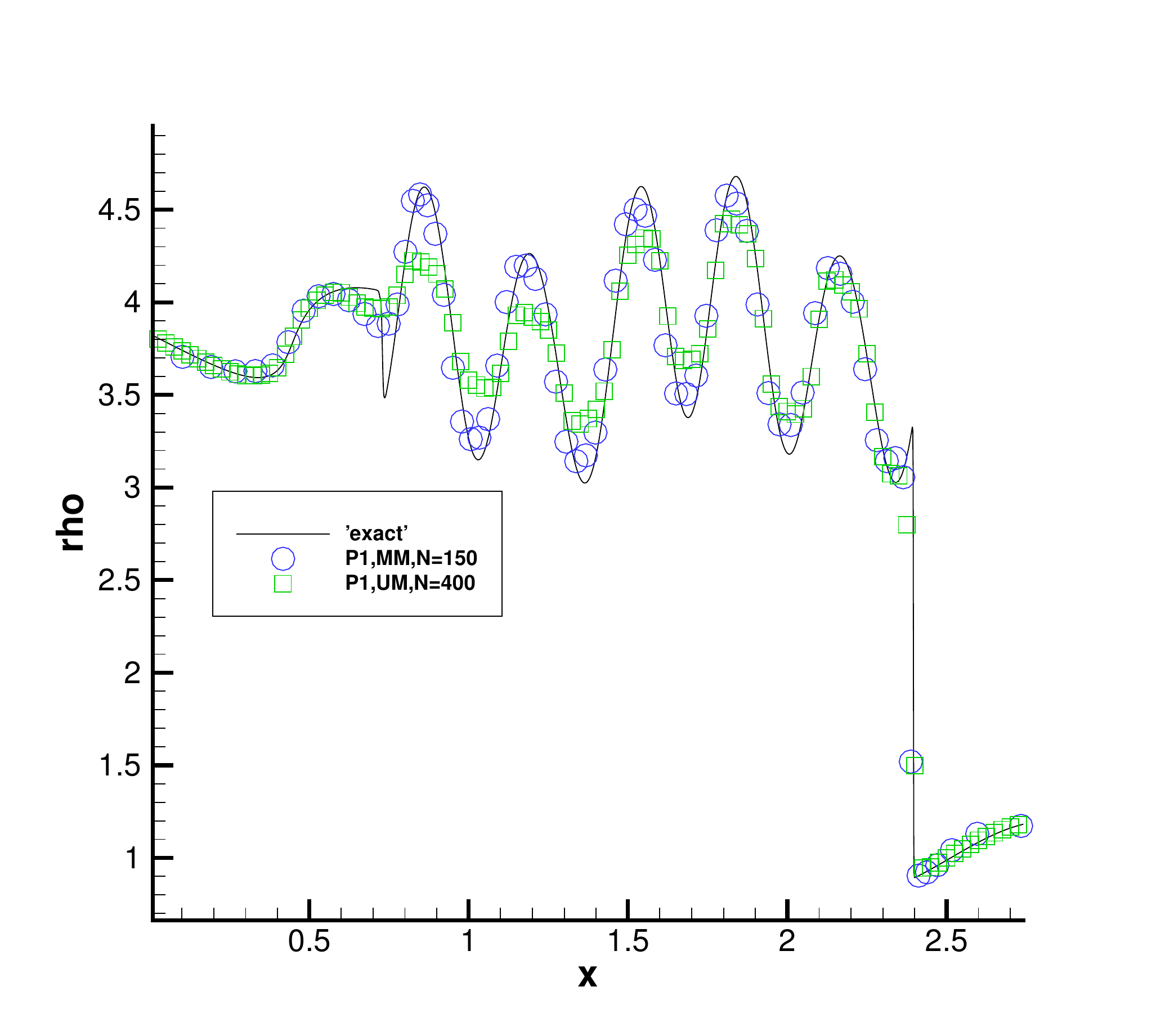}}
  }
   \mbox{\subfigure[MM: $N=150$, UM: $N=600$]
 {\includegraphics[width=8cm]{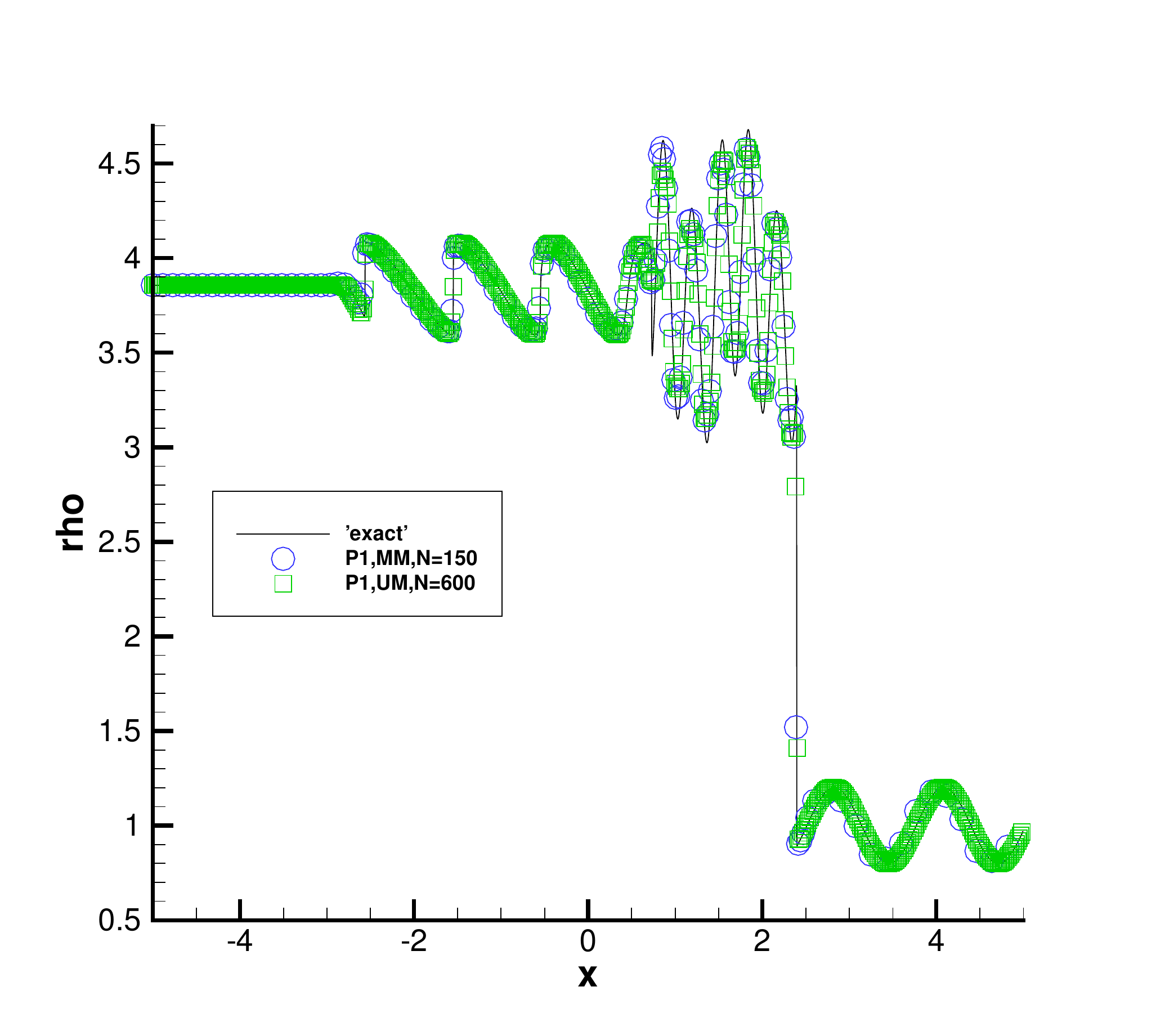}}\quad
   \subfigure[close view of (e)]
   {\includegraphics[width=8cm]{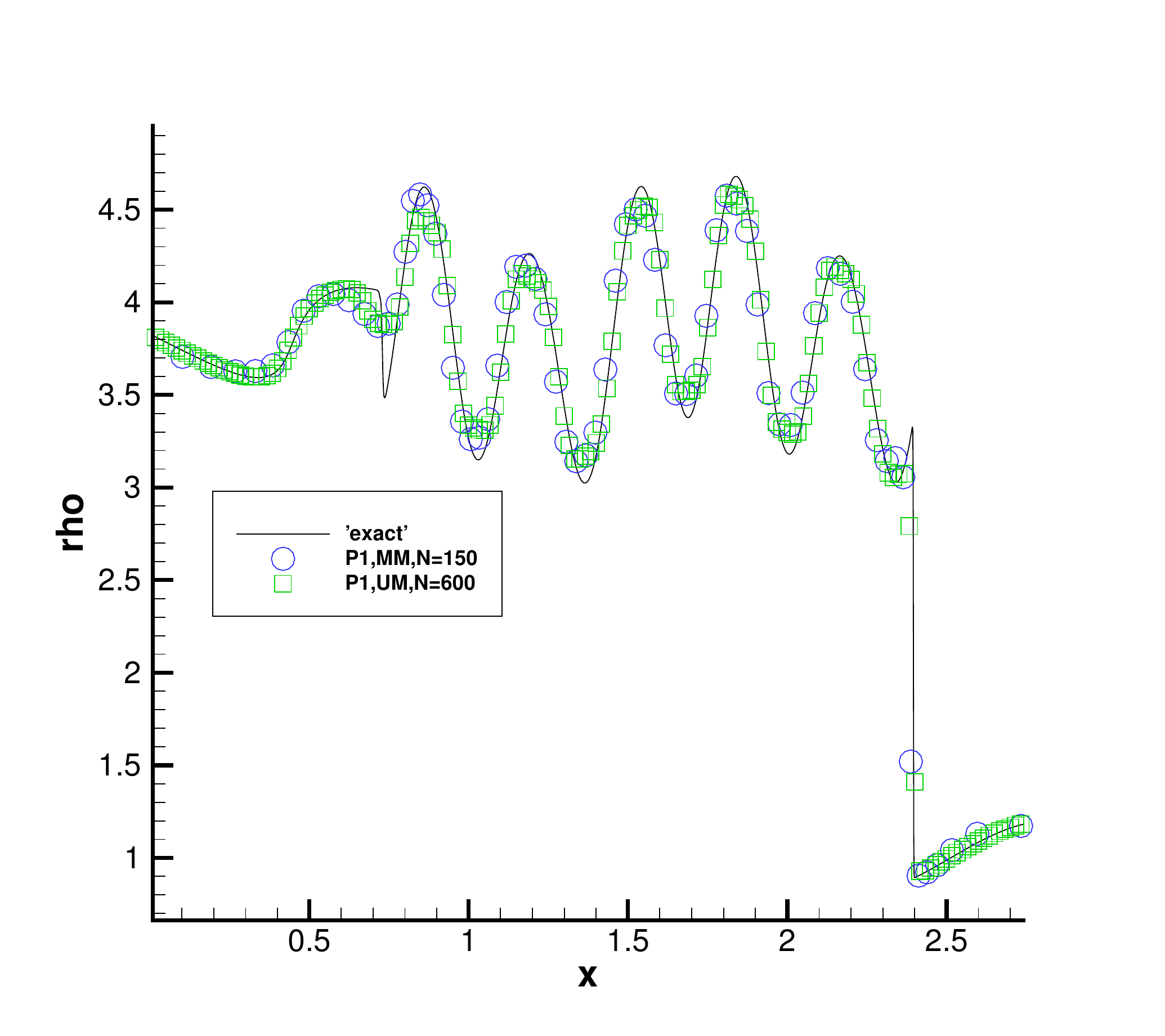}}
   }

   \caption{Example~\ref{exam4.5} (Shu-Osher Problem). The moving mesh solution (density) with $N=150$ is compared with the uniform mesh solutions with $ N=150$, $400$, and $600$. $P^1$ elements are used.}
   \label{fig:edge5}
   \end{center}
   \end{figure}

\begin{figure}[hbtp]
 \begin{center}
 \mbox{\subfigure[MM: $N=150$, UM: $N=150$]
 {\includegraphics[width=8cm]{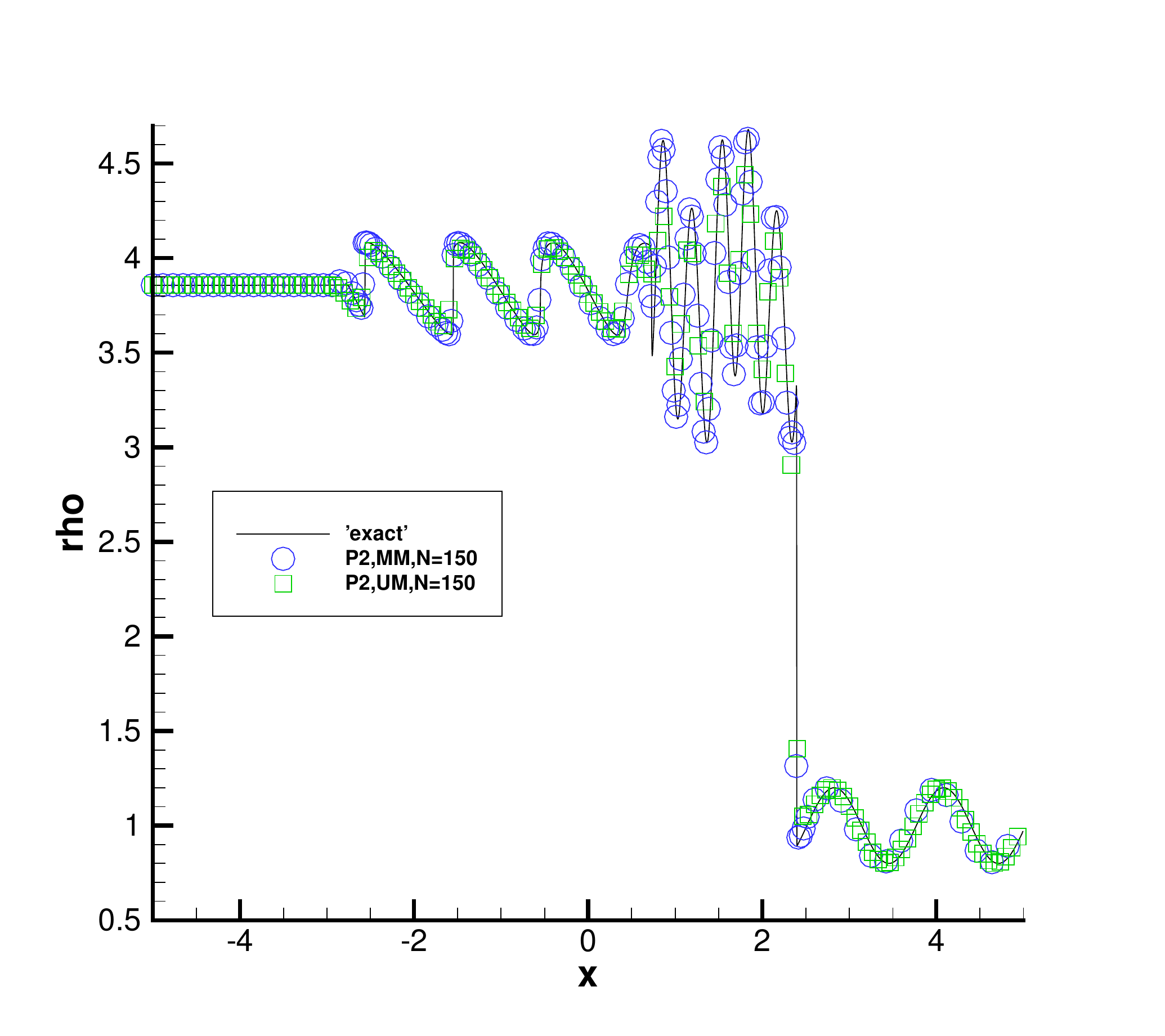}}\quad
   \subfigure[close view of (a) ]
   {\includegraphics[width=8cm]{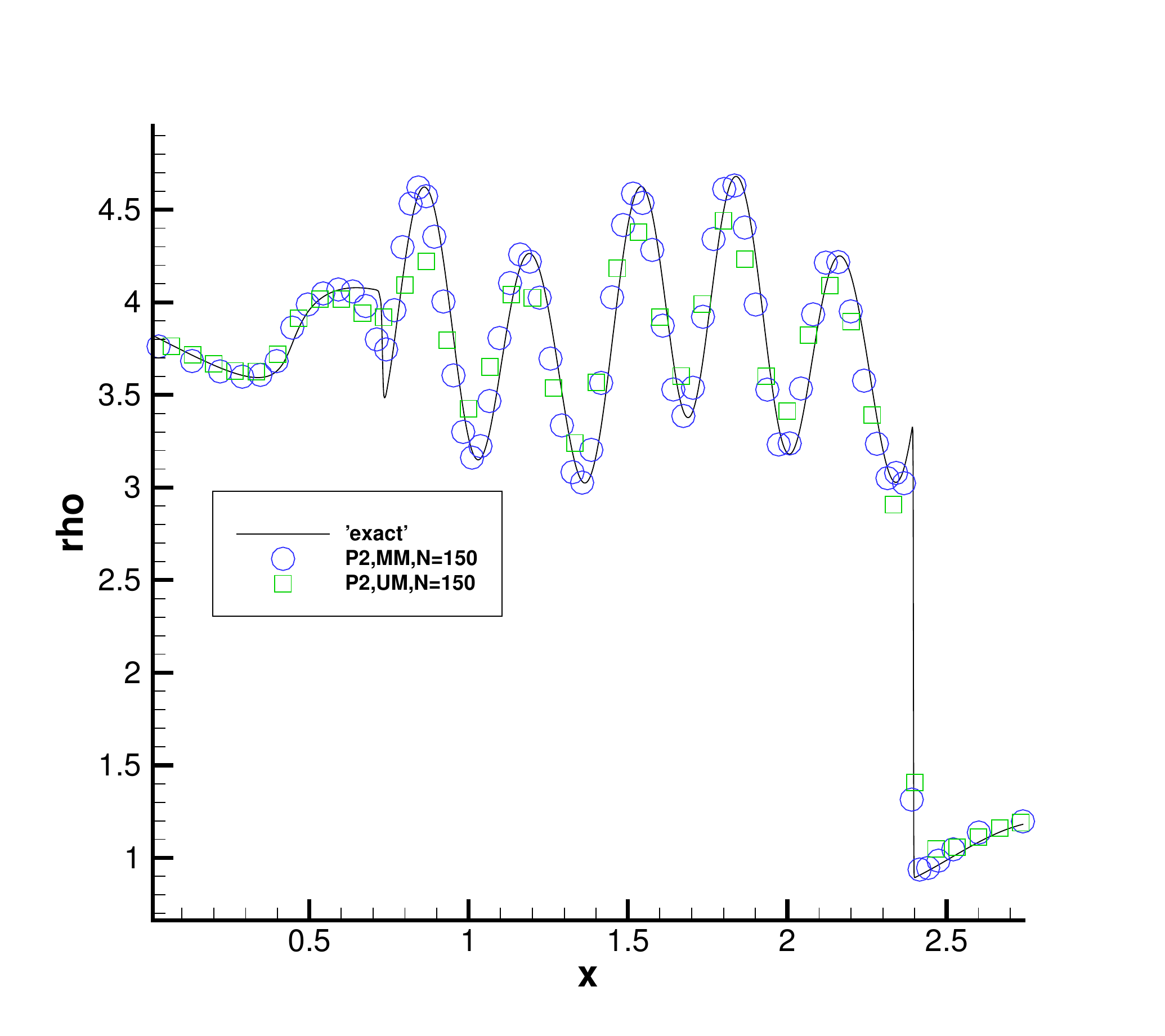}}
   }
 \mbox{\subfigure[MM: $N=150$, UM: $N=400$]
 {\includegraphics[width=8cm]{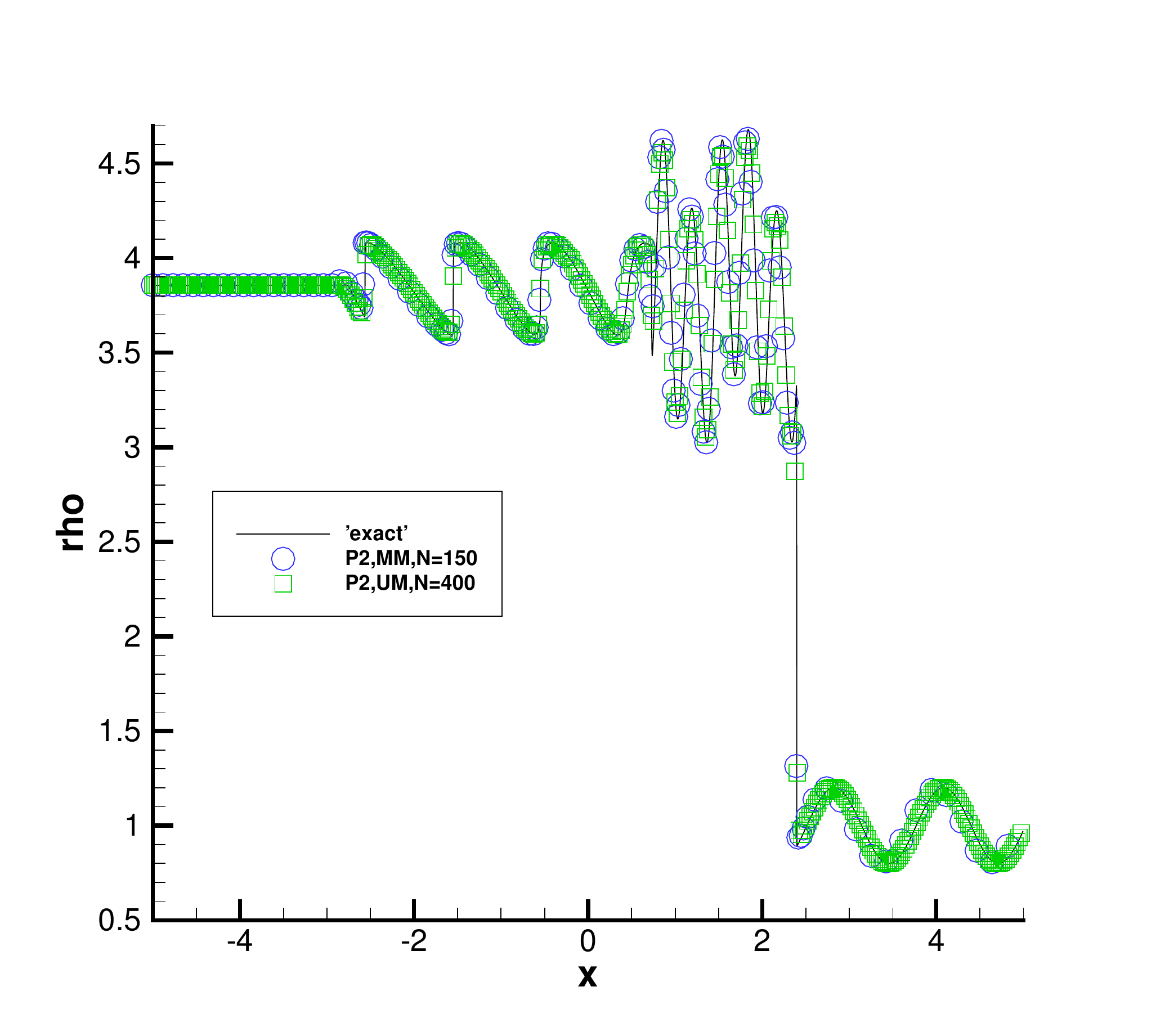}}\quad
   \subfigure[close view of (c) ]
   {\includegraphics[width=8cm]{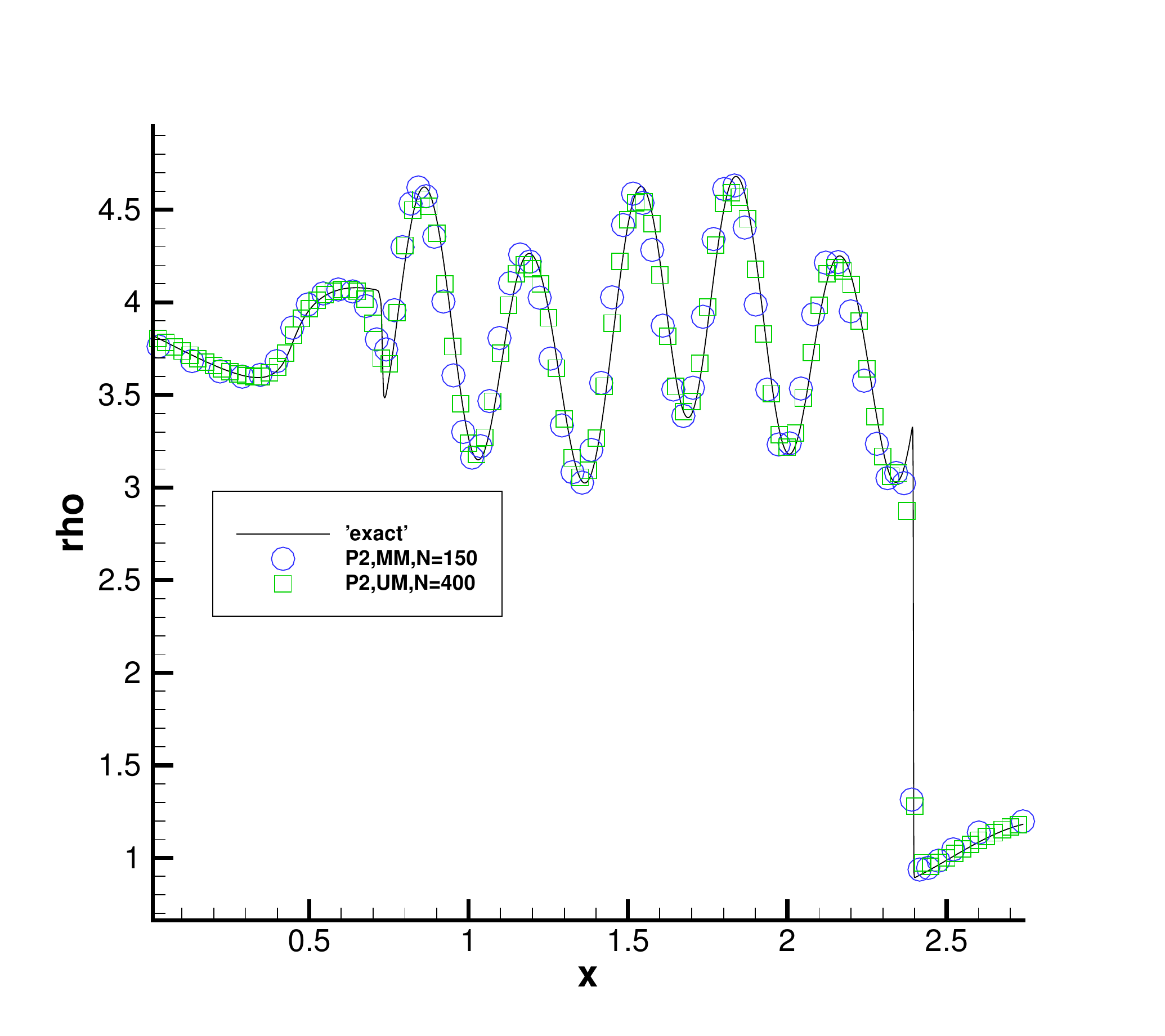}}
   }
   \mbox{\subfigure[MM: $N=150$, UM: $N=600$]
 {\includegraphics[width=8cm]{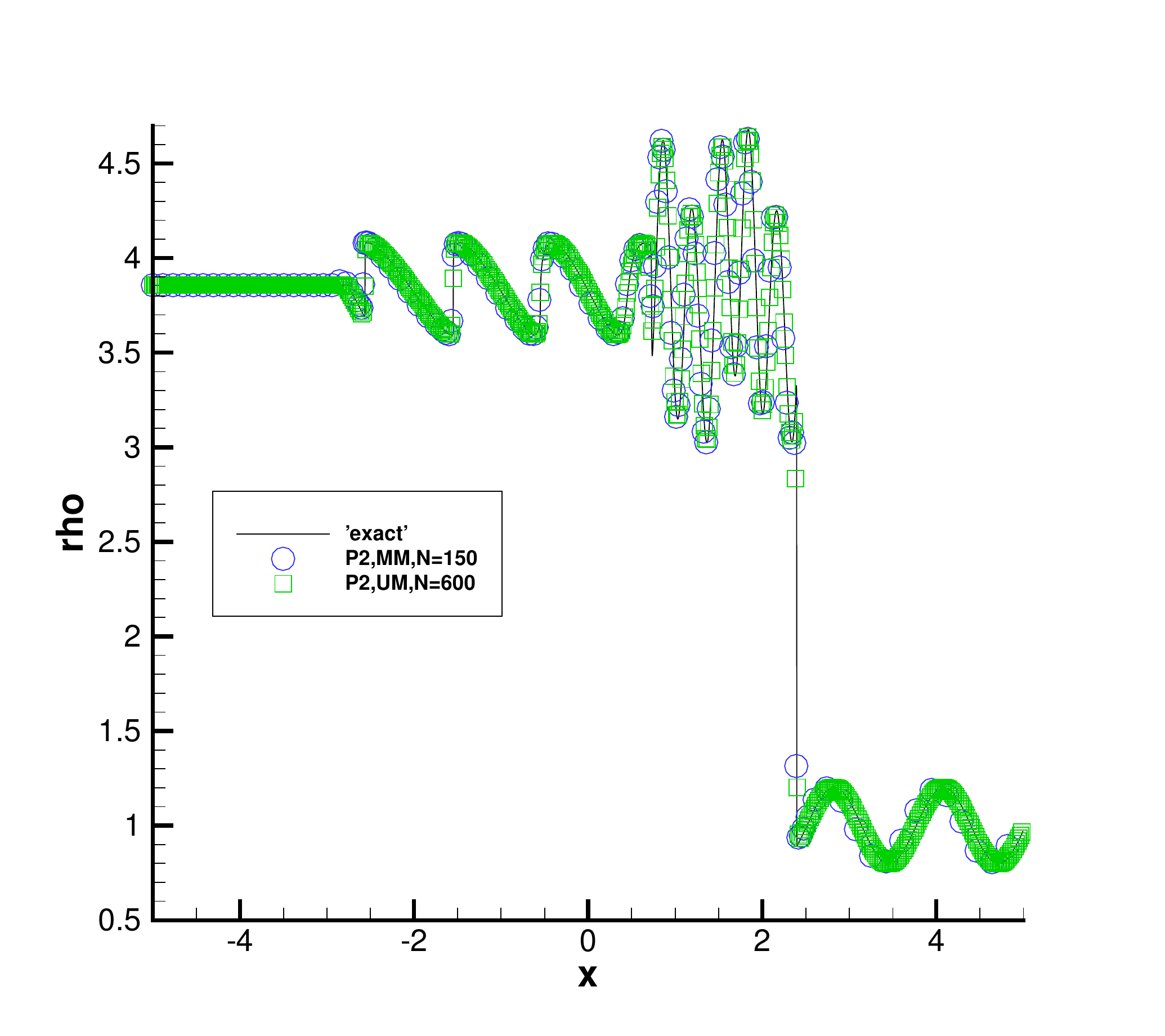}}\quad
   \subfigure[close view of (e) ]
   {\includegraphics[width=8cm]{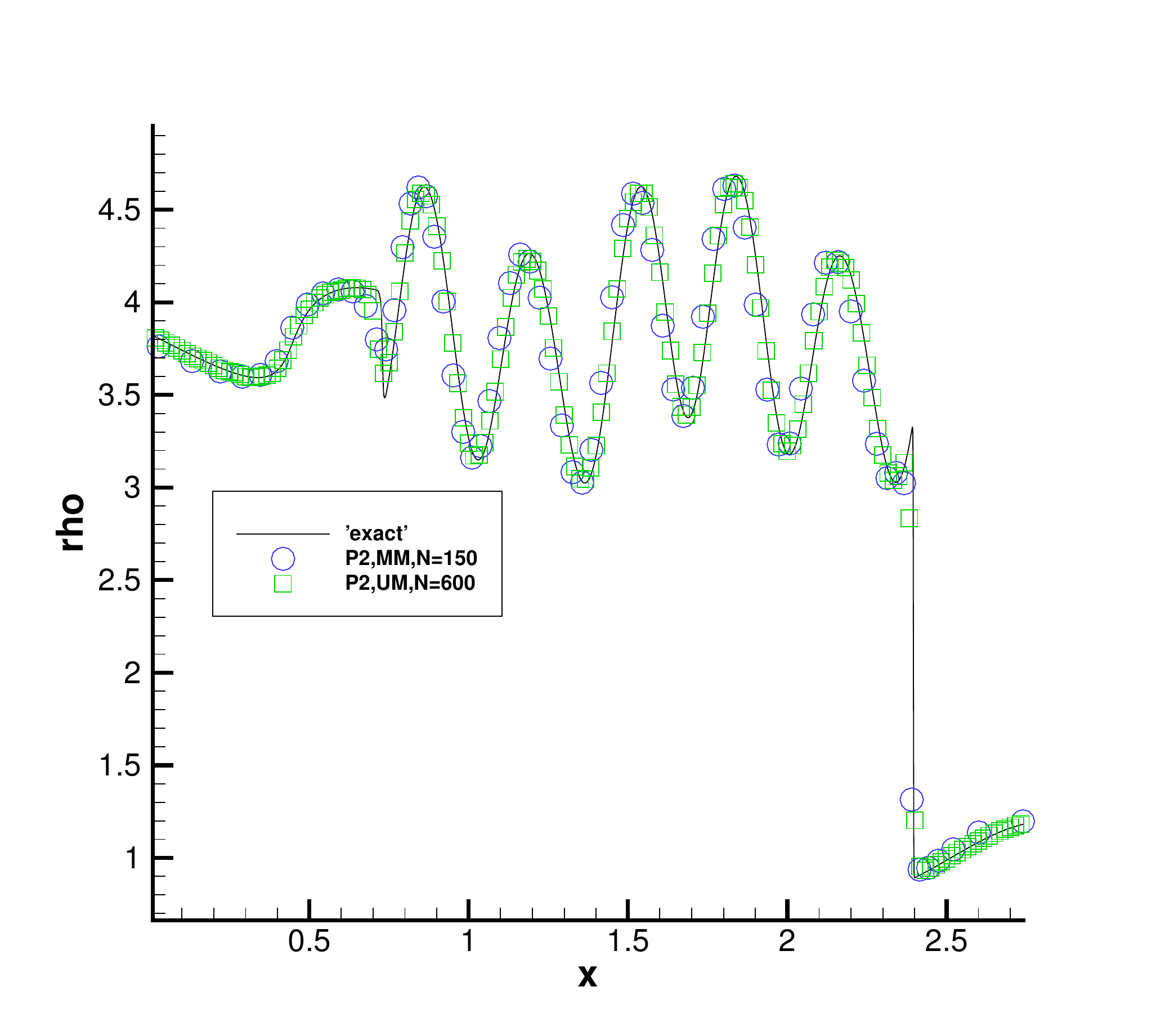}}
   }

   \caption{Example~\ref{exam4.5} (Shu-Osher Problem). The moving mesh solution (density) with $N=150$ is compared with the uniform mesh solutions with $ N=150$, $400$, and $600$. $P^2$ elements are used.}
   \label{fig:edge6}
   \end{center}
   \end{figure}

\begin{figure}[hbtp]
 \begin{center}
 \mbox{\subfigure[$P^1$ elements]
 {\includegraphics[width=8cm]{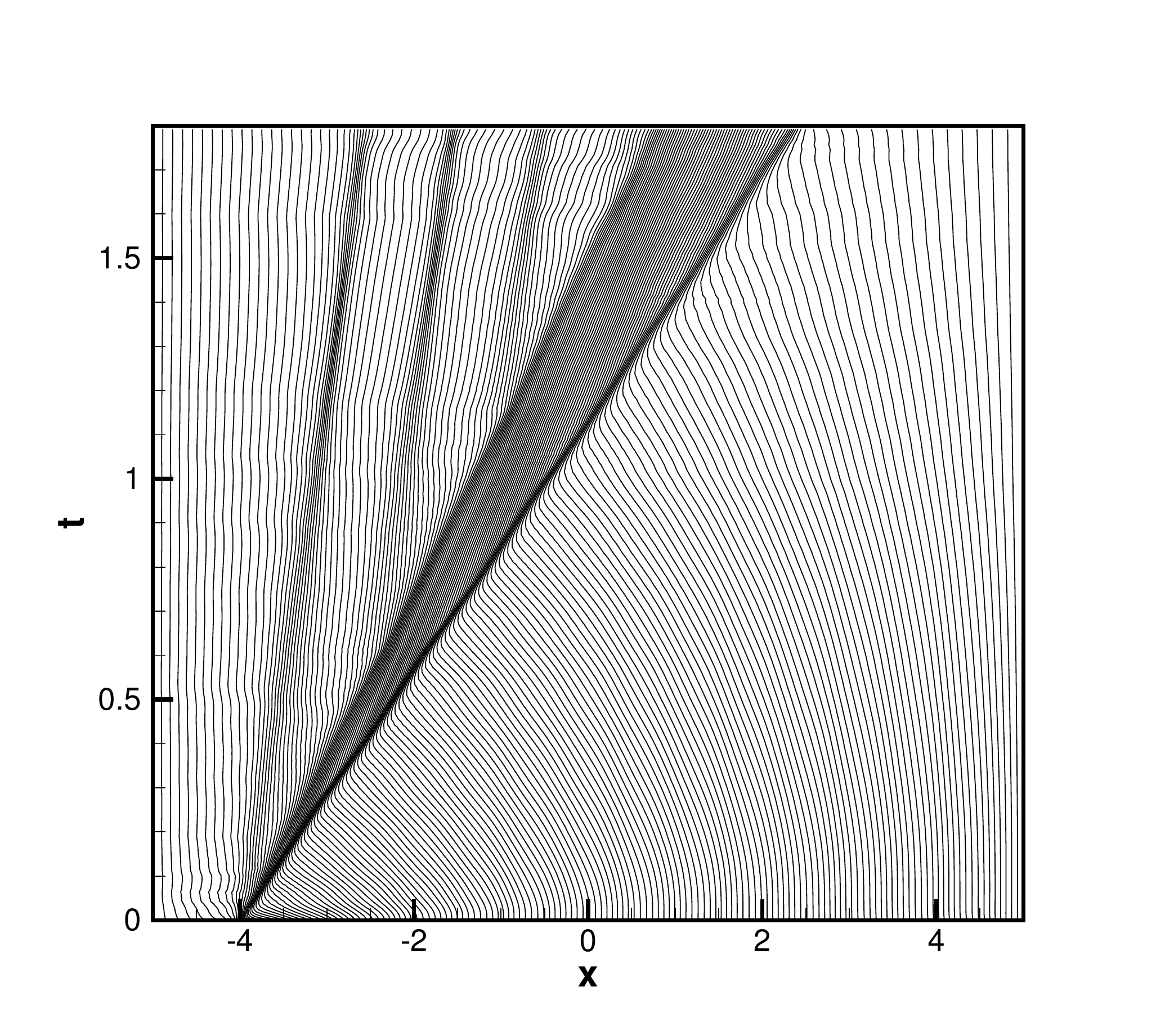}}
 \subfigure[$P^2$ elements]
 {\includegraphics[width=8cm]{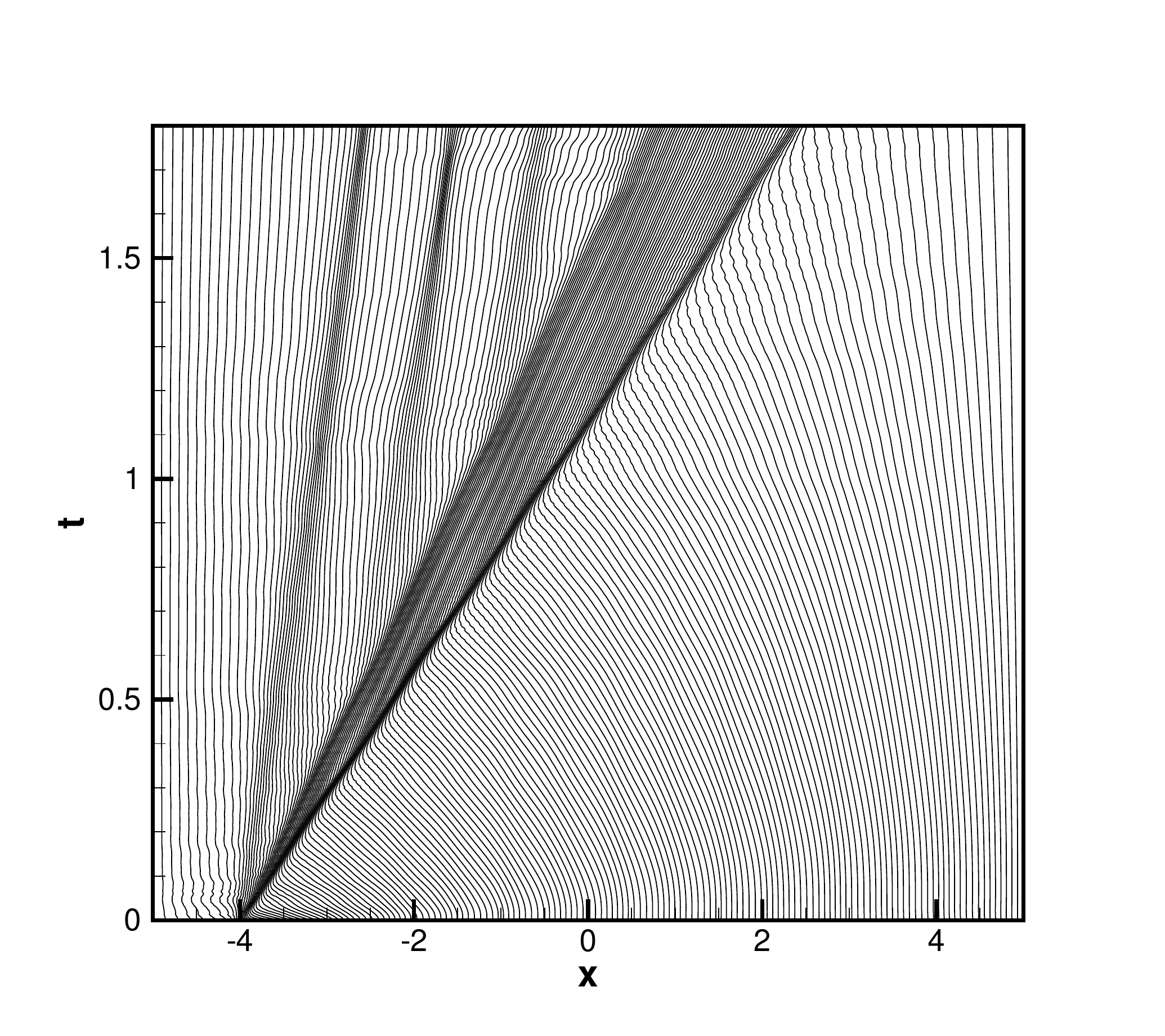}}
 }

 \caption{Example~\ref{exam4.5} (Shu-Osher Problem). The trajectories of a moving mesh with $N=150$ arte plotted.}
 \label{trfig3}
   \end{center}
   \end{figure}

\begin{exam}{\em
\label{exam4.5}
The Shu-Osher problem \cite{EH25} is considered in this example, which contains both shocks and complex smooth region structures. We solve the Euler equations (\ref{Euler}) with a moving shock ($\hbox{Mach}=3$) interacting with  a sine wave in density. The initial condition is
\begin{equation}
(\rho,u,p)=
\begin{cases}
(3.857143,2.629369,10.333333), \quad &\text{for}\quad x<-4\\
(1+0.2\hbox{sin}(5x),0,1), \quad &\text{for}\quad x>-4.\notag
\end{cases}
\end{equation}
The physical domain is taken as $(-5,5)$ in the computation. The computed density is shown at $T=1.8$ against an ``exact solution'' obtained by a fifth-order finite volume WENO scheme with 10,000 uniform points.

The trajectories are plotted in Fig. \ref{trfig3}. The moving mesh solution (density) with $N=150$ is compared with the uniform mesh solutions with $N=150$, $400$, and $600$ in Figs. \ref{fig:edge5} and \ref{fig:edge6}. From the figures, one can observe that for the same number of mesh points, the moving mesh results are clearly better than the uniform mesh ones. Moreover, the moving mesh solution with $N=150$ is comparable with the uniform mesh solution with $N=400$ when $P^2$ elements are used in Fig. \ref{fig:edge6} and is better than the uniform mesh solution with $N=600$ for $P^1$ elements in Fig. \ref{fig:edge5}. These results demonstrate the advantage of the moving mesh method in improving the computational accuracy.

}\end{exam}

\begin{exam}{\em
\label{exam4.6}
We consider the interaction of blast waves of the Euler equations (\ref{Euler}), which was first used by Woodward and Colella \cite{EH26} as a test problem for various numerical schemes.  The initial condition is given by
\begin{equation}
(\rho,u,p)=
\left
\{
\begin{array}{ll}
(1.0,0,1000), \quad& \text{for}\quad 0\leq x<0.1\\
(1.0,0,0.01), \quad& \text{for}\quad 0.1\leq x<0.9\\
(1.0,0,100),\quad&  \text{for}\quad 0.9\leq x\leq 1 .\notag
\end{array}
\right.
\end{equation}
The physical domain is taken as $(0,1)$ and a reflective boundary condition is applied to both ends. The results at time $T=0.038$ are plotted against an ``exact solution" computed by a fifth-order finite difference WENO scheme \cite{EH06} with 81,920 uniform mesh points.

In this example, the parameter $\beta$ in (\ref{ent}) is taken as 1. The trajectories of a moving mesh method are plotted in Fig. \ref{trfig5} which show that the two blast waves propagate in the right direction and finally collide. From Figs. \ref{fig:edge8} and \ref{fig:edge10}, we can see that for the same number of mesh points, the solutions with a moving mesh are much better than those with a uniform mesh. In addition, the moving mesh solution obtained with $N=150$ is comparable with the uniform mesh solution obtained with $N=600$.

\begin{figure}[hbtp]
 \begin{center}
 \mbox{\subfigure[MM: $N=150$, UM: $N=150$]
 {\includegraphics[width=8cm]{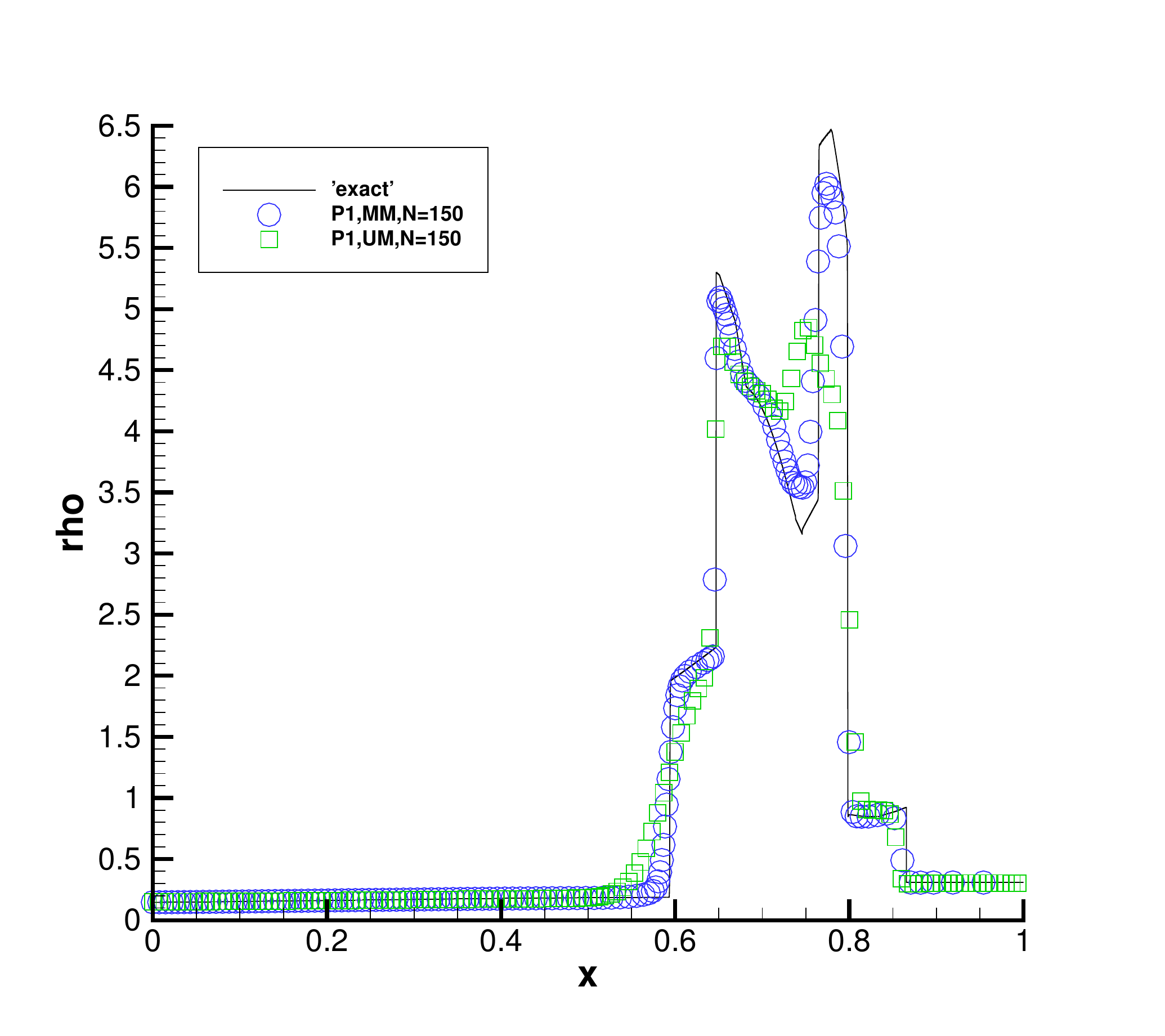}}\quad
   \subfigure[close view of (a) near shock]
   {\includegraphics[width=8cm]{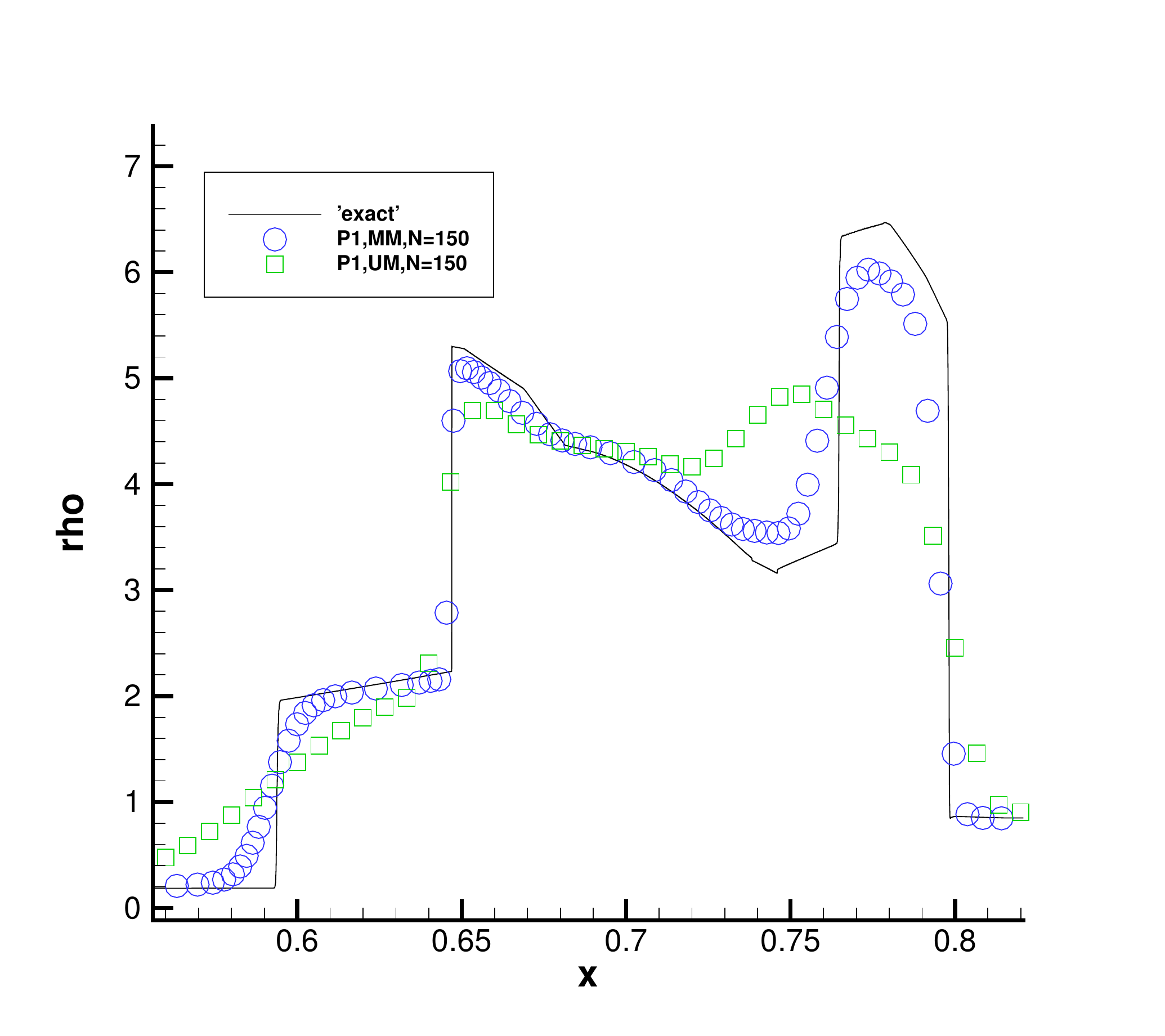}}
   }
 \mbox{\subfigure[MM: $N=150$, UM: $N=600$]
 {\includegraphics[width=8cm]{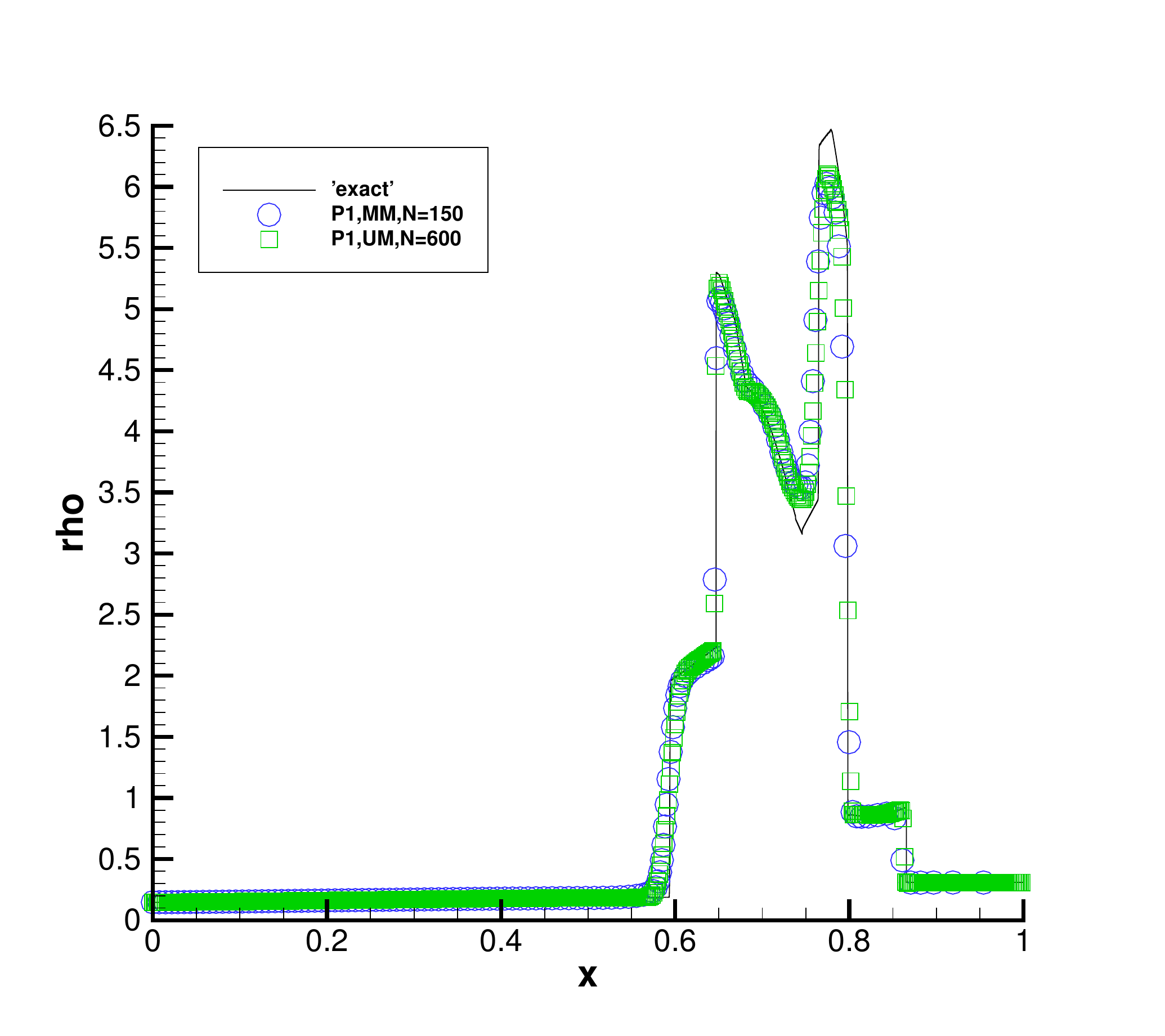}}\quad
   \subfigure[close view of (c) near shock]
   {\includegraphics[width=8cm]{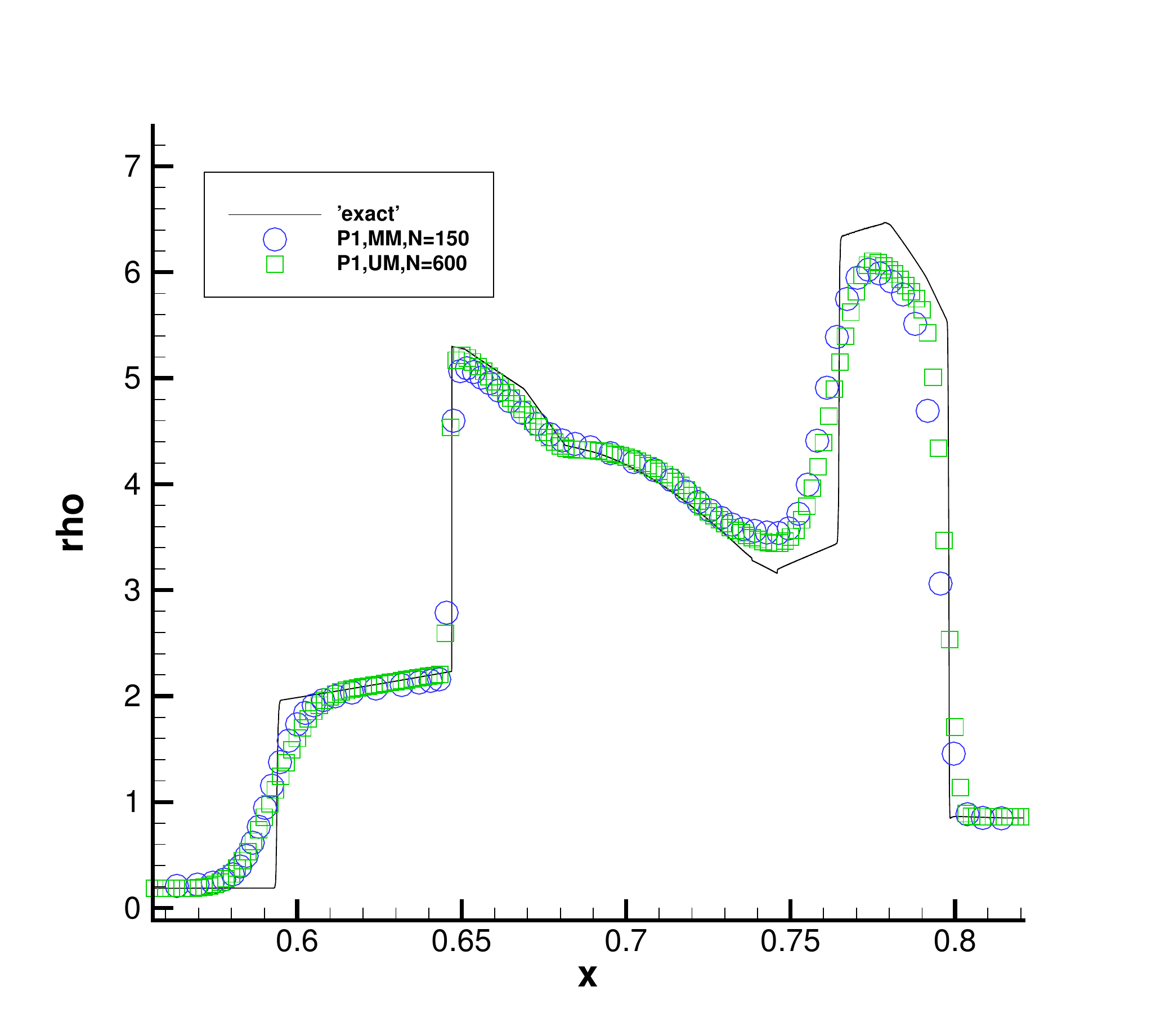}}
   }

   \caption{Example~\ref{exam4.6} (Blastwave Problem). The moving mesh solution (density) with $N=150$ is compared with the uniform mesh solutions with  $N=150$ and $600$. $P^1$ elements are used.}
   \label{fig:edge8}
   \end{center}
   \end{figure}

   \begin{figure}[hbtp]
 \begin{center}
 \mbox{\subfigure[MM: $N=150$, UM: $N=150$]
 {\includegraphics[width=8cm]{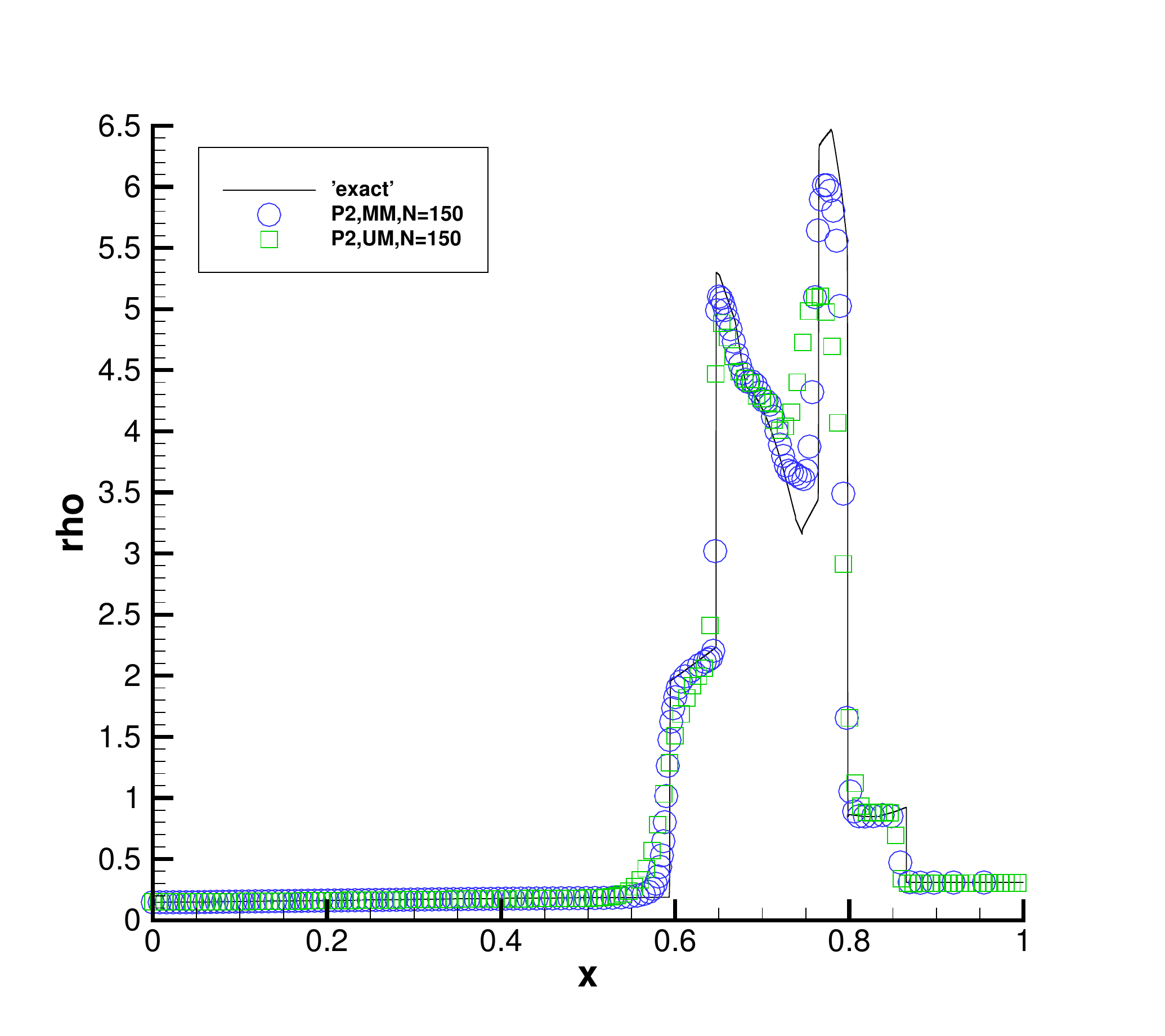}}\quad
   \subfigure[close view of (a) near shock]
   {\includegraphics[width=8cm]{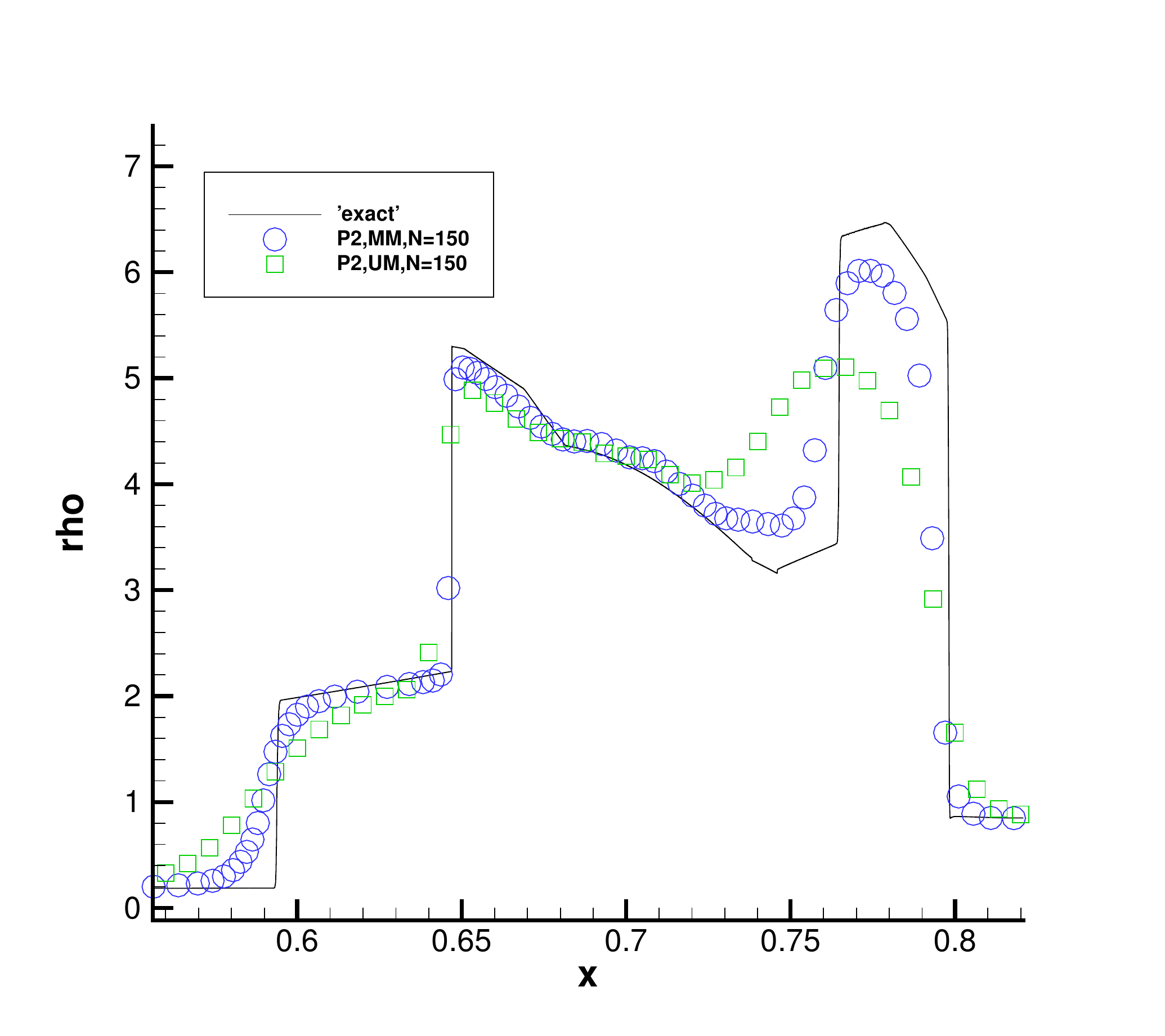}}
   }
 \mbox{\subfigure[MM: $N=150$, UM: $N=600$]
 {\includegraphics[width=8cm]{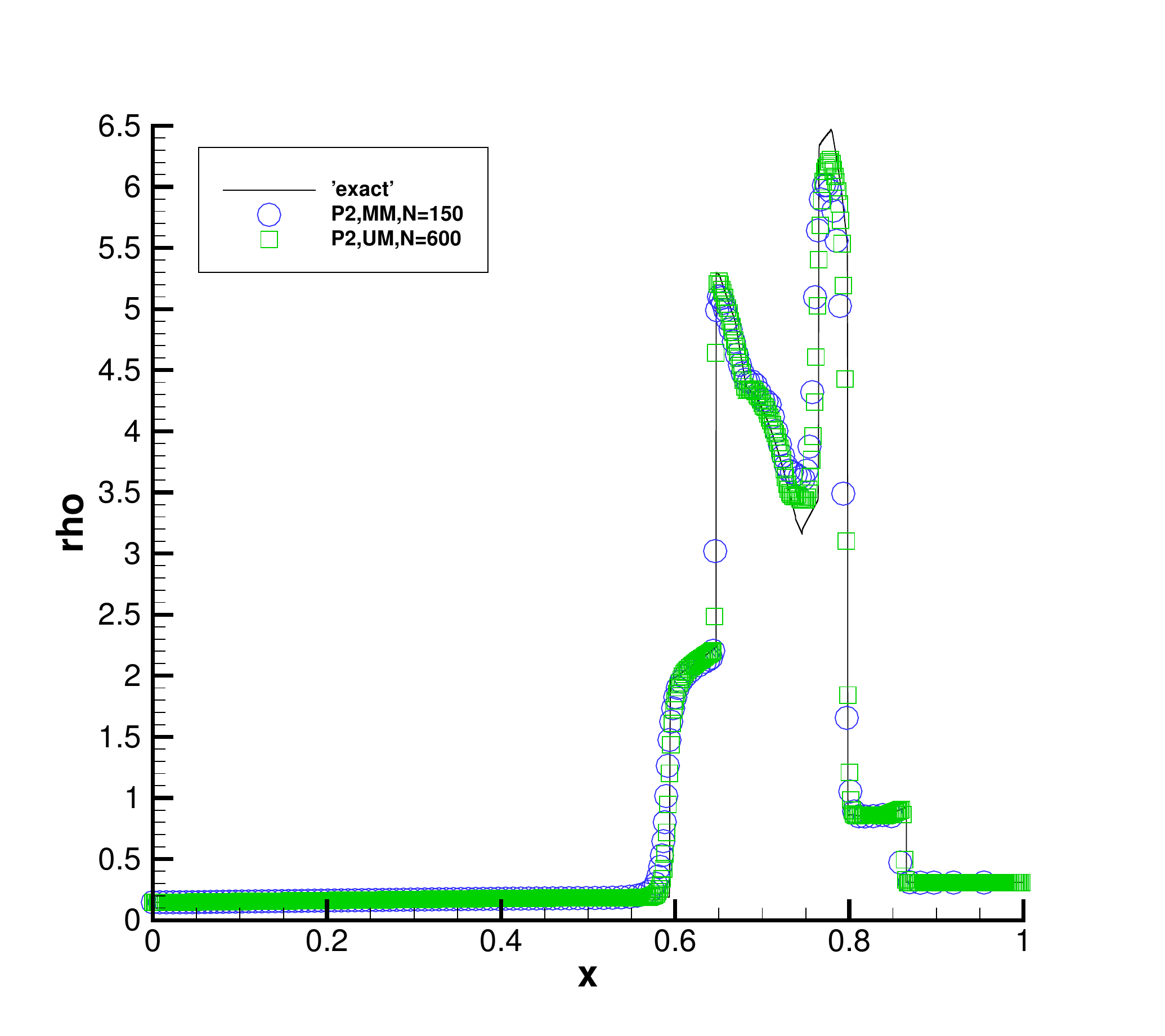}}\quad
   \subfigure[close view of (c) near shock]
   {\includegraphics[width=8cm]{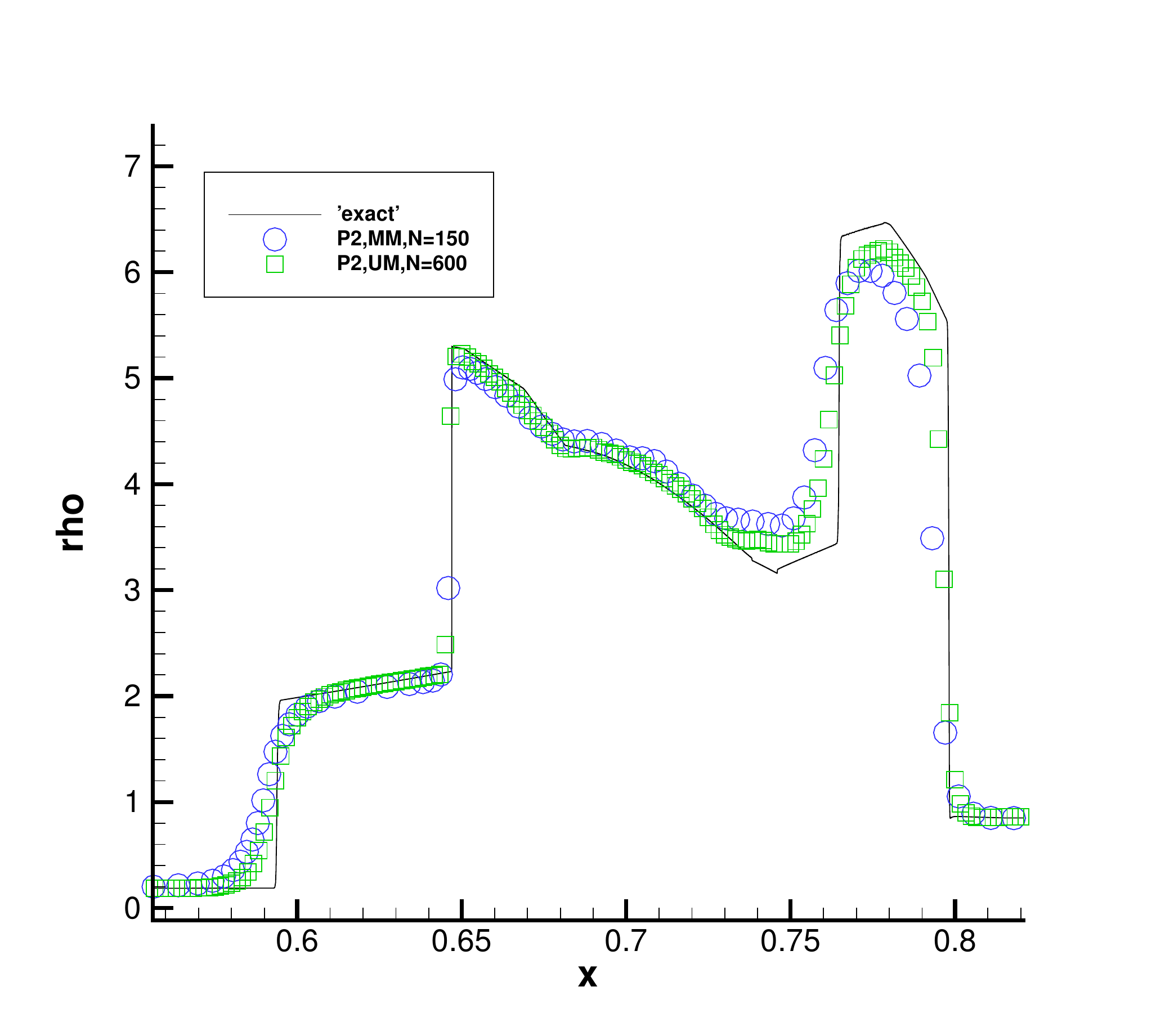}}
   }

   \caption{Example~\ref{exam4.6} (Blastwave Problem). The moving mesh solution (density) with $N=150$ is compared with the uniform mesh solutions with  $N=150$ and $600$. $P^2$ elements are used.}
   \label{fig:edge10}
   \end{center}
   \end{figure}

   \begin{figure}[hbtp]
 \begin{center}
 \mbox{\subfigure[$P^1$ elements]
 {\includegraphics[width=8cm]{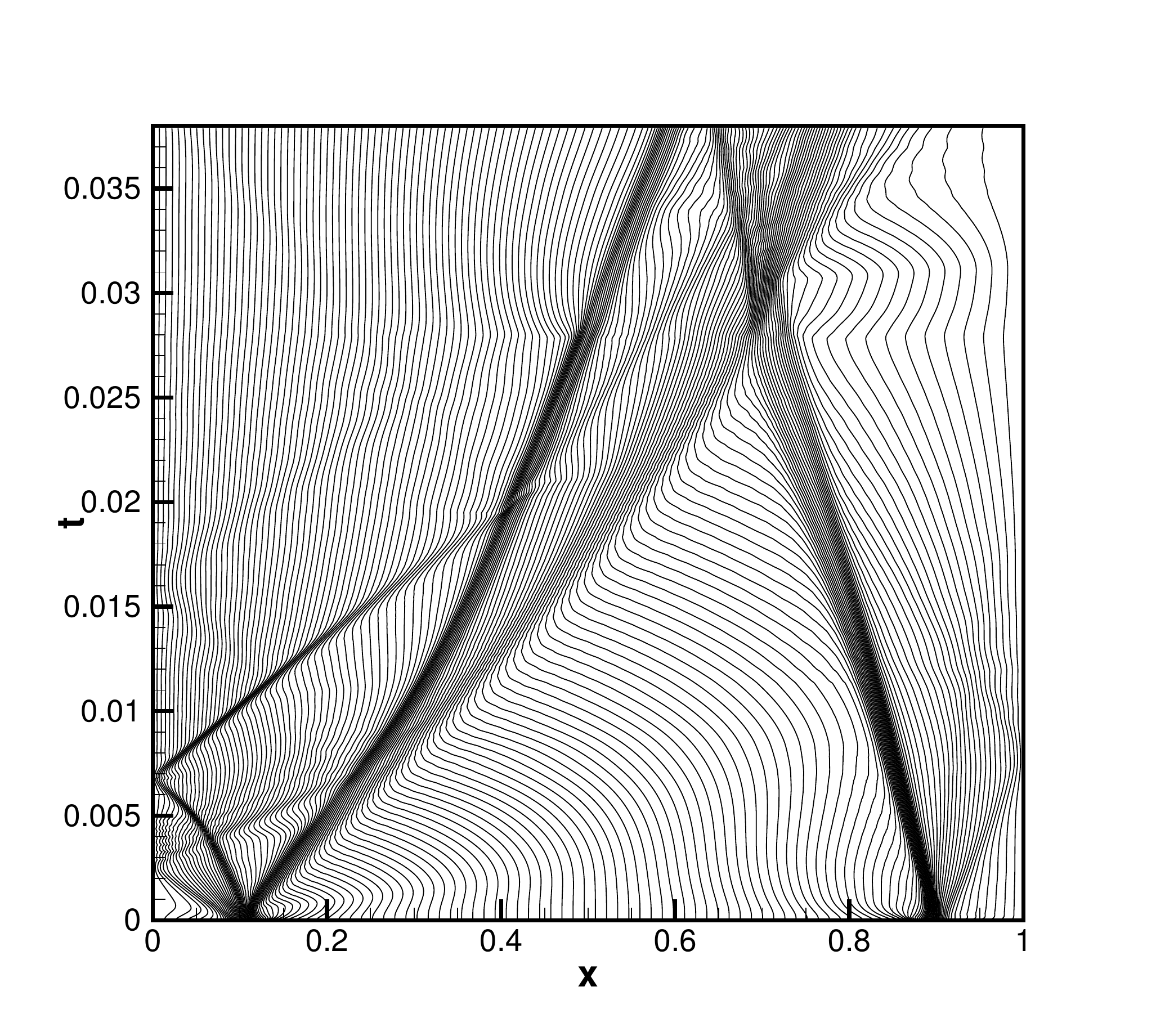}}
 \subfigure[$P^2$ elements]
 {\includegraphics[width=8cm]{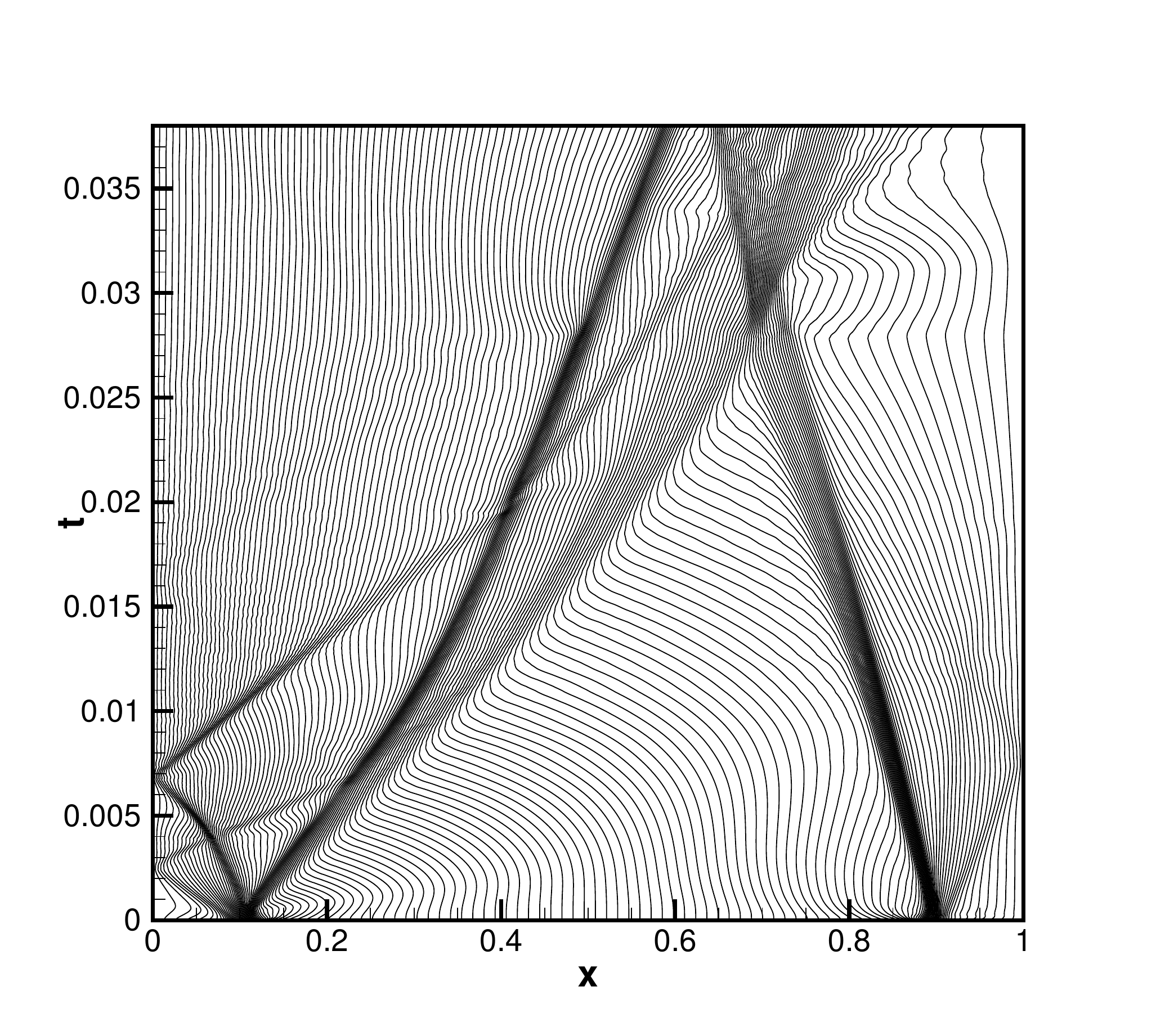}}
 }

 \caption{Example~\ref{exam4.6}  (Blastwave Problem). The trajectories of a moving mesh with $N=150$ are plotted.}
 \label{trfig5}
  \end{center}
   \end{figure}

}\end{exam}

\subsection{Two-dimensional examples}

\begin{figure}[hbtp]
  \begin{center}
  {\includegraphics[scale=0.3]{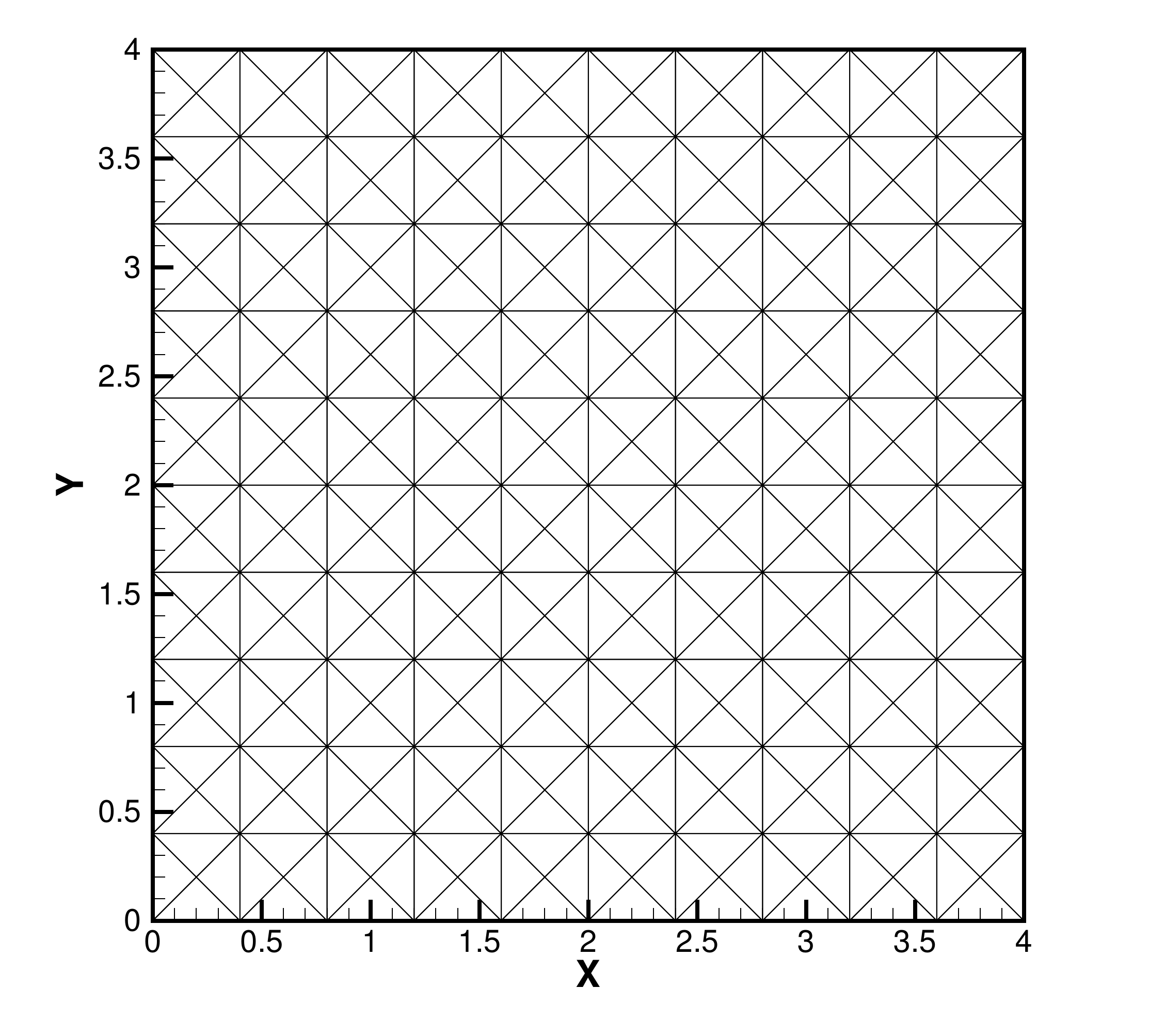}}
    \caption{An initial triangular mesh used in two dimensional computation. A moving mesh associated with
    this initial mesh is denoted by $N = 10\times 10 \times 4$.}
    \label{sample}
    \end{center}
    \end{figure}

For two-dimensional examples, an initial triangular mesh is obtained by dividing any rectangular element into
four triangular elements; see Fig. \ref{sample}. A moving mesh associated with the initial mesh
in Fig. \ref{sample} is denoted by $N = 10\times 10 \times 4$. Other meshes will be denoted similarly.

\begin{exam}{\em
\label{exam4.7}
We solve Burgers' equation in two dimensions,
\[
u_t+\left (\frac{u^2}{2}\right )_x+\left (\frac{u^2}{2}\right )_y=0,\quad (x,y)\in(0,4)\times(0,4)
\]
subject to the initial condition $u(x,y,0)=0.5+\hbox{sin}(\frac{\pi (x+y)}{2})$ and a periodic boundary condition in both directions. 

We compute the solution up to $ T = \frac{0.5}{\pi}$ when the solution is still smooth. The error is listed in Table \ref{ex4.7}. The results show the anticipated $(k+1)$th order convergence of the moving mesh DG method when $P^k$ elements ($k=1$ and $2$) are used.

To see how the smoothness of the mesh affects the accuracy of the method, we list the $L^1$ error
in Tables~\ref{ex4.7nc} for the discontinuous solutions computed with $ T = \frac{1.5}{\pi}$ and different numbers
of sweeps of the low-pass filter applied to the metric tensor. As for Example~\ref{ex4.1}, one can see that the error
and convergence order are comparable although the error is slightly worse for the case with more sweeps.

}\end{exam}

\begin{table}
\caption{Example~\ref{exam4.7} (Burgers' equation): Solution error with periodic boundary conditions and $T=\frac{0.5}{\pi}$.}
\renewcommand{\multirowsetup}{\centering}
\begin{center}
\begin{tabular}{|c|c|c|c|c|c|c|c|c|c|c|c|c|}
\hline
    $k$ & $N$       & {\scriptsize $2\times 2 \times 4$} & {\scriptsize $4\times 4 \times 4$} & {\scriptsize $8\times 8 \times 4$} & {\scriptsize $16\times 16 \times 4$} & {\scriptsize $32\times 32 \times 4$} & {\scriptsize $64\times 64 \times 4$} & {\scriptsize $128\times 128 \times 4$} \\
\hline
\multirow{6}{1cm}{1}
& $L^1$       & 6.419e-1 & 1.743e-1 & 4.231e-2 & 9.996e-3 & 2.575e-3  & 6.557e-4 & 1.639e-4\\
&  Order      & \quad    & 1.88     & 2.04     & 2.08     &  1.96        & 1.97  &2.00 \\
&$L^2$       & 5.714e-1 & 1.711e-1 & 4.361e-2 & 1.054e-2 &  2.735e-3  & 7.041e-4 & 1.781e-4\\
 &Order       & \quad    & 1.74     & 1.97     & 2.05    &   1.95      & 1.96 & 1.98  \\
&$L_{\infty}$&2.661e-1  &1.025e-1  & 3.199e-2 & 9.034e-3  &  2.366e-3 & 6.175e-4 & 1.508e-4\\
&Order    & \quad       & 1.38     & 1.68     & 1.82    &   1.93      & 1.94 & 2.03  \\
\hline
\multirow{6}{1cm}{2}

& $L^1$       & 1.518e-1 & 2.300e-2 &4.185e-3 & 5.490e-4 &  6.936e-5 & 8.397e-6 & 1.014e-6 \\
&  Order      & \quad    & 2.72    & 2.46    & 2.93    & 2.98     & 3.05 & 3.05  \\
&$L^2$        & 1.498e-1 & 2.228e-2 & 5.626e-3 & 8.441e-4 &  1.148e-4 & 1.466e-5 & 1.732e-6\\
&Order       & \quad    & 2.75     & 1.99     & 2.74       & 2.88    & 2.97 & 3.08\\
&$L_{\infty}$&7.936e-2  &1.218e-2  & 6.167e-3 & 1.211e-3  & 1.901e-4 & 2.765e-5 & 3.489e-6 \\
&Order    & \quad       & 2.70     & 0.98    & 2.35      & 2.67   & 2.78 &2.99   \\
\hline

\end{tabular}
\end{center}
\label{ex4.7}
\end{table}

\begin{table}
\caption{Example~\ref{exam4.7}: Solution error in $L^1$ norm with periodic boundary conditions and $T=\frac{1.5}{\pi}$. Various numbers of sweeps of a low-pass filter have been applied to the smoothing of the metric tensor.
}
\renewcommand{\multirowsetup}{\centering}
\begin{center}
\begin{tabular}{|c|c|c|c|c|c|c|c|c|c|c|c|}
\hline
$k$ & Sweeps$\backslash$ $N$     & {\scriptsize $2\times 2 \times 4$} & {\scriptsize $4\times 4 \times 4$} & {\scriptsize $8\times 8 \times 4$} & {\scriptsize $16\times 16 \times 4$} & {\scriptsize $32\times 32 \times 4$} & {\scriptsize $64\times 64 \times 4$} & {\scriptsize $128\times 128 \times 4$} \\\hline
  & 3       & 2.083e0 & 8.779e-1  & 3.224e-1 & 1.144e-1 &4.371e-2 & 1.945e-2 & 9.826e-3\\
 &  Order      & \quad      &  1.25        & 1.45 &1.49 & 1.39 & 1.17 & 0.99\\
\multirow{2}{1cm}{1}
  & 30       & 2.084e0 & 8.702e-1  & 3.496e-1 & 1.302e-1 &5.001e-2 & 2.236e-2 & 1.052e-2\\
 &  Order      & \quad      &  1.26     & 1.32 &1.42 & 1.38 & 1.16 & 1.09\\
   & 100       & 2.084e0 & 8.724e-1  & 3.533e-1 & 1.395e-1 &5.304e-2 & 2.377e-2 & 1.100e-2\\
 &  Order      & \quad      &  1.26     & 1.30 &1.34 & 1.39 & 1.16 & 1.11\\
 \hline
   & 3       & 1.074e0 & 4.423e-1 & 1.570e-1 & 6.979e-2 & 3.097e-2 & 1.628e-2 & 9.115e-3\\
 &  Order      & \quad     & 1.28    & 1.49 & 1.17 & 1.17 & 0.93 & 0.84 \\
 \multirow{2}{1cm}{2}
  & 30       & 1.077e0 & 4.370e-1 & 1.851e-1 & 8.258e-2 & 3.719e-2 & 1.903e-2 & 9.760e-3\\
 &  Order      & \quad    & 1.30    & 1.24 & 1.16 & 1.15 & 0.97 & 0.96 \\
  & 100       & 1.077e0 & 4.378e-1 & 1.907e-1 & 8.870e-2 & 4.003e-2 & 2.031e-2 & 1.020e-2\\
 &  Order      & \quad     & 1.30    & 1.20 & 1.10 & 1.15 & 0.98 & 0.99 \\
 \hline
\end{tabular}
\end{center}
\label{ex4.7nc}
\end{table}

\begin{exam}{\em
\label{exam4.8}
The Euler equations (\ref{2d}) are solved
subject to the initial condition $\rho(x,y,0)=1+0.2\hbox{sin}(\frac{\pi (x+y)}{2}),\; \mu(x,y,0)=0.7,\; \nu(x,y,0)=0.3,\; p(x,y,0)=1$ and a periodic boundary condition in both directions.
The computational domain is $(0,2)\times (0,2)$ and the final time is $T=1$. We take the parameter $\beta$ in \eqref{ent} as 100. The results in Table~\ref{ex4.8} show the convergence order of the second order
for $k=1$ and the third order for $k=2$ for the moving mesh DG method for the Euler system in two dimensions.

}\end{exam}

\begin{table}
\caption{Example~\ref{exam4.8}: Solution error (in density) for periodic boundary conditions and $T=1$, $\beta=100$.}
\renewcommand{\multirowsetup}{\centering}
\begin{center}
\begin{tabular}{|c|c|c|c|c|c|c|c|c|c|c|c|c|}
\hline
    $k$ & $N$       & {\scriptsize $2\times 2 \times 4$} & {\scriptsize $4\times 4 \times 4$} & {\scriptsize $8\times 8 \times 4$} & {\scriptsize $16\times 16 \times 4$} & {\scriptsize $32\times 32 \times 4$} & {\scriptsize $64\times 64 \times 4$} & {\scriptsize $128\times 128 \times 4$} \\ \hline
\multirow{6}{1cm}{1}
& $L^1$       & 2.590e-1 & 5.764e-2 & 1.260e-2 & 3.071e-3 &8.804e-4  & 2.079e-4 & 4.687e-5\\
&  Order      & \quad    & 2.17     & 2.18     & 2.06      &  1.81        & 2.05  &2.15 \\
&$L^2$       & 1.509e-1 & 3.752e-2 & 8.455e-3 & 2.128e-3 & 6.387e-4  & 1.525e-4 & 3.329e-5\\
 &Order       & \quad    & 2.01     & 2.15     & 1.99     &   1.74      & 2.07 & 2.20  \\
&$L_{\infty}$&1.542e-1   &8.699e-2  & 2.224e-2 & 5.584e-3  &  1.934e-3 & 4.537e-4 & 9.286e-5\\
&Order    & \quad       & 0.83     & 1.97     & 1.99     &   1.53      & 2.09 & 2.29   \\
\hline
\multirow{6}{1cm}{2}

& $L^1$       & 5.155e-2 & 1.059e-2 &1.710e-3 & 2.406e-4 &  2.888e-5 & 3.425e-6 & 4.203e-7 \\
&  Order      & \quad    & 2.28    & 2.63     & 2.83     & 3.06     & 3.08 & 3.03  \\
&$L^2$        & 3.495e-2 &7.505e-3 & 1.273e-3 & 1.931e-4 &  2.370e-5 & 2.771e-6 & 3.382e-7\\
&Order       & \quad    & 2.22     & 2.56     & 2.72     & 3.03    & 3.10 & 3.03 \\
&$L_{\infty}$&5.637e-2  &2.073e-2  & 4.382e-3 & 7.621e-4  & 1.054e-4 & 1.301e-5 & 1.553e-6 \\
&Order    & \quad       & 1.44     & 2.24    & 2.52      & 2.85   & 3.02 &3.07   \\
\hline
\end{tabular}
\end{center}
\label{ex4.8}
\end{table}

\begin{figure}[hbtp]
  \begin{center}
  \mbox{\subfigure[Fixed Mesh, $N = 120 \times 30 \times 4$] 
  {\includegraphics[width=8cm]{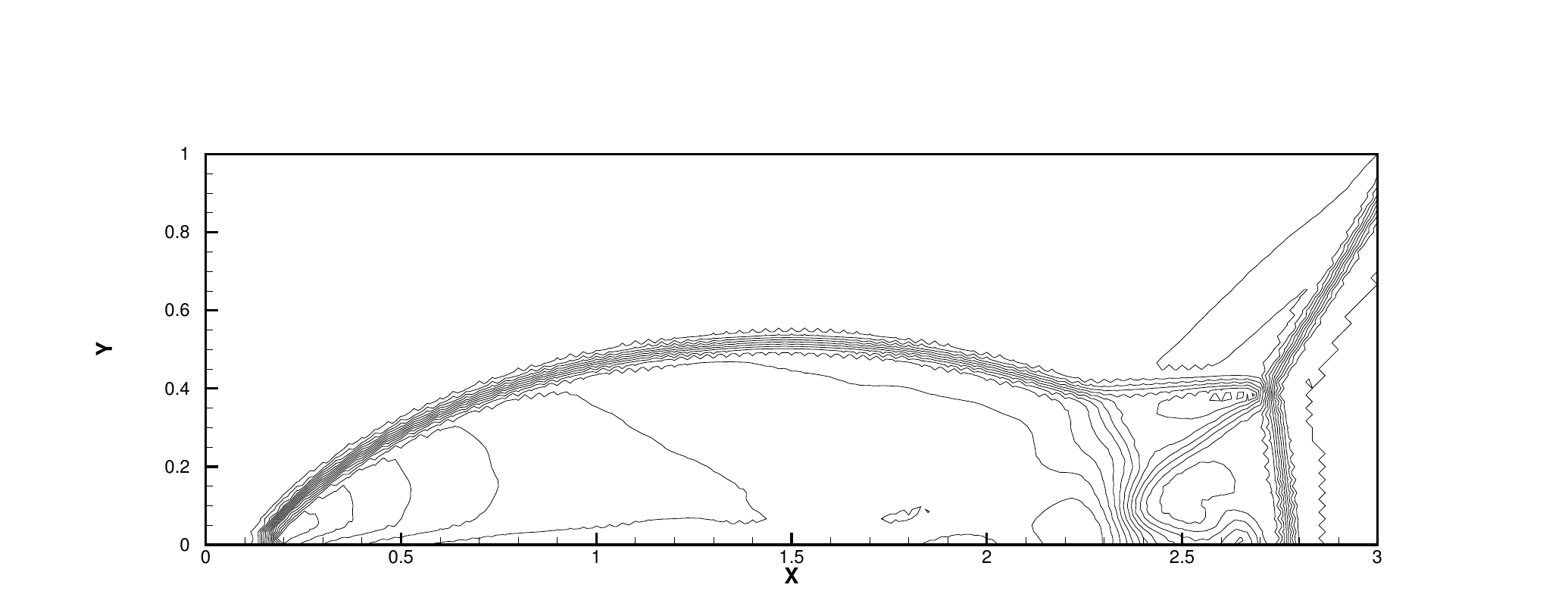}}\quad
    \subfigure[Fixed Mesh, $N = 180 \times 45 \times 4$] 
    {\includegraphics[width=8cm]{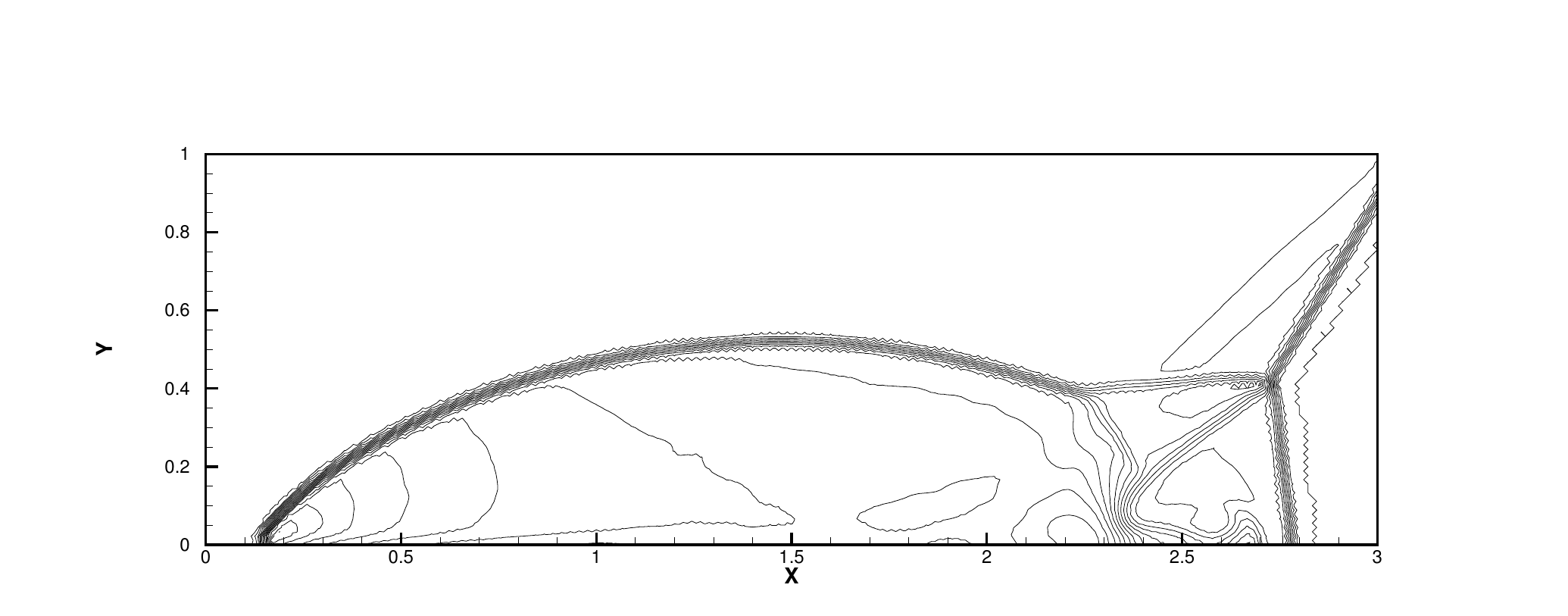}}
    }

  \mbox{\subfigure[Moving Mesh, $N = 120\times 30 \times 4$] 
  {\includegraphics[width=8cm]{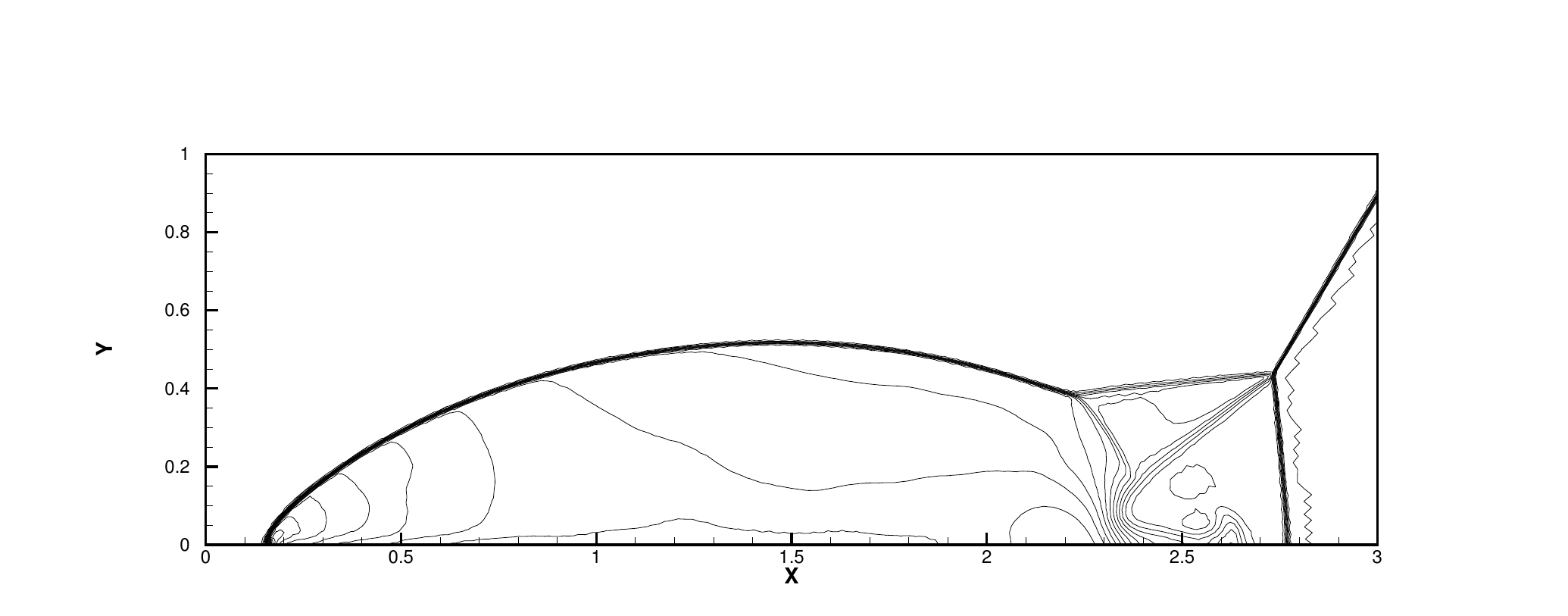}}\quad
    \subfigure[Moving mesh]
    {\includegraphics[width=8cm]{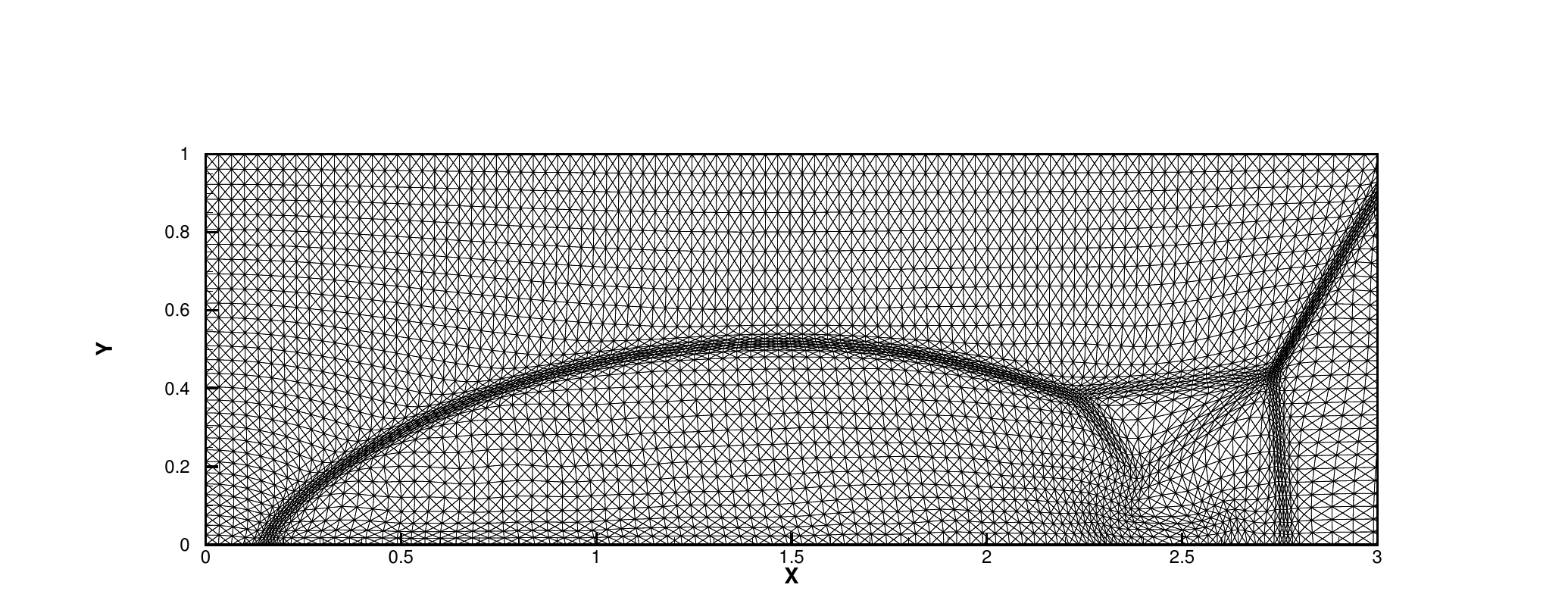}}
    }

    \caption{$P^1$ elements are used. 30 equally spaced density contours from 1.4 to 22.1183 are used in the contour plots.}
    \label{fig:edge11}
    \end{center}
    \end{figure}

\begin{figure}[hbtp]
  \begin{center}
  \mbox{\subfigure[Fixed Mesh, $N = 120\times 30 \times 4$] 
  {\includegraphics[width=8cm]{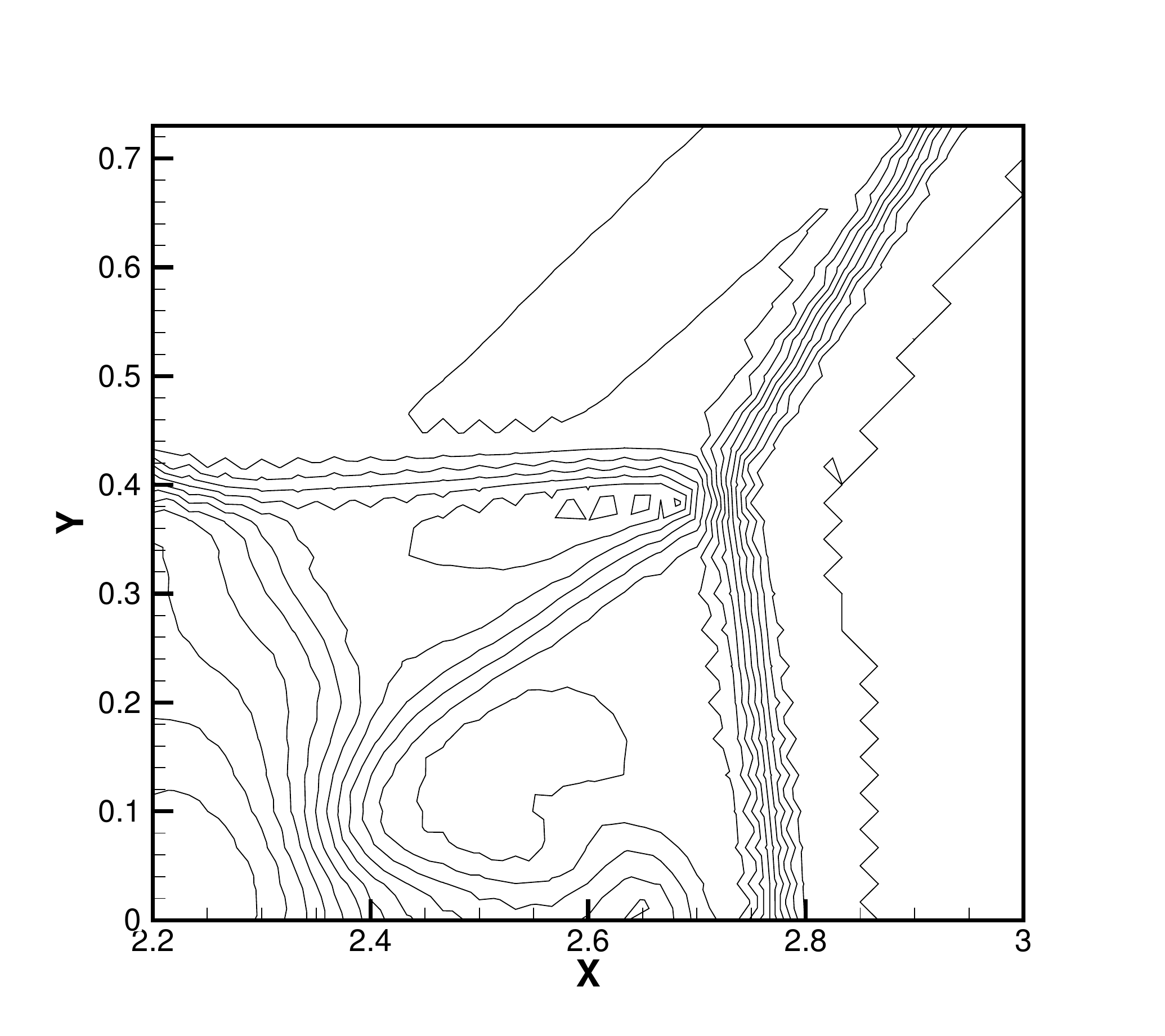}}\quad
    \subfigure[Fixed Mesh, $N = 180\times 45 \times 4$] 
    {\includegraphics[width=8cm]{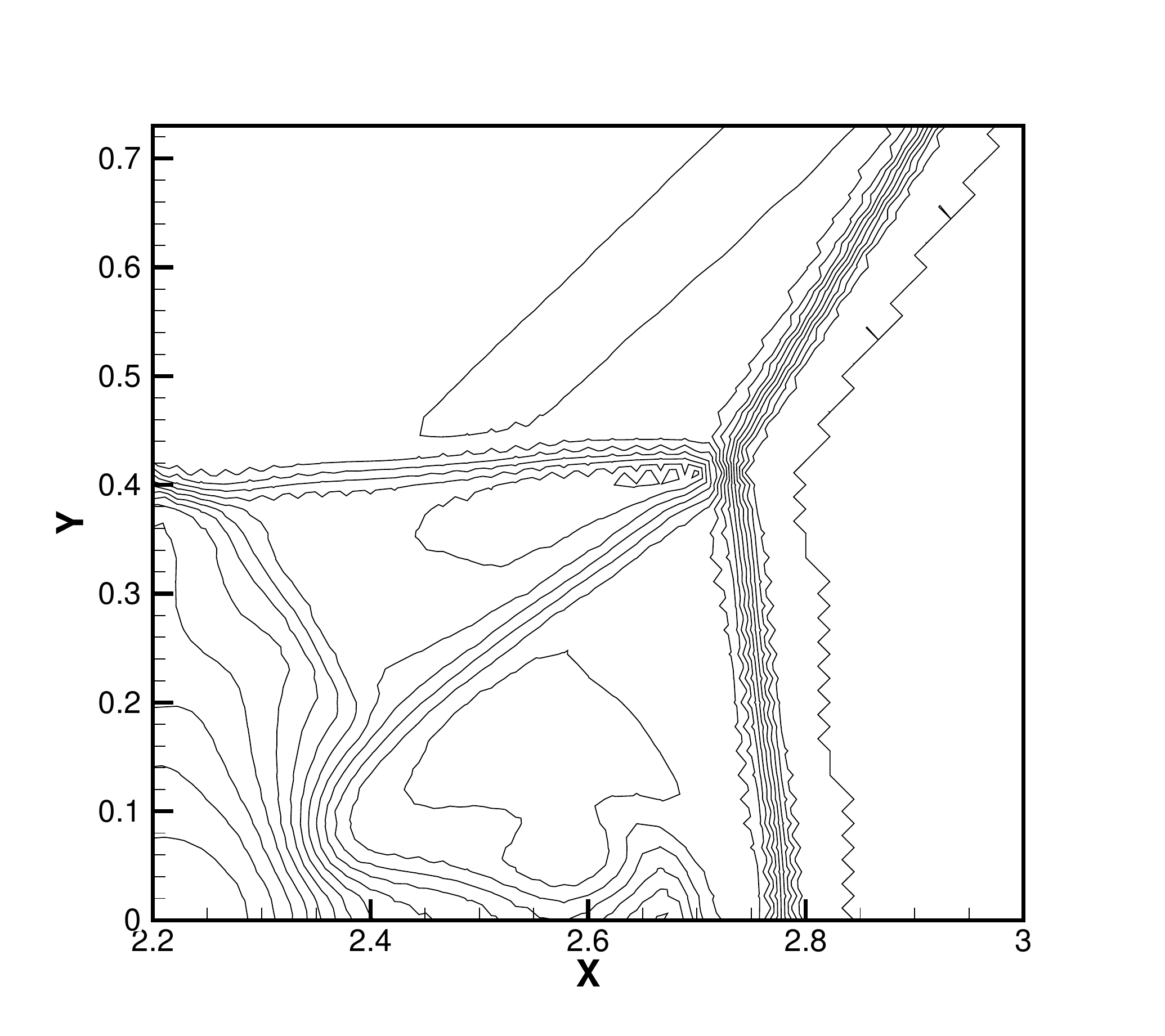}}
    }

  \mbox{\subfigure[Moving Mesh, $N = 120\times 30 \times 4$] 
  {\includegraphics[width=8cm]{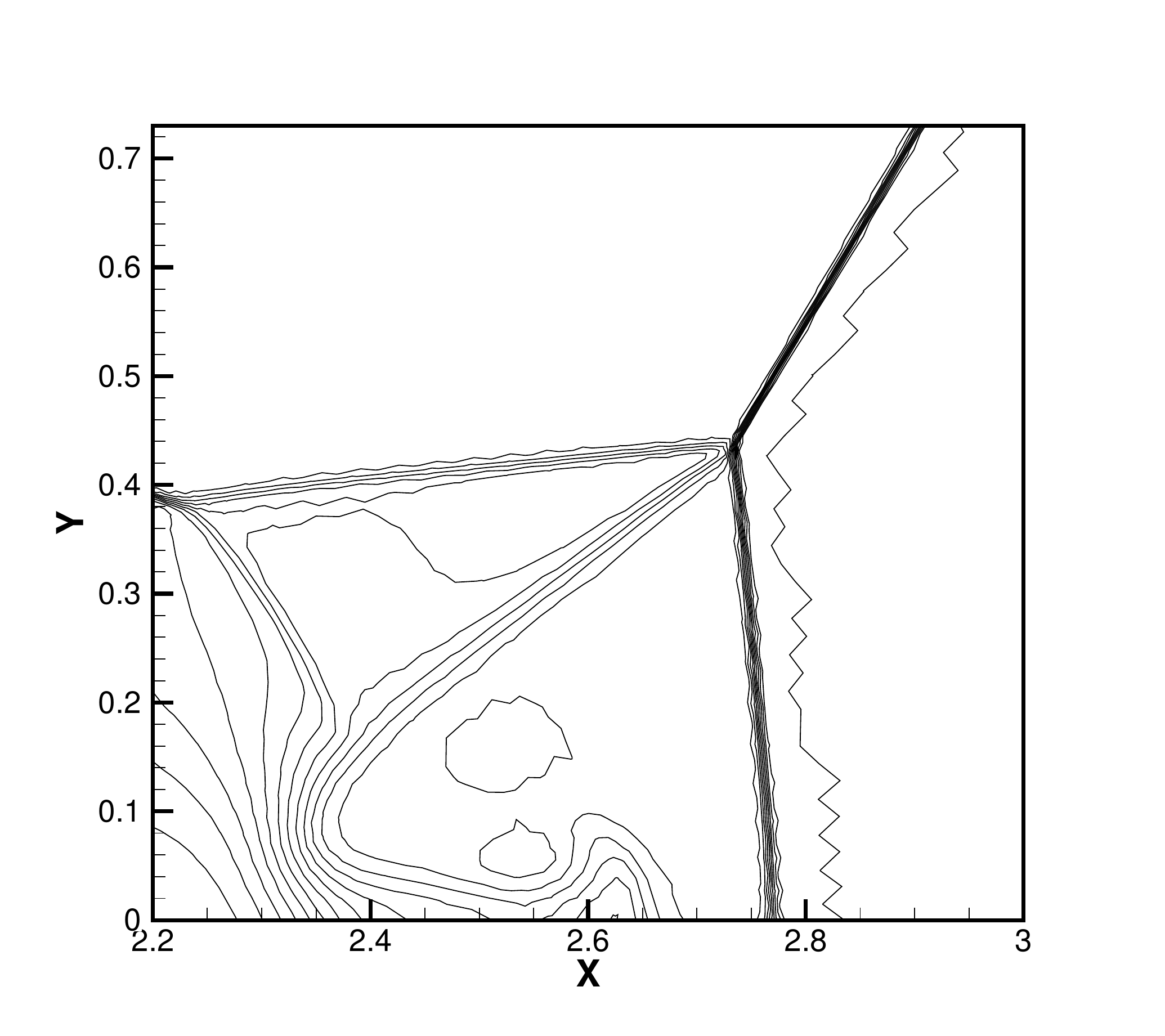}}\quad
    \subfigure[Moving mesh]
    {\includegraphics[width=8cm]{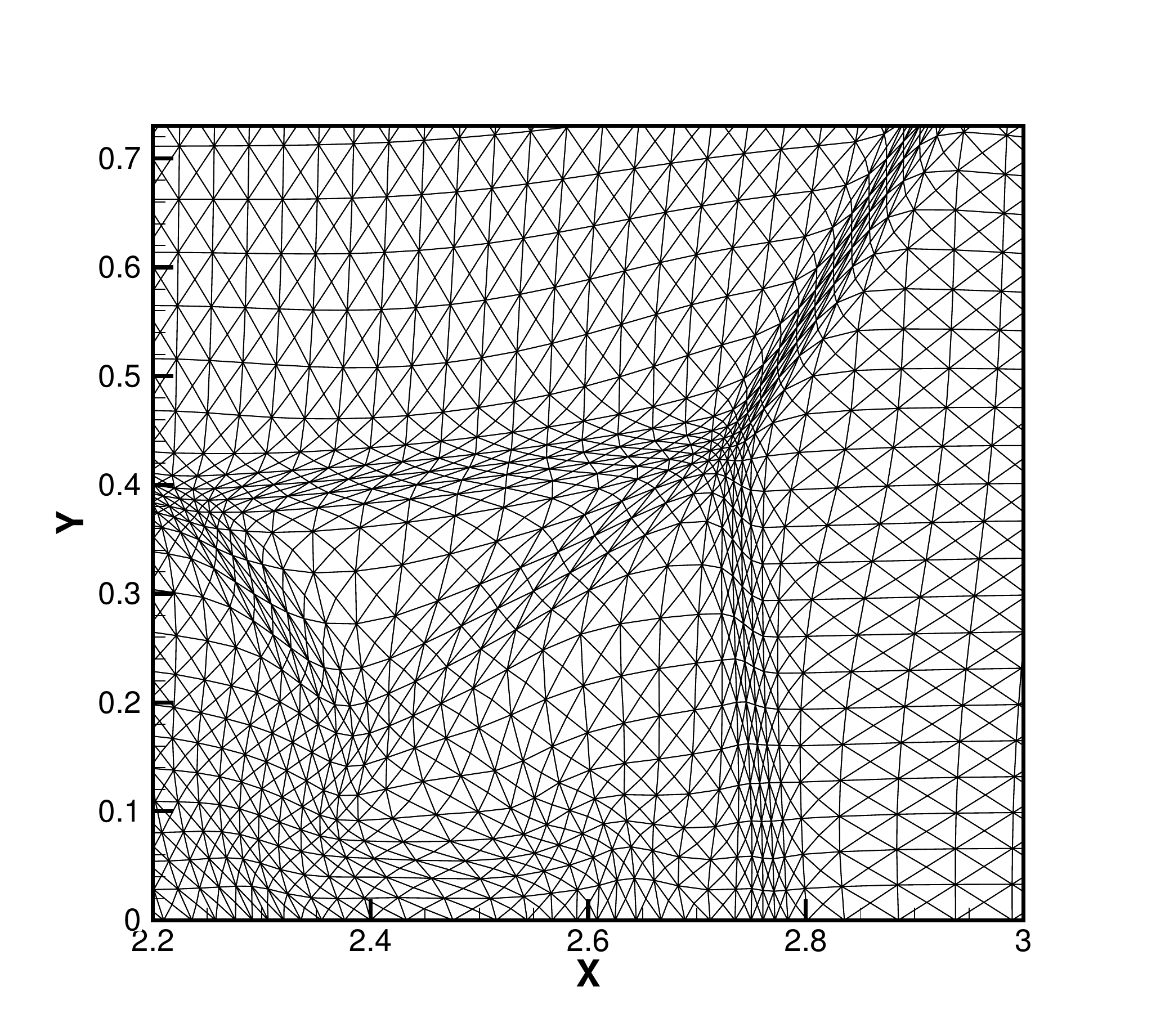}}
    }

    \caption{$P^1$ elements are used. Close view of the complex zone in Fig. \ref{fig:edge11}.}
    \label{fig:edge11zoom}
    \end{center}
    \end{figure}

\begin{figure}[hbtp]
  \begin{center}
  \mbox{\subfigure[Fixed Mesh, $N = 240\times 60 \times 4$] 
  {\includegraphics[width=8cm]{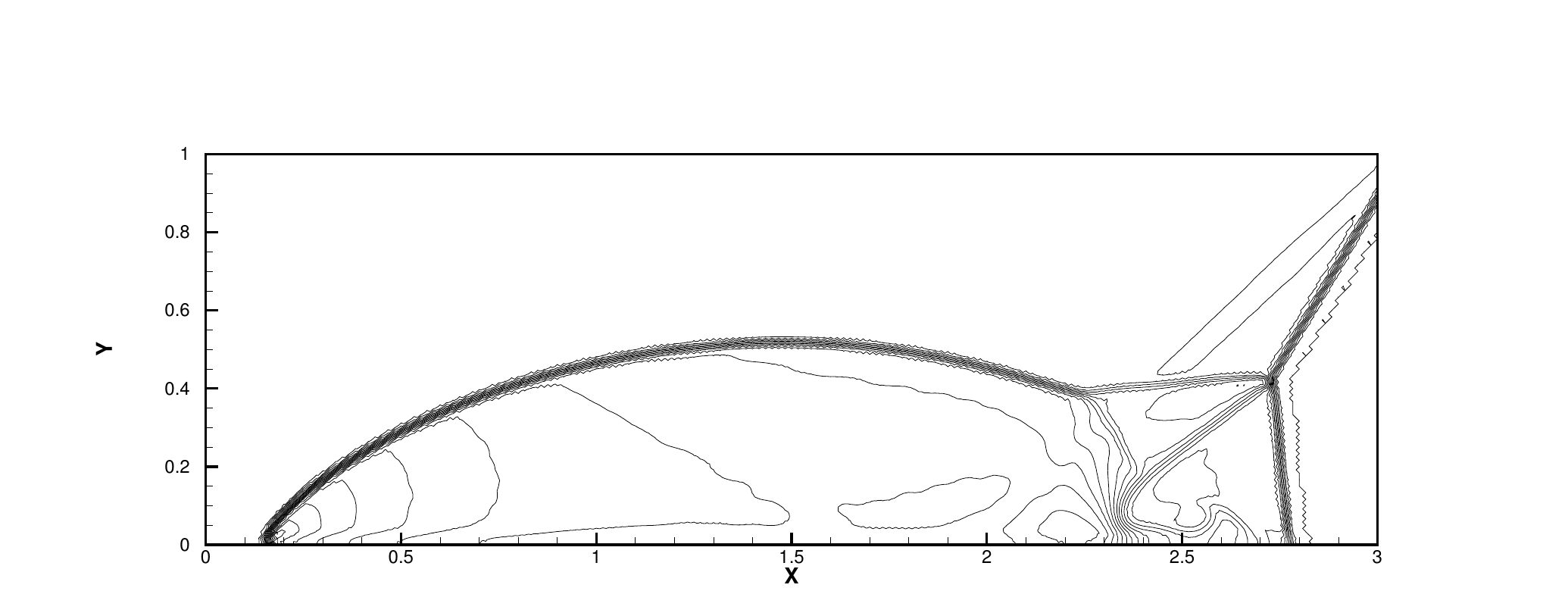}}\quad
    \subfigure[Fixed Mesh, $N = 360\times 90 \times 4$] 
    {\includegraphics[width=8cm]{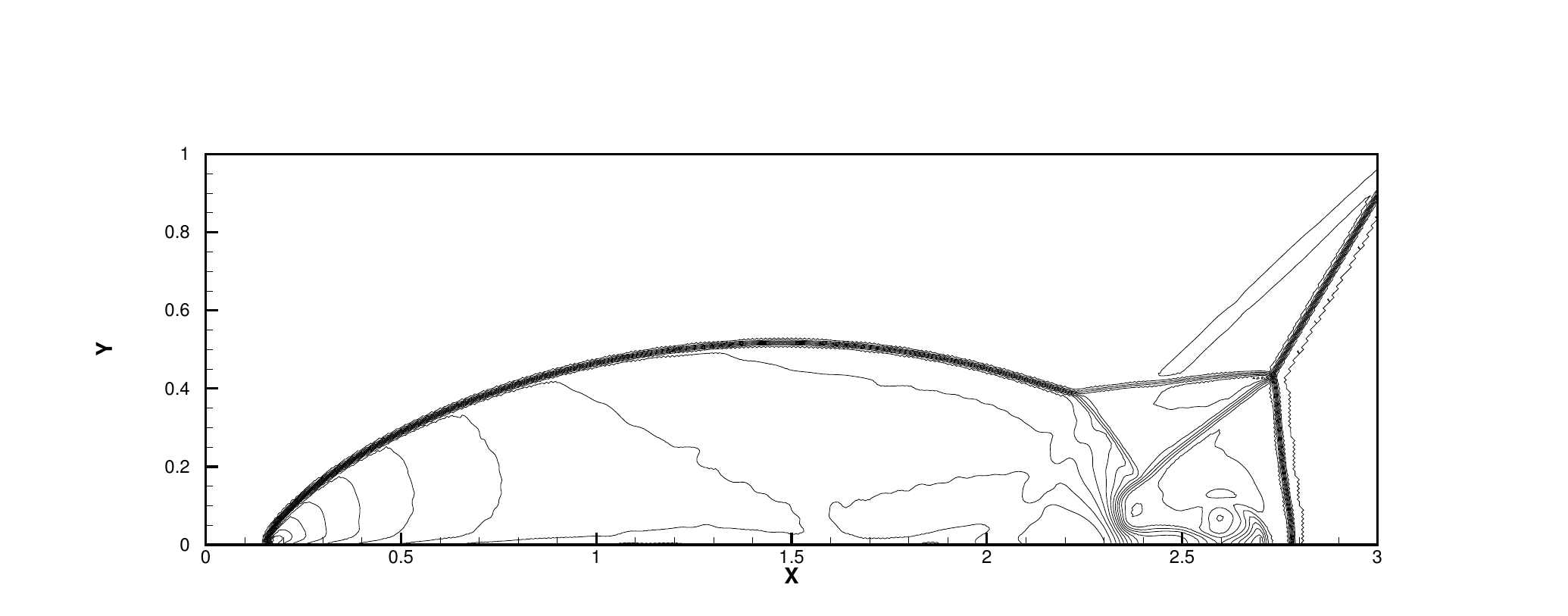}}
    }

  \mbox{\subfigure[Moving Mesh, $N = 240\times 60 \times 4$] 
  {\includegraphics[width=8cm]{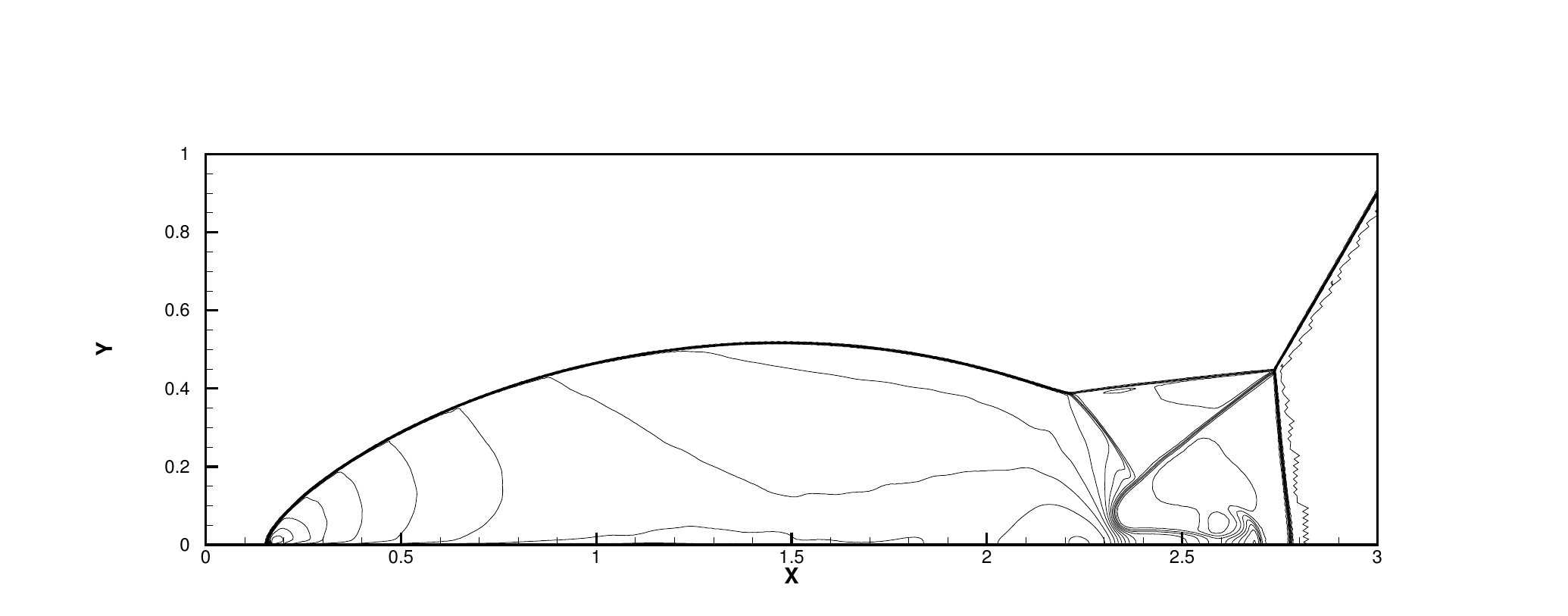}}\quad
    \subfigure[Moving mesh]
    {\includegraphics[width=8cm]{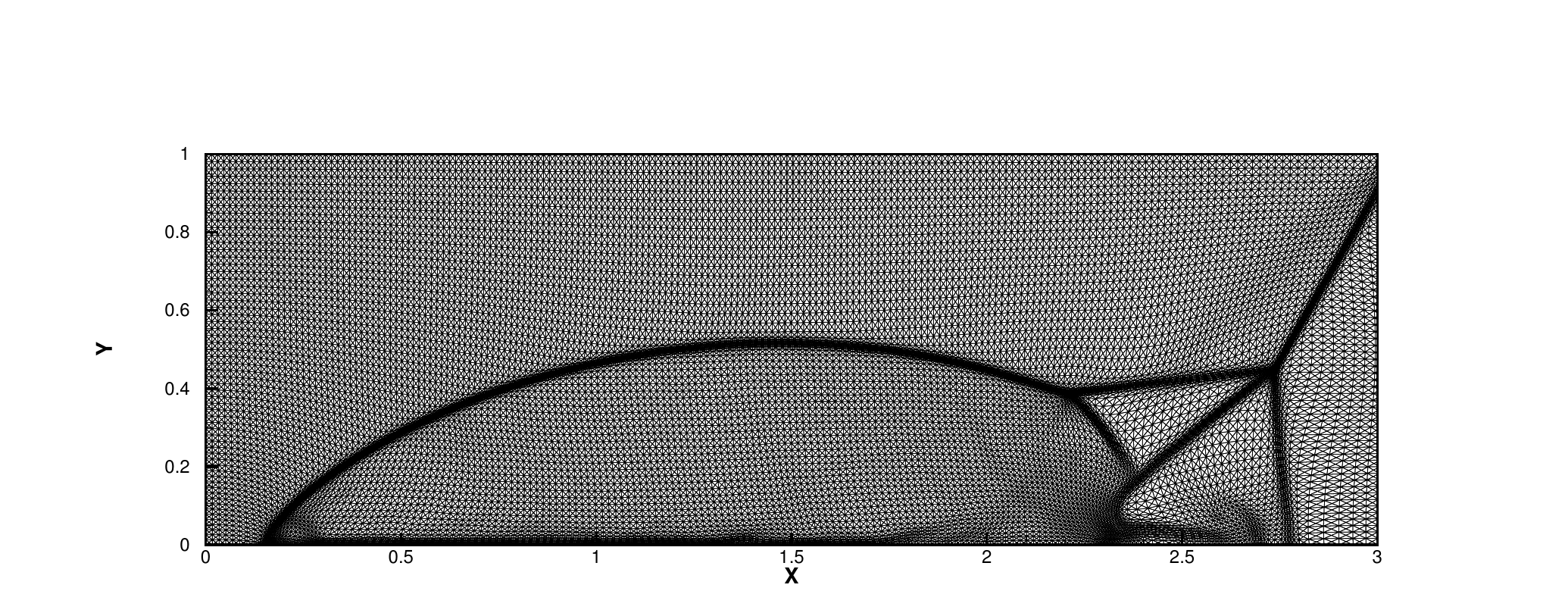}}
    }

    \caption{$P^1$ elements are used. 30 equally spaced density contours from 1.4 to 22.1183 are used
    in the contour plots.}
    \label{fig:edge12}
    \end{center}
    \end{figure}

\begin{figure}[hbtp]
  \begin{center}
  \mbox{\subfigure[Fixed Mesh, $N = 240\times 60 \times 4$] 
  {\includegraphics[width=8cm]{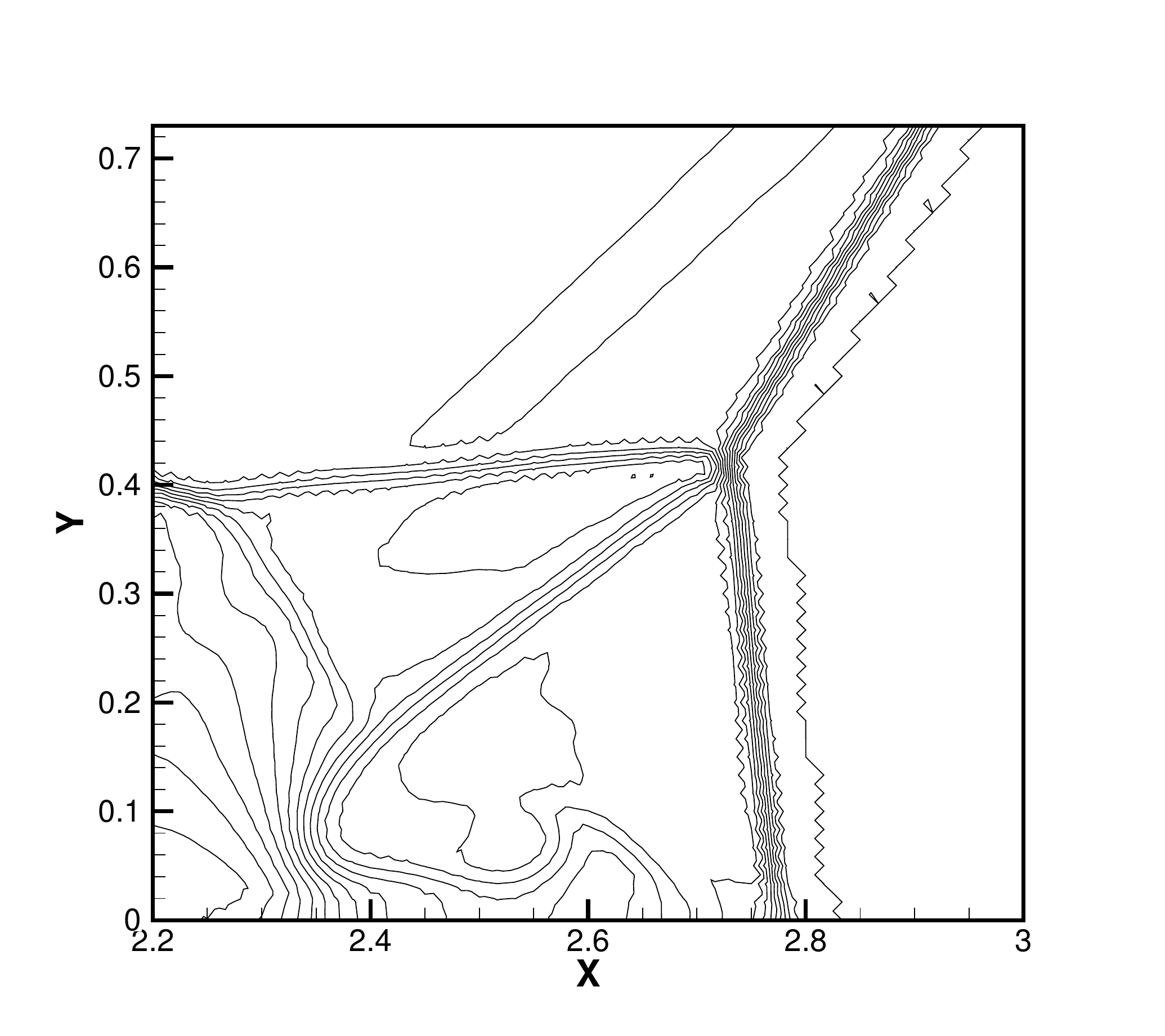}}\quad
    \subfigure[Fixed Mesh, $N = 360\times 90 \times 4$] 
    {\includegraphics[width=8cm]{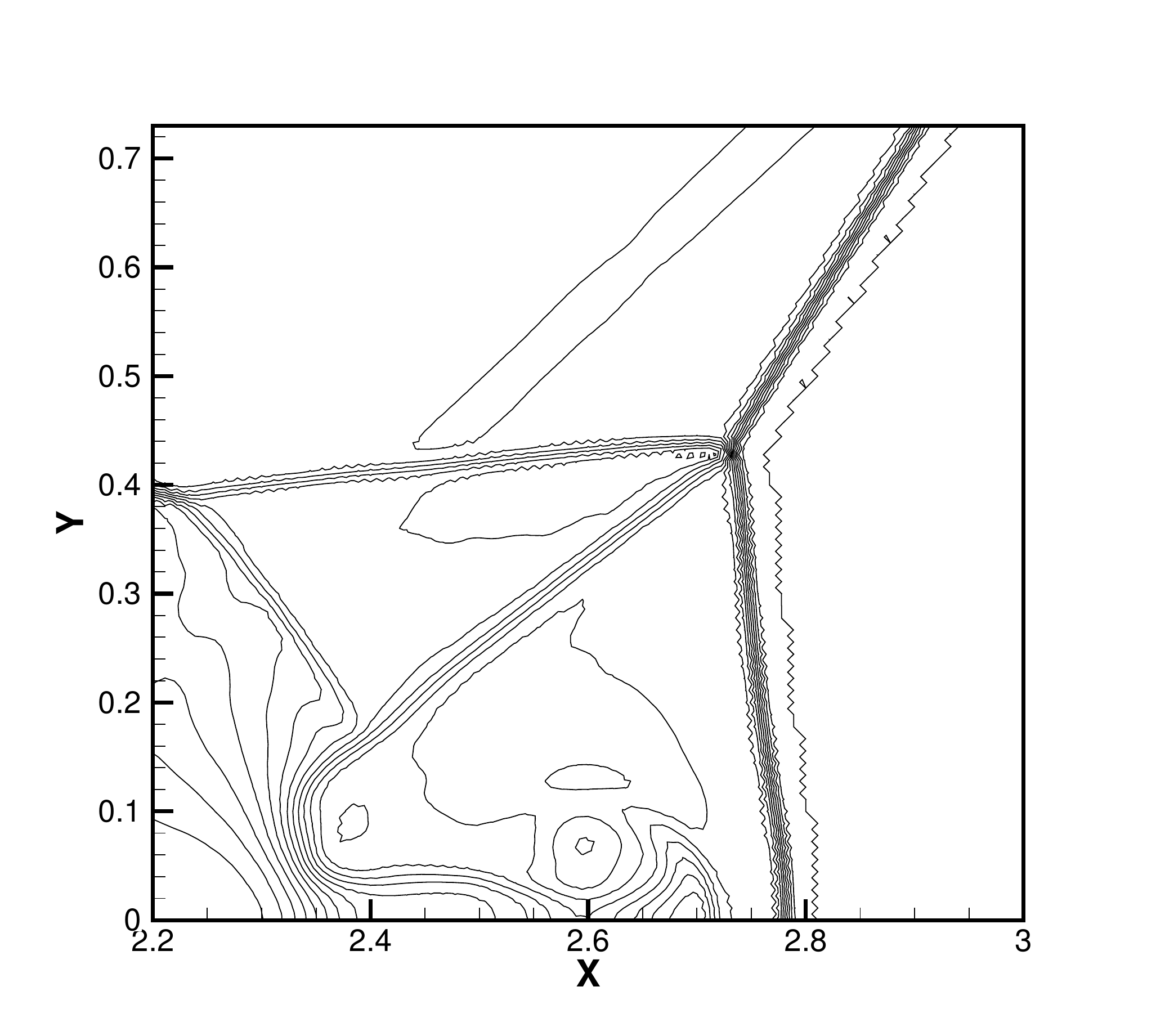}}
    }

  \mbox{\subfigure[Moving Mesh, $N = 240\times 60 \times 4$] 
  {\includegraphics[width=8cm]{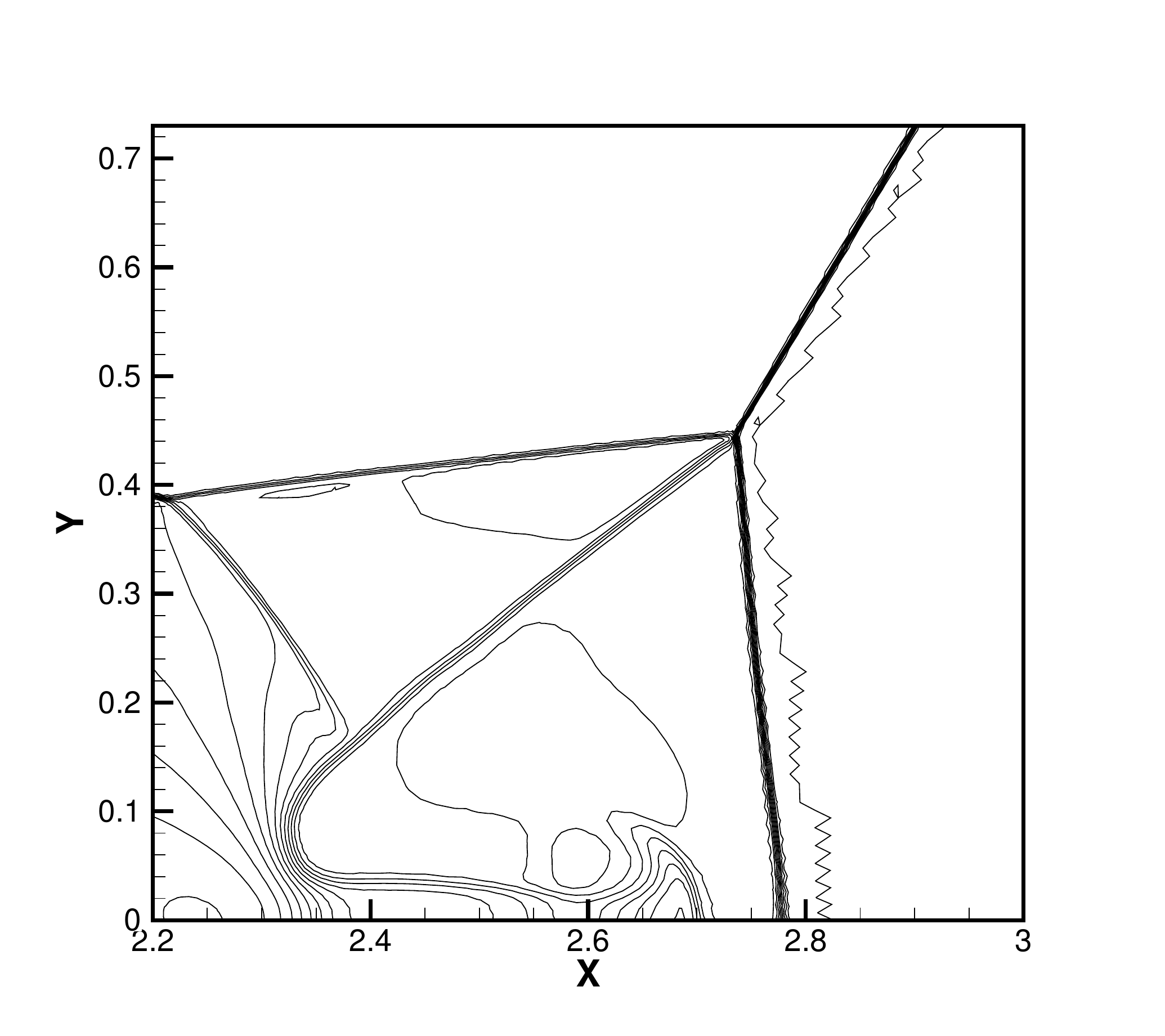}}\quad
    \subfigure[Moving mesh]
    {\includegraphics[width=8cm]{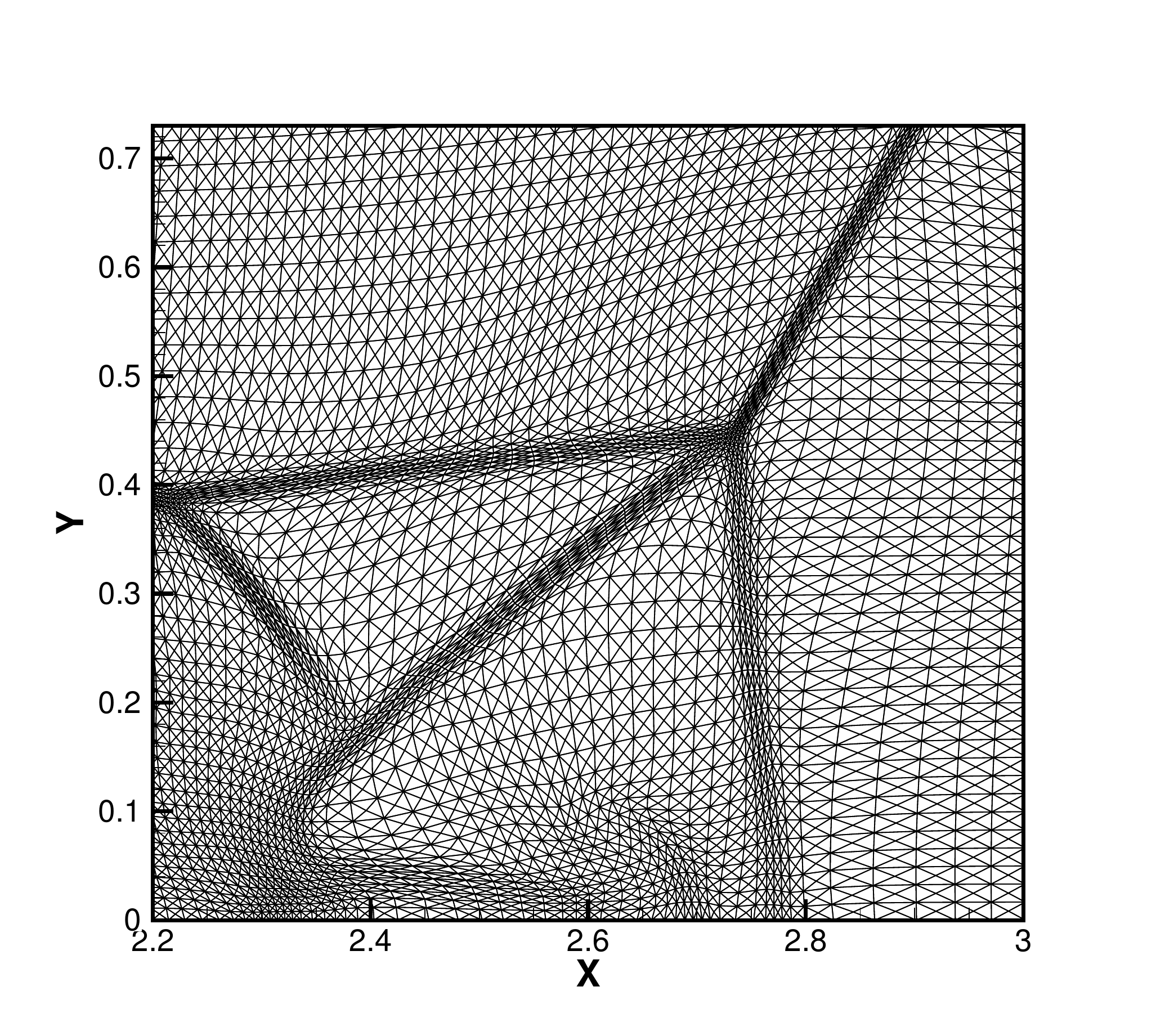}}
    }

    \caption{$P^1$ elements are used. Close view of the complex zone in Fig. \ref{fig:edge12}.}
    \label{fig:edge12zoom}
    \end{center}
    \end{figure}

\begin{figure}[hbtp]
  \begin{center}
  \mbox{\subfigure[Fixed Mesh, $N = 120\times 30 \times 4$] 
  {\includegraphics[width=8cm]{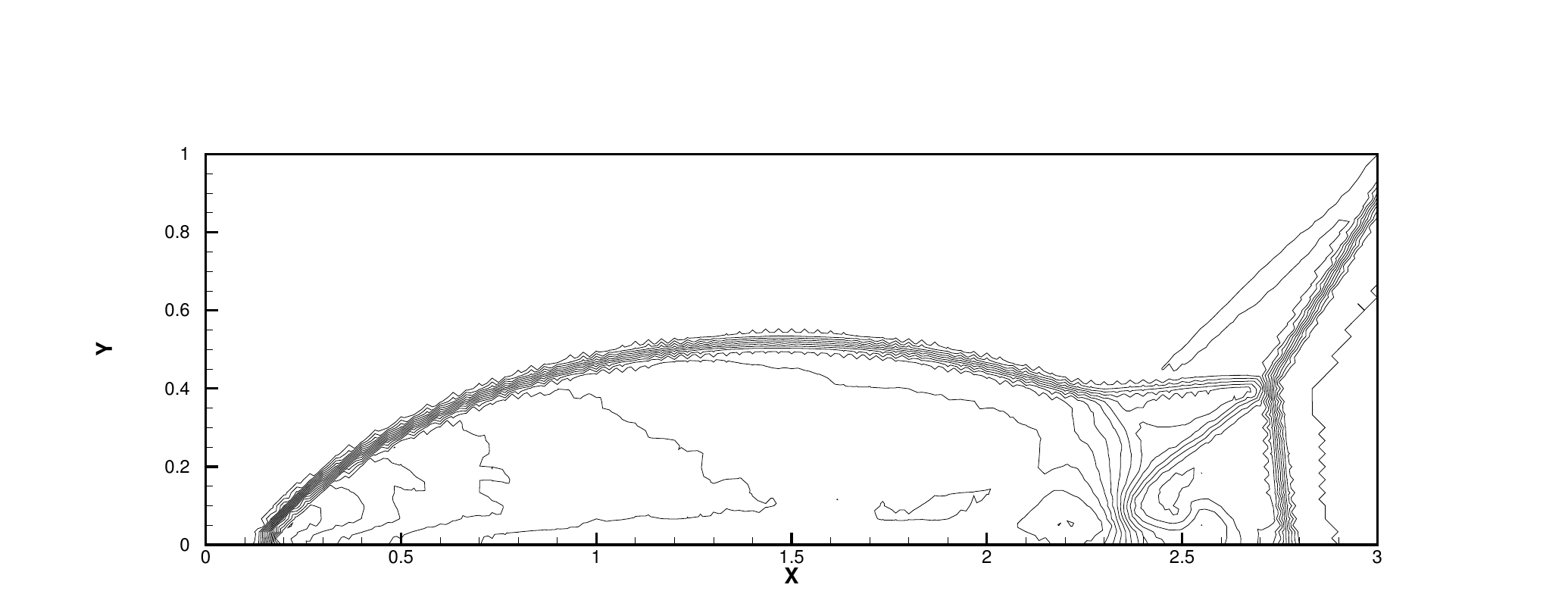}}\quad
    \subfigure[Fixed Mesh, $N = 180\times 45 \times 4$] 
    {\includegraphics[width=8cm]{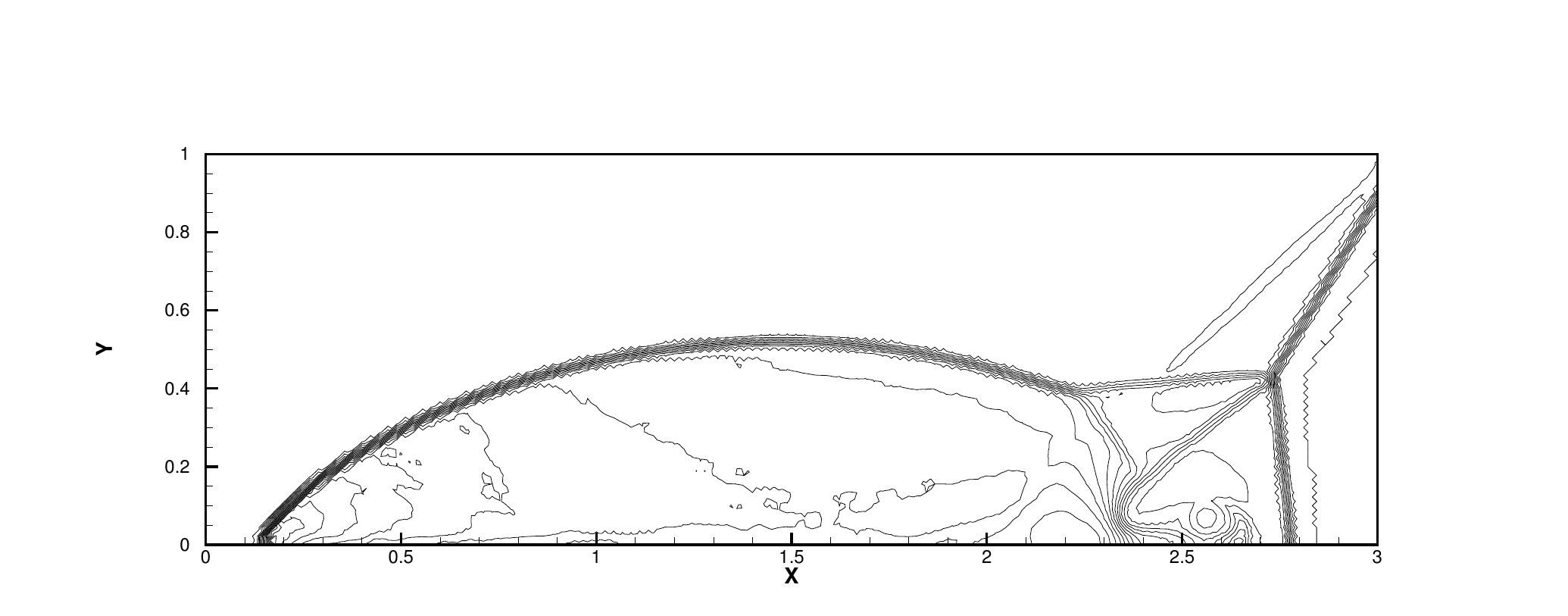}}
    }

  \mbox{\subfigure[Moving Mesh, $N = 120\times 30 \times 4$] 
  {\includegraphics[width=8cm]{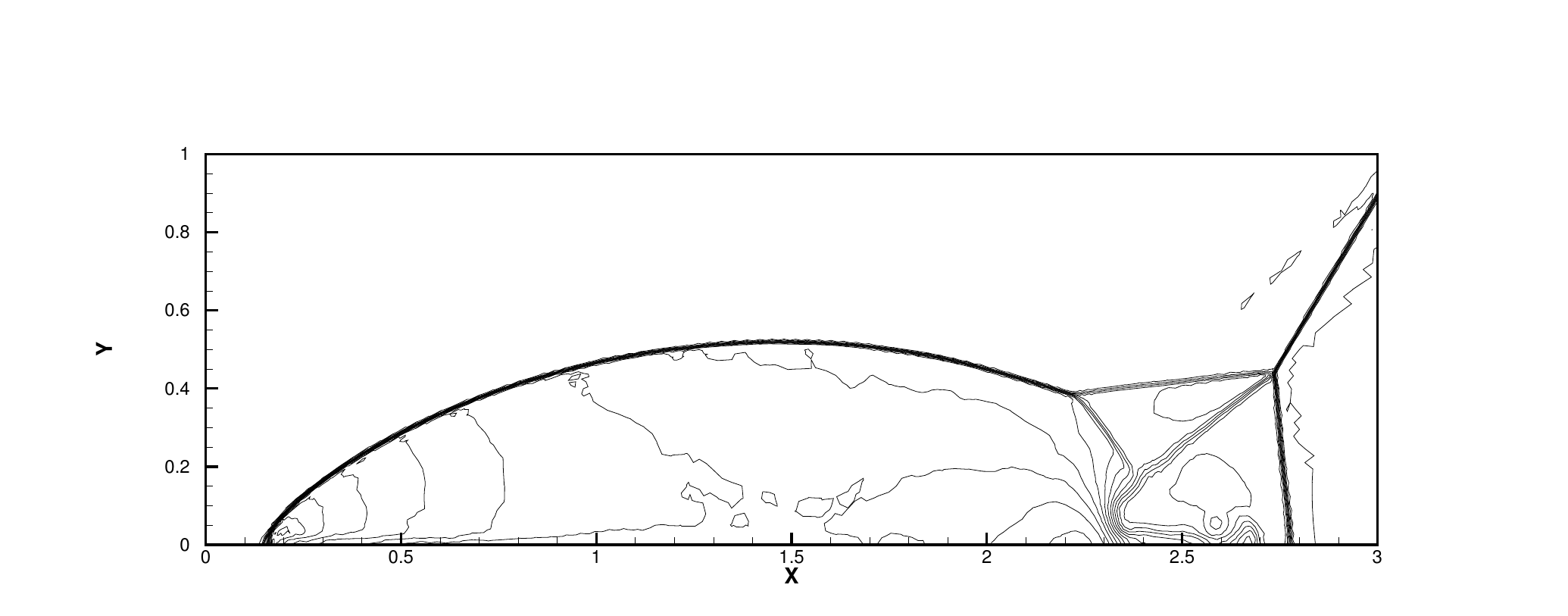}}\quad
    \subfigure[Moving mesh]
    {\includegraphics[width=8cm]{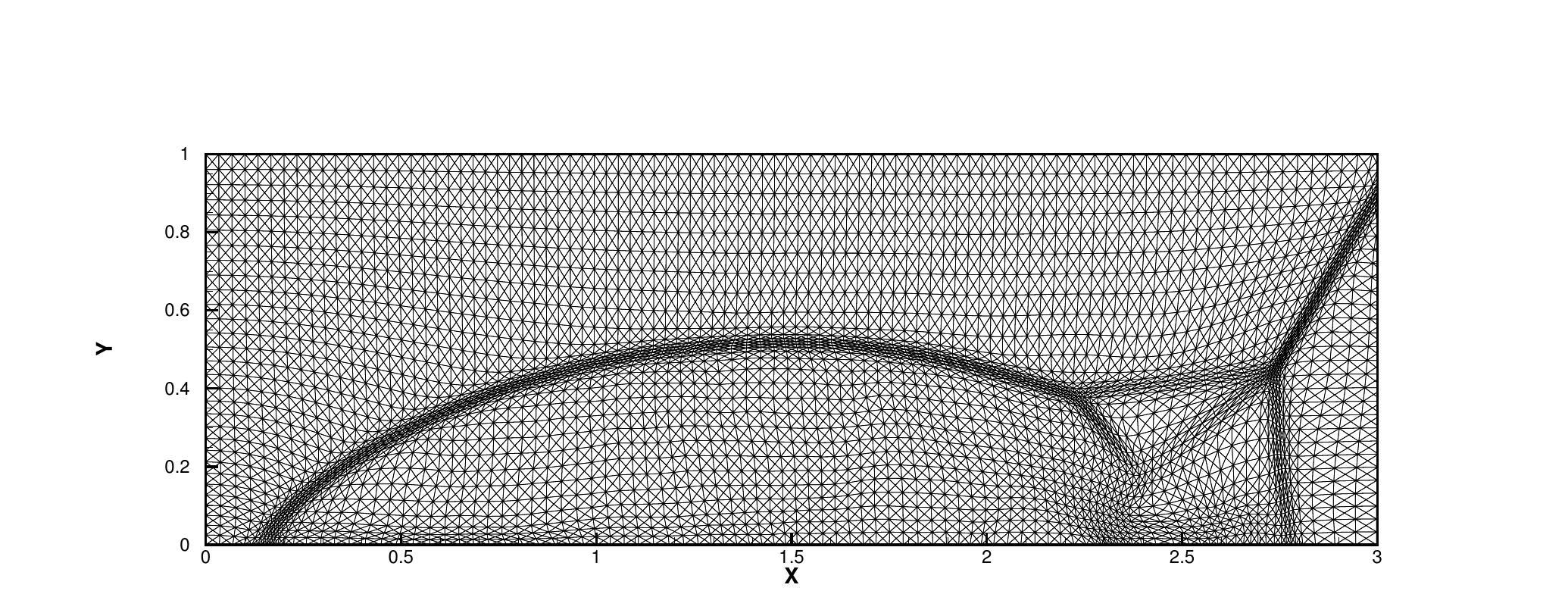}}
    }

    \caption{$P^2$ elements are used. 30 equally spaced density contours from 1.4 to 22.1183 are used in the contour plots.}
    \label{fig:edge13}
    \end{center}
    \end{figure}

\begin{figure}[hbtp]
  \begin{center}
  \mbox{\subfigure[Fixed Mesh, $N = 120\times 30 \times 4$]
  {\includegraphics[width=8cm]{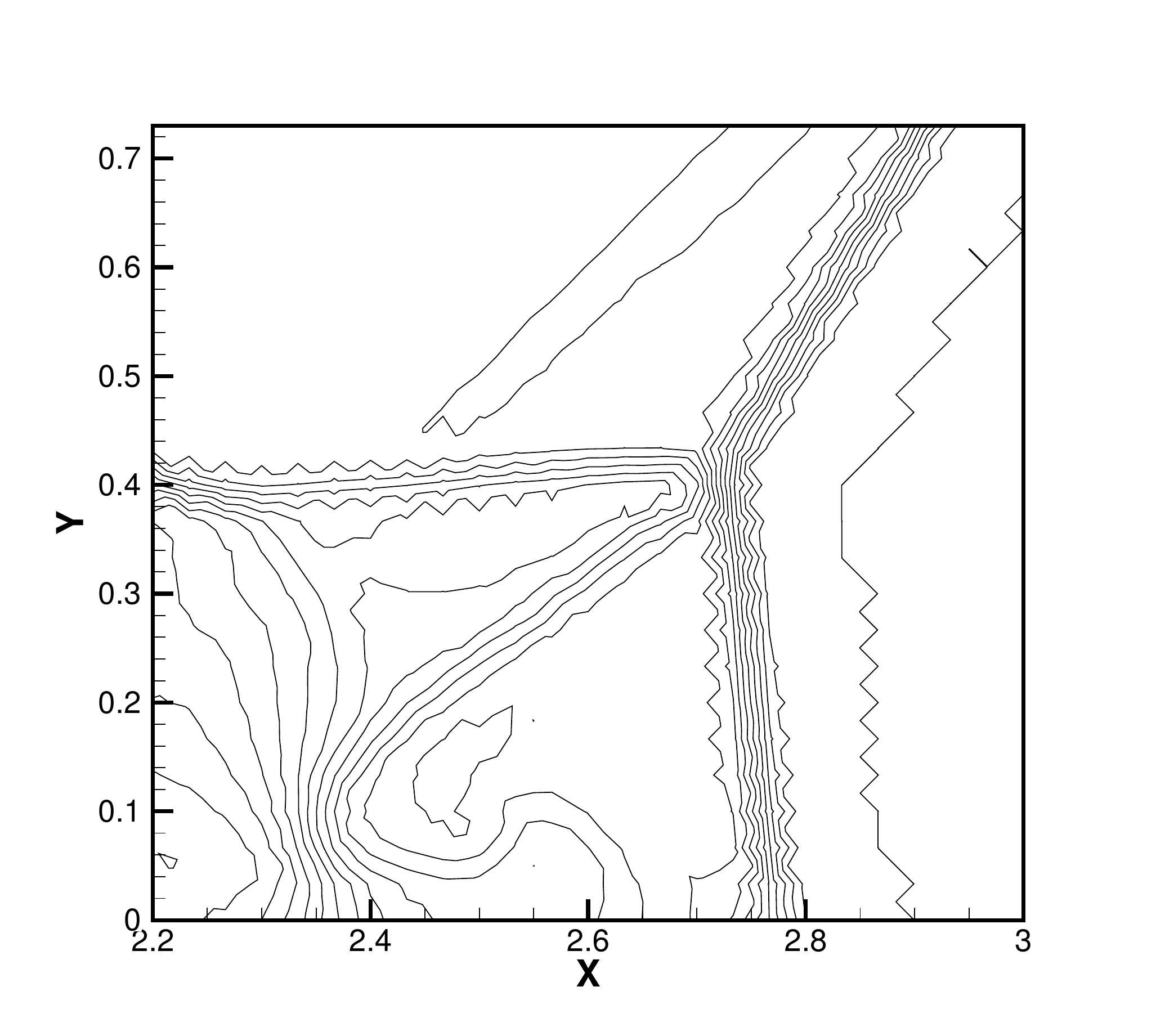}}\quad
    \subfigure[Fixed Mesh, $N = 180\times 45 \times 4$] 
    {\includegraphics[width=8cm]{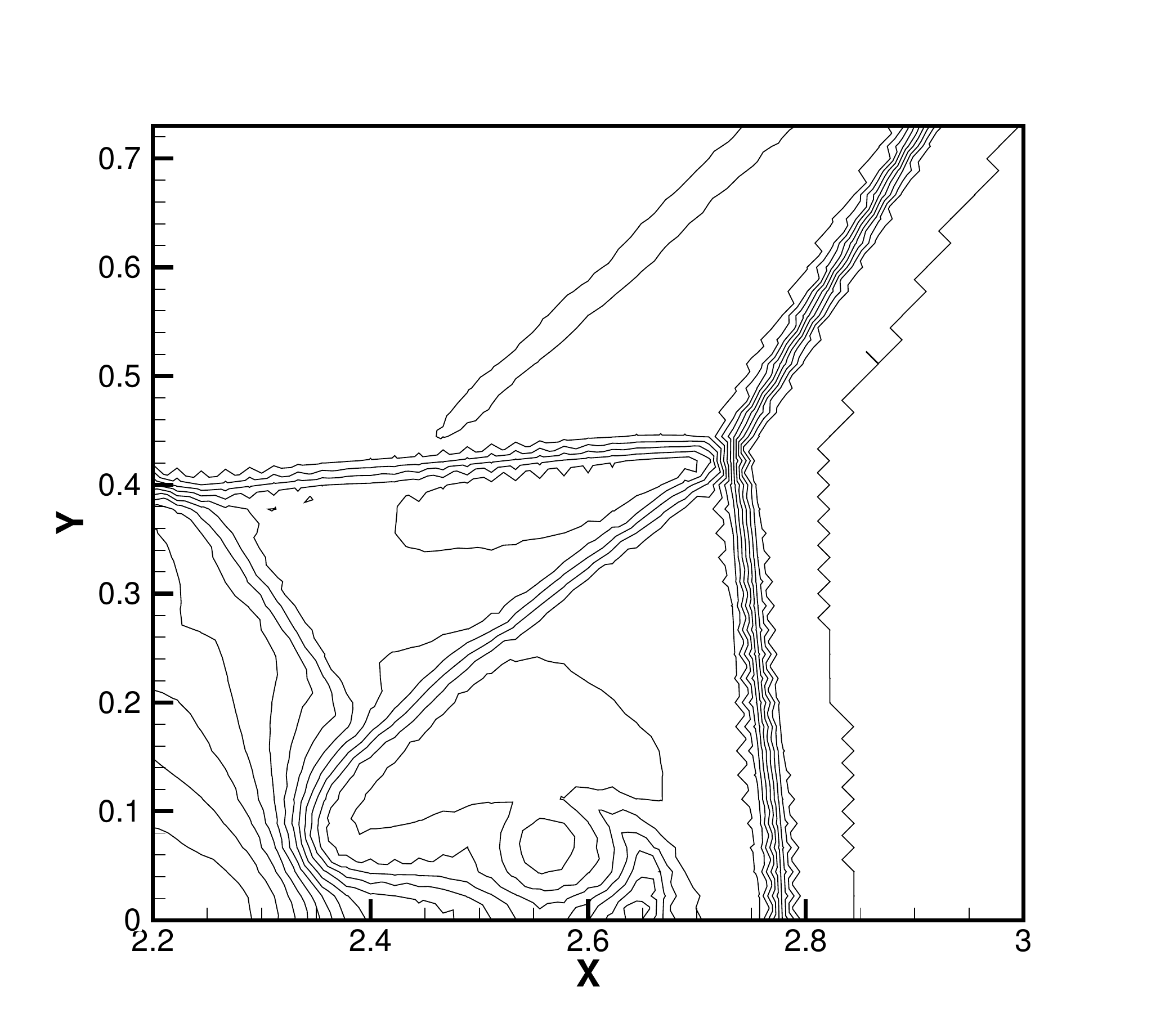}}
    }

  \mbox{\subfigure[Moving Mesh, $N = 120\times 30 \times 4$]
  {\includegraphics[width=8cm]{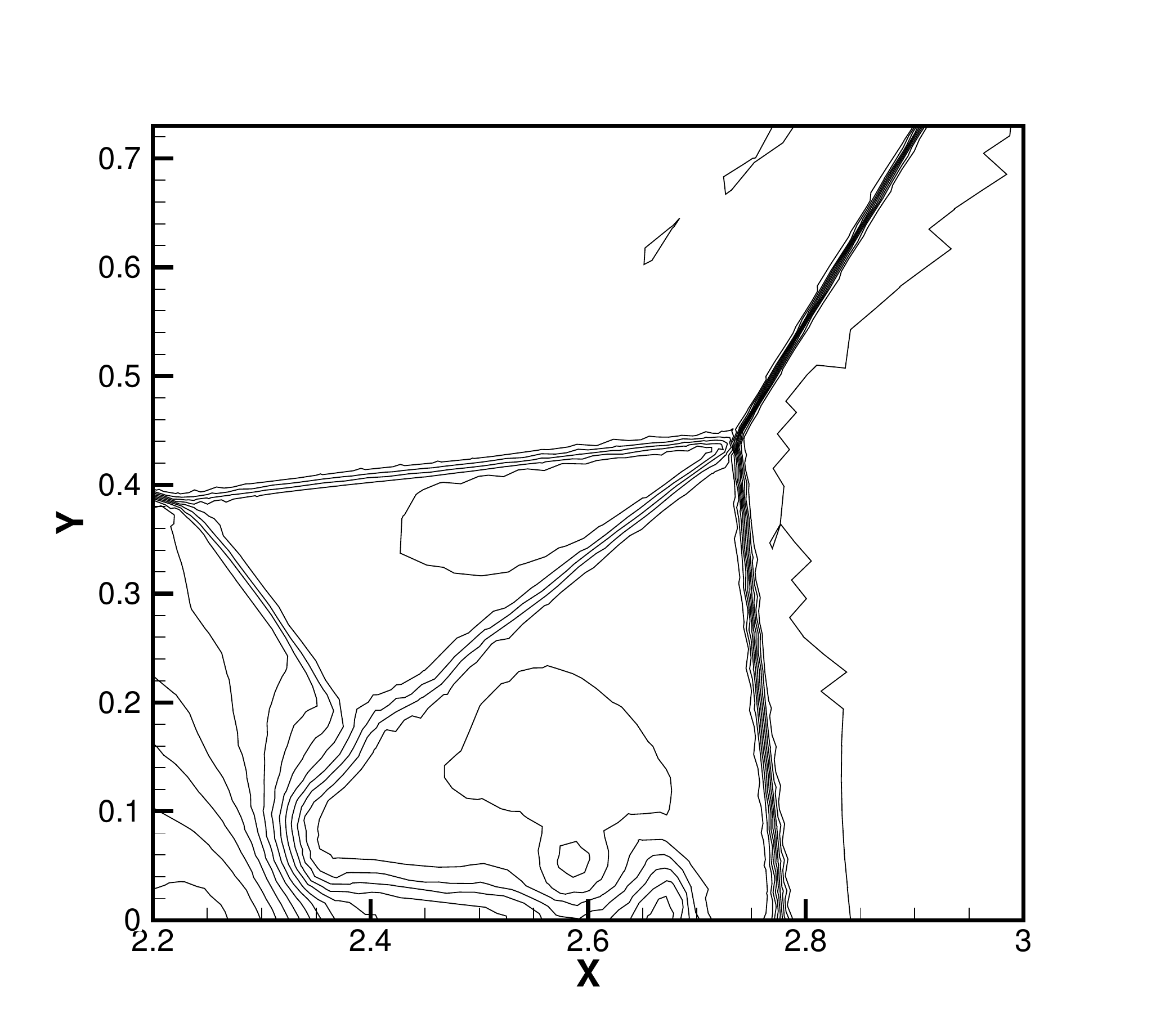}}\quad
    \subfigure[Moving mesh]
    {\includegraphics[width=8cm]{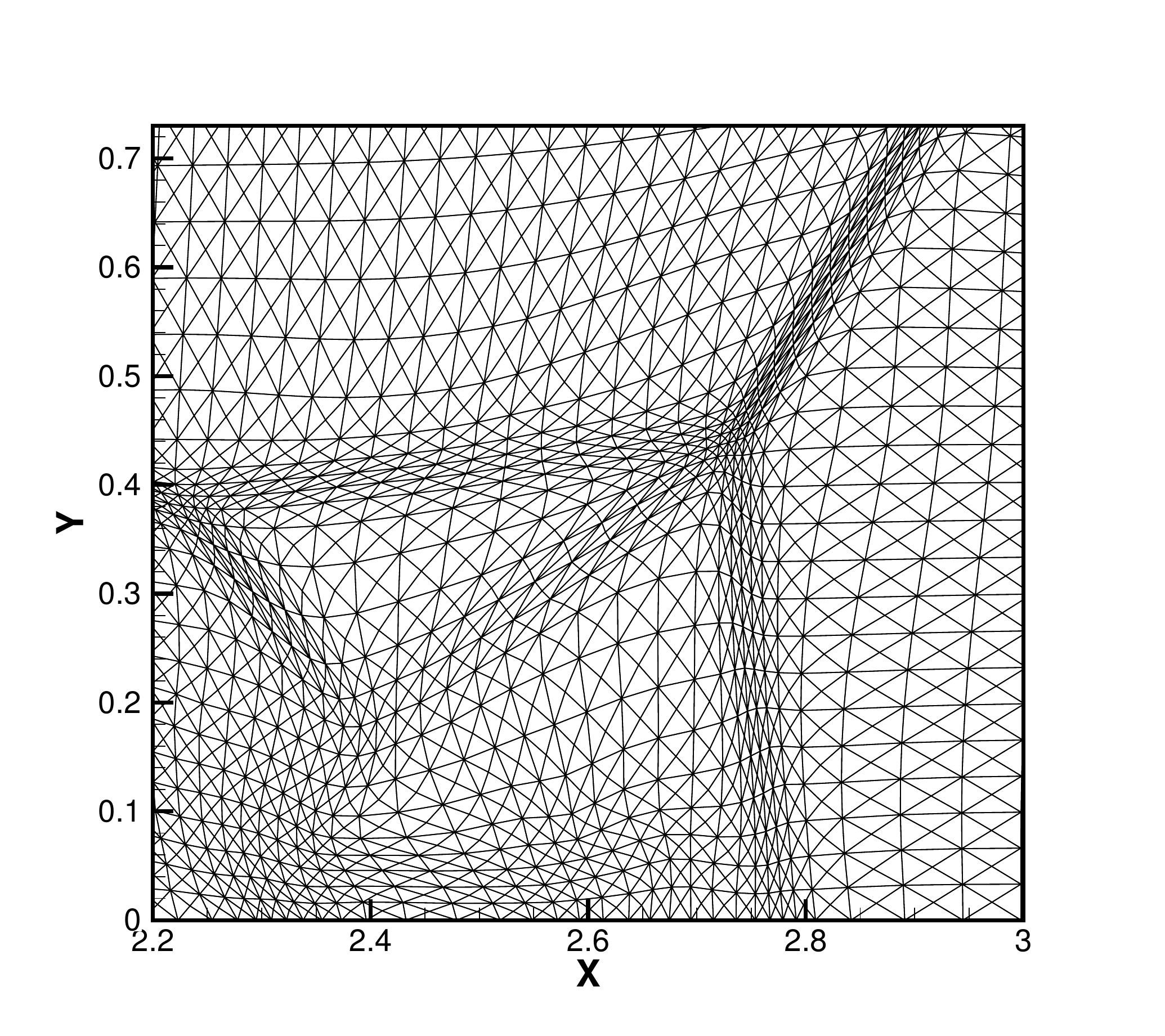}}
    }

    \caption{$P^2$ elements are used. Close view of the complex zone in Fig. \ref{fig:edge13}.}
    \label{fig:edge13zoom}
    \end{center}
    \end{figure}

\begin{figure}[hbtp]
  \begin{center}
  \mbox{\subfigure[Fixed Mesh, $N = 240\times 60 \times 4$]
  {\includegraphics[width=8cm]{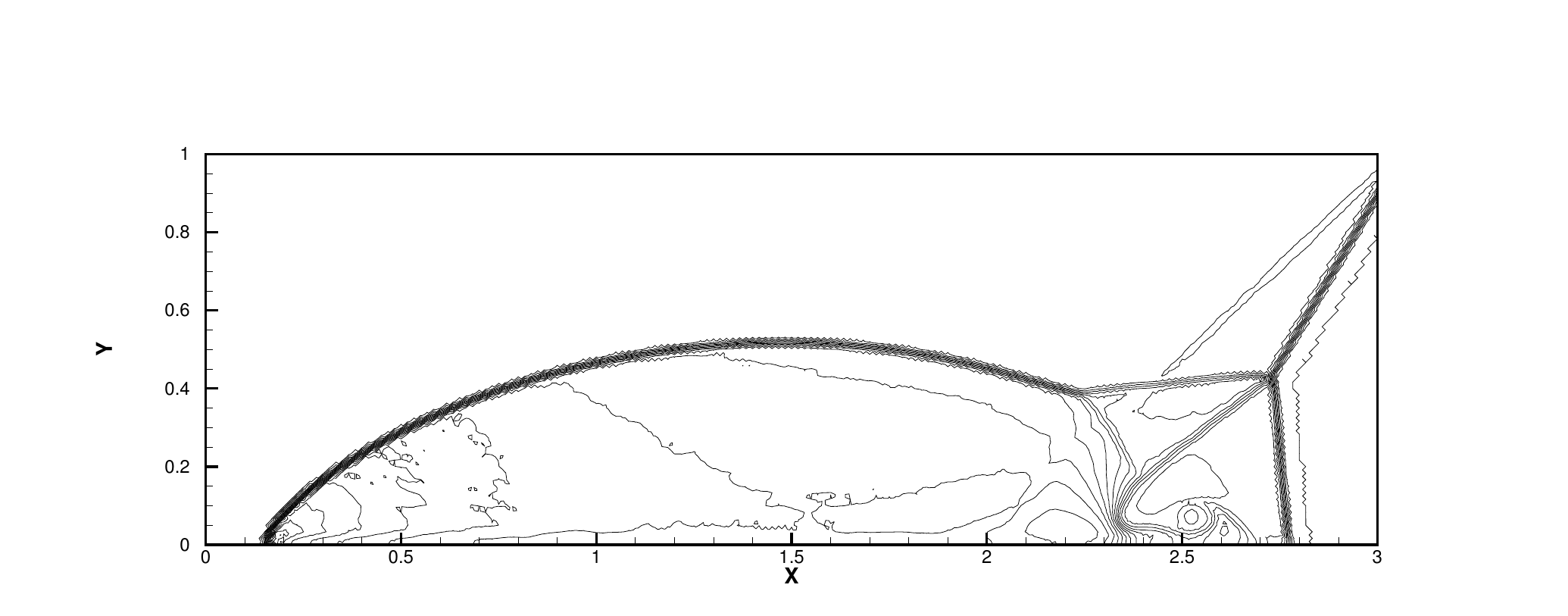}}\quad
    \subfigure[Fixed Mesh, $N = 360\times 90 \times 4$]
    {\includegraphics[width=8cm]{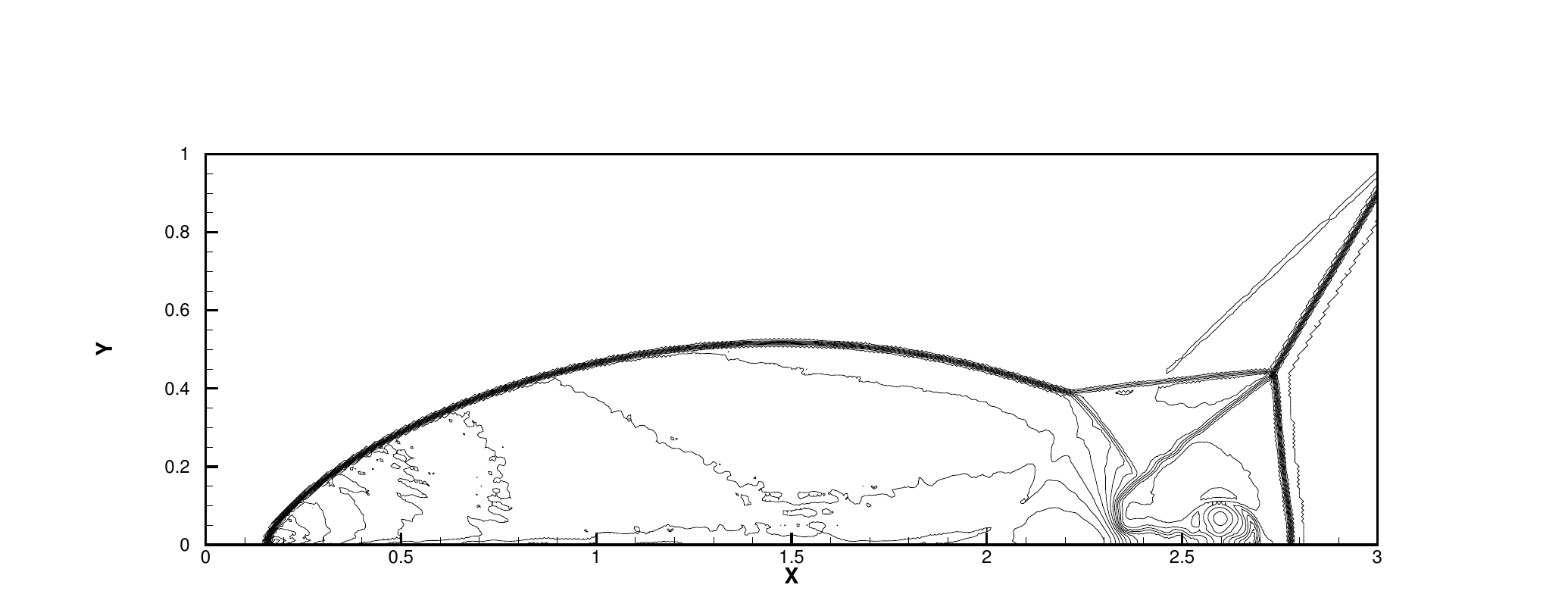}}
    }

  \mbox{\subfigure[Moving Mesh, $N = 240\times 60 \times 4$]
  {\includegraphics[width=8cm]{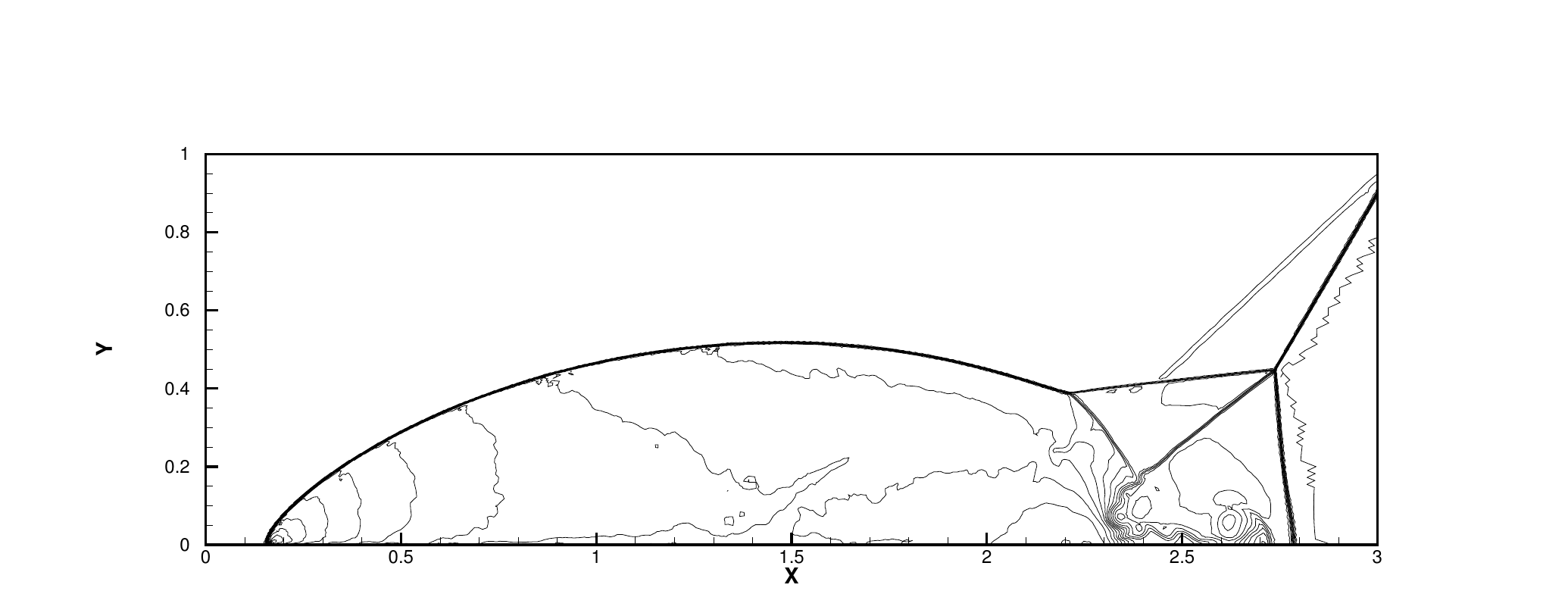}}\quad
    \subfigure[Moving mesh]
    {\includegraphics[width=8cm]{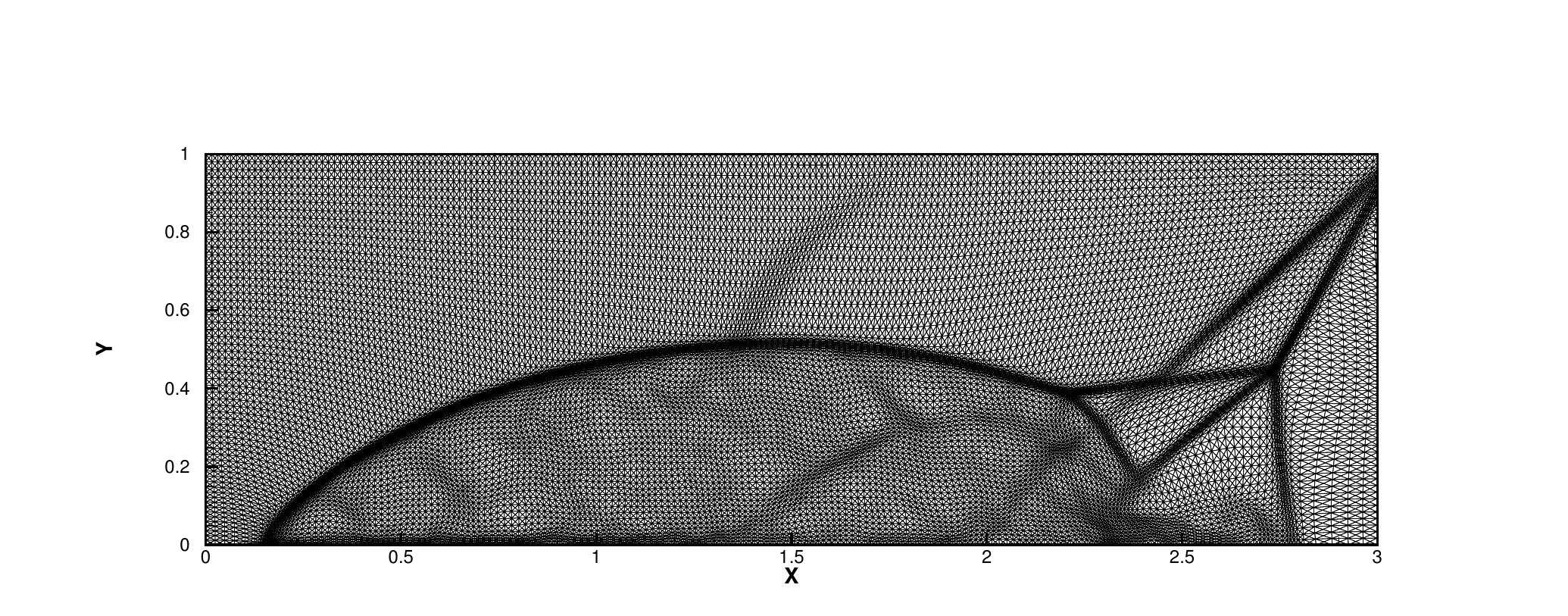}}
    }

    \caption{$P^2$ elements are used. 30 equally spaced density contours from 1.4 to 22.1183 are used in the contour plots.}
    \label{fig:edge14}
    \end{center}
    \end{figure}

\begin{figure}[hbtp]
  \begin{center}
  \mbox{\subfigure[Fixed Mesh, $N = 240\times 60 \times 4$]
  {\includegraphics[width=8cm]{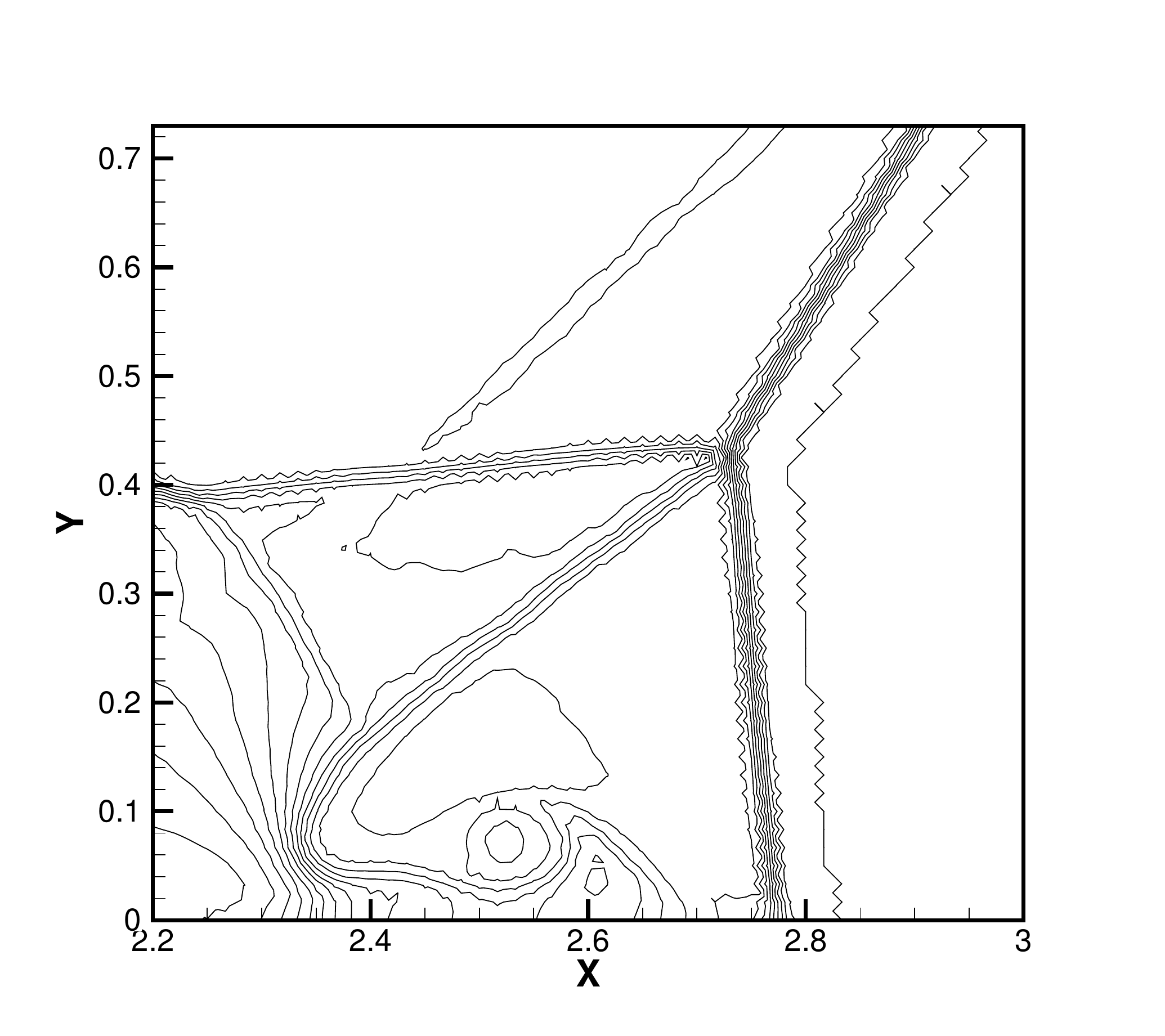}}\quad
    \subfigure[Fixed Mesh, $N = 360\times 90 \times 4$]
    {\includegraphics[width=8cm]{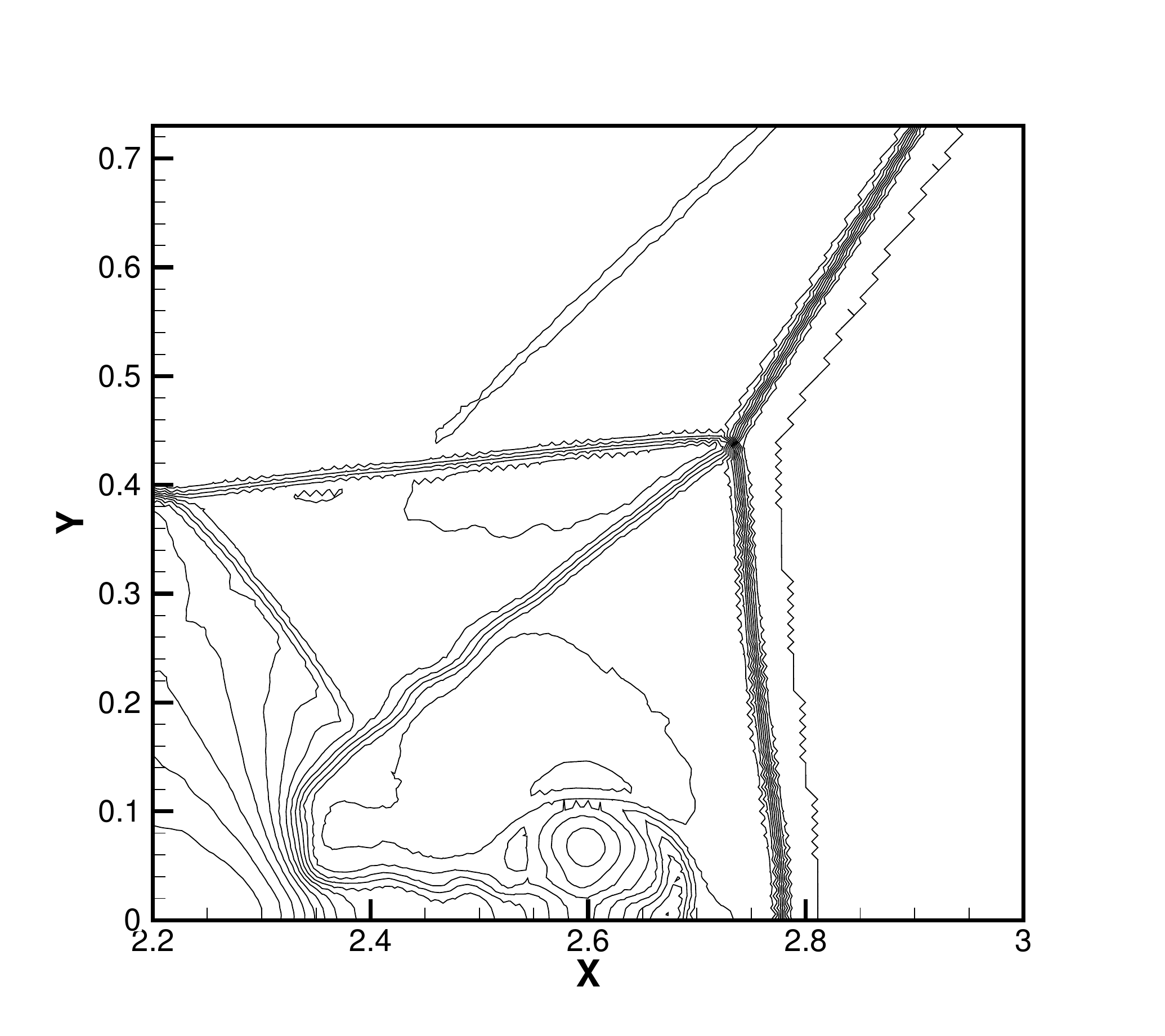}}
    }

  \mbox{\subfigure[Moving Mesh, $N = 240\times 60 \times 4$]
  {\includegraphics[width=8cm]{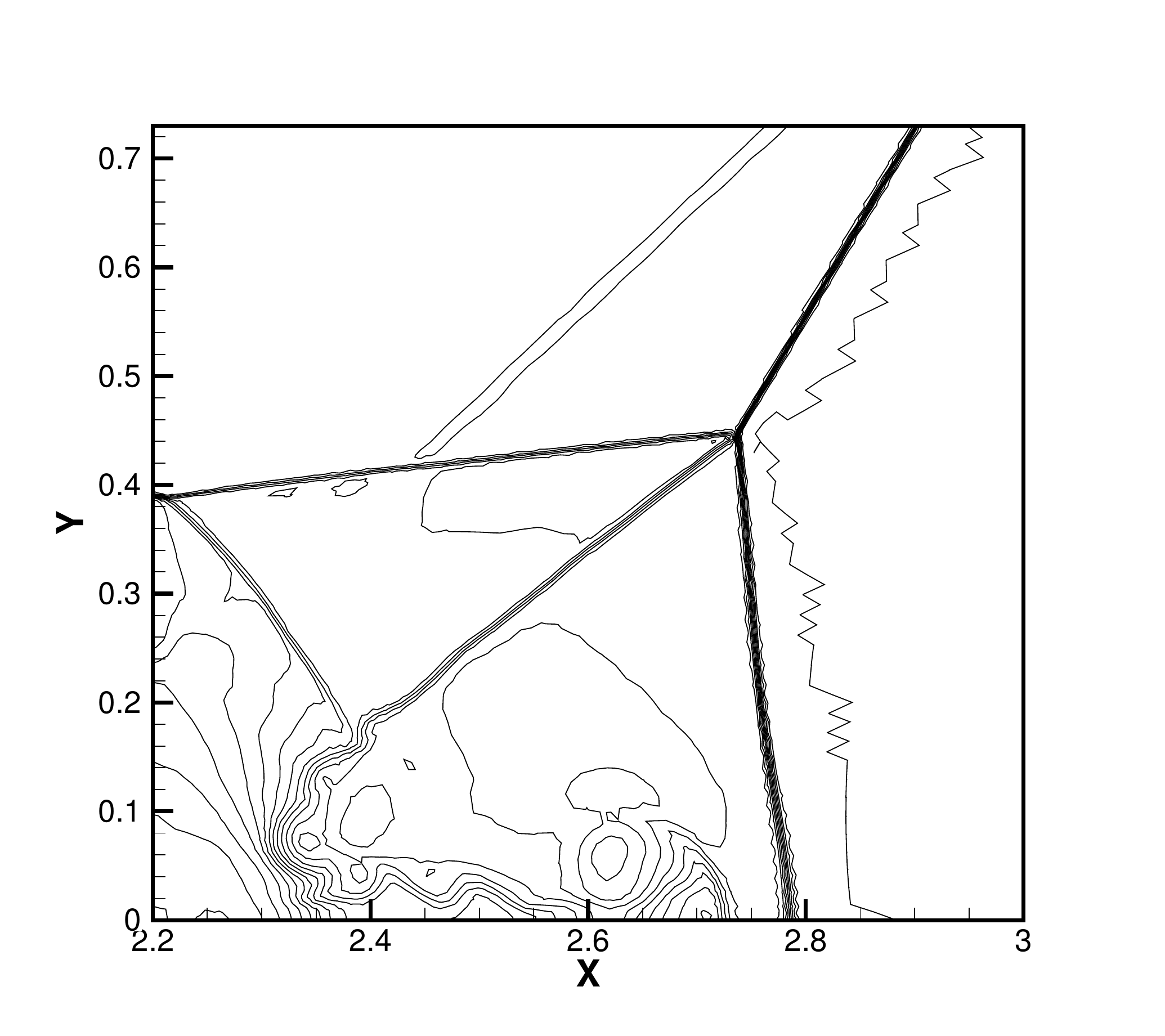}}\quad
    \subfigure[Moving mesh]
    {\includegraphics[width=8cm]{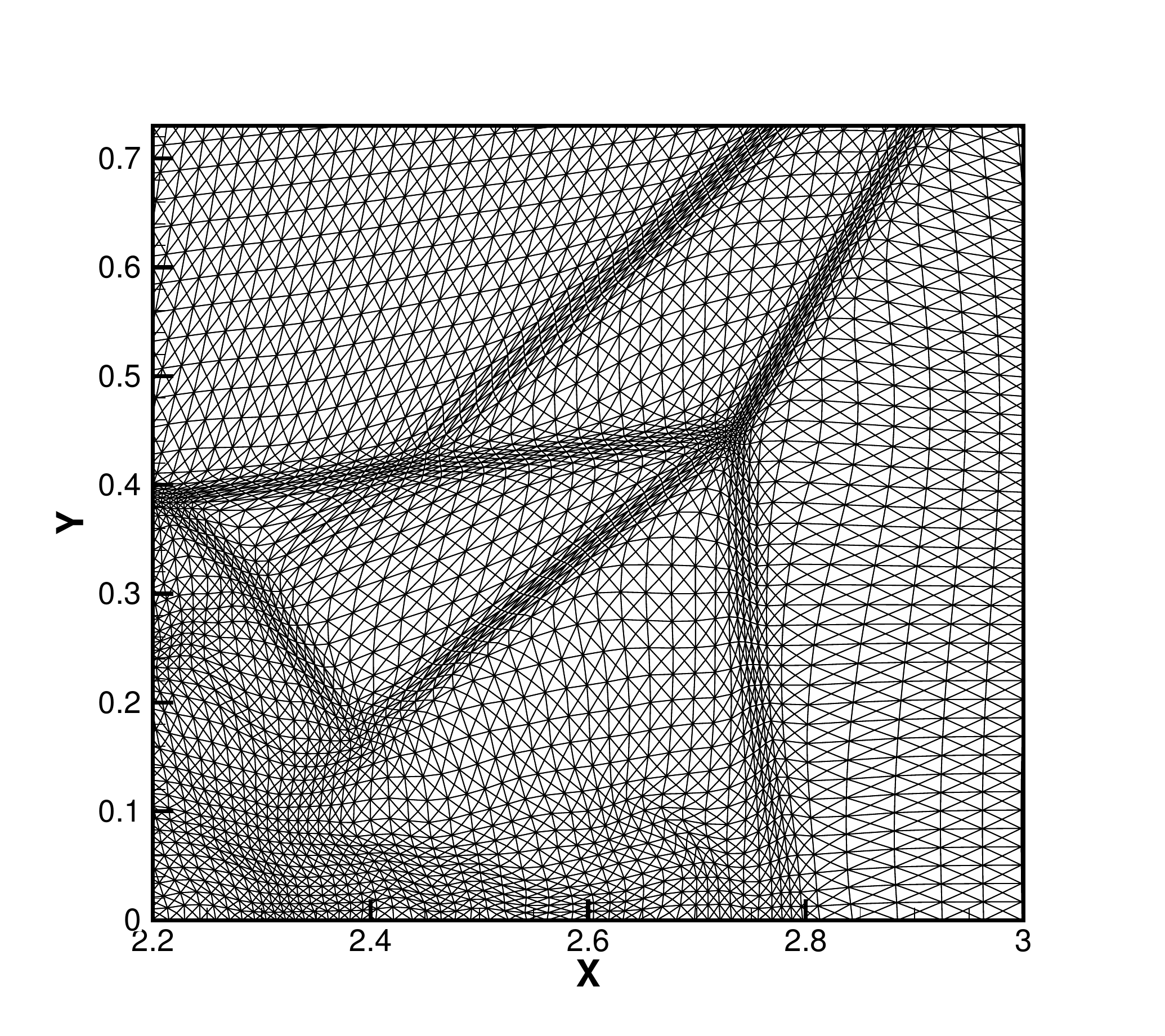}}
    }

    \caption{$P^2$ elements are used. Close view of the complex zone in Fig. \ref{fig:edge14}.}
    \label{fig:edge14zoom}
    \end{center}
    \end{figure}

\begin{exam}{\em
\label{exam4.9}
This is the double Mach reflection problem \cite{EH26}. We solve the Euler equations (\ref{2d}) in a computational domain of $(0,4)\times (0,1)$. The initial condition is given by 
\begin{equation}
u=
\left\{
\begin{array}{ll}
(8,57.1597,-33.0012,563.544)^T, \quad &\text{for} \quad  y\geq h(x,0)\\
(1.4,0,0,2.5)^T, \quad &\text{otherwise}
\end{array}
\right.
\notag
\end{equation}
where $h(x,t)=\sqrt{3}(x-\frac{1}{6})-20t$.
The exact post shock condition is imposed from $0$ to $\frac{1}{6}$ at the bottom while the reflection boundary condition for the rest of the bottom boundary. At the top, the boundary condition is the values that describe the exact motion of the Mach $10$ shock. On the left and right boundaries, the inflow and outflow boundary conditions are used, respectively. The final time is $T=0.2$. 

The density contours are shown in Figs. \ref{fig:edge11} - \ref{fig:edge14zoom} on $(0,3)\times (0,1)$. From these figures, one can see that more elements are concentrated in the regions where the shock and the complex structures are located. The result of the moving mesh method with $N = 120\times 30 \times 4$ 
is comparable with that obtained with $1.5$ times more uniform mesh points. The same is observed
with a moving mesh of $N = 240\times 60 \times 4$. 
 

 }\end{exam}

%
%

\begin{figure}[hbtp]
  \begin{center}
  \mbox{\subfigure[Fixed Mesh, $N = 120\times 40 \times 4$]
  {\includegraphics[width=8cm]{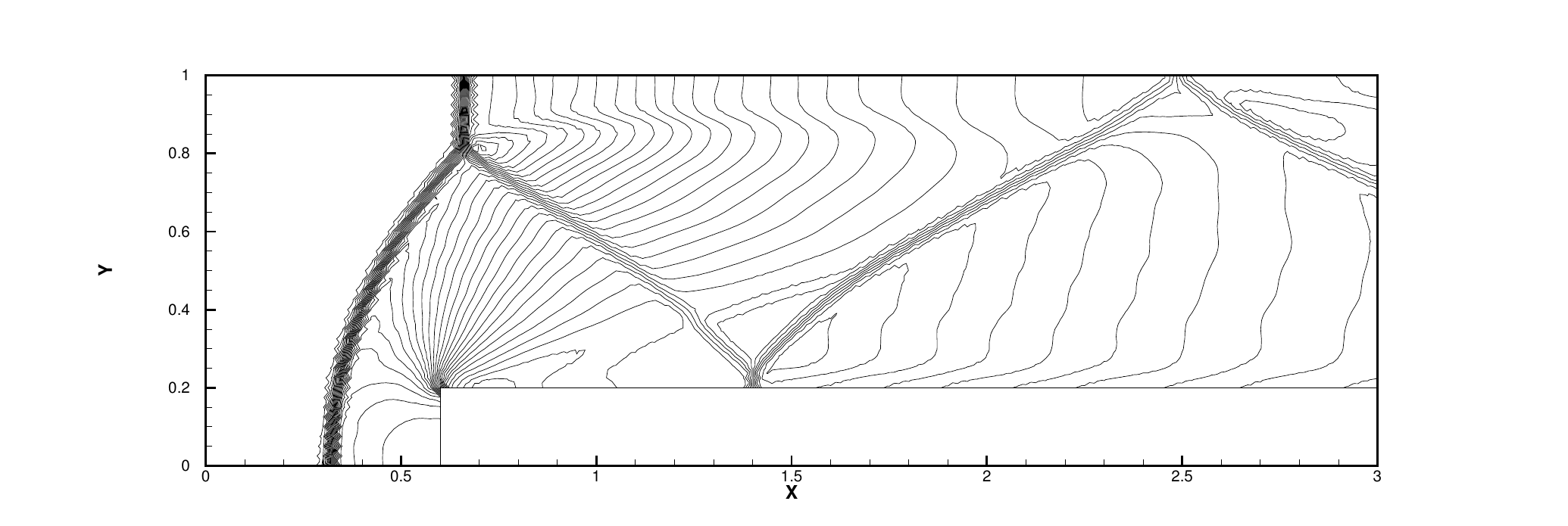}}\quad
    \subfigure[Fixed Mesh, $N = 240\times 80 \times 4$]
    {\includegraphics[width=8cm]{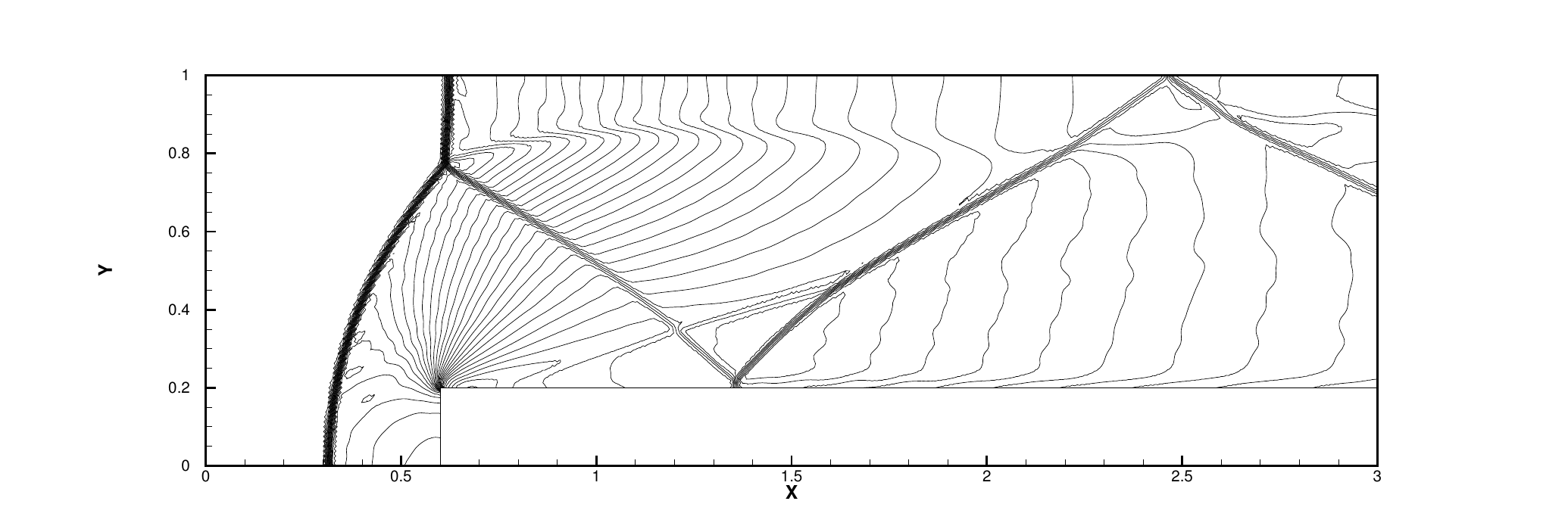}}
    }

  \mbox{\subfigure[Moving Mesh, $N = 120\times 40 \times 4$]
  {\includegraphics[width=8cm]{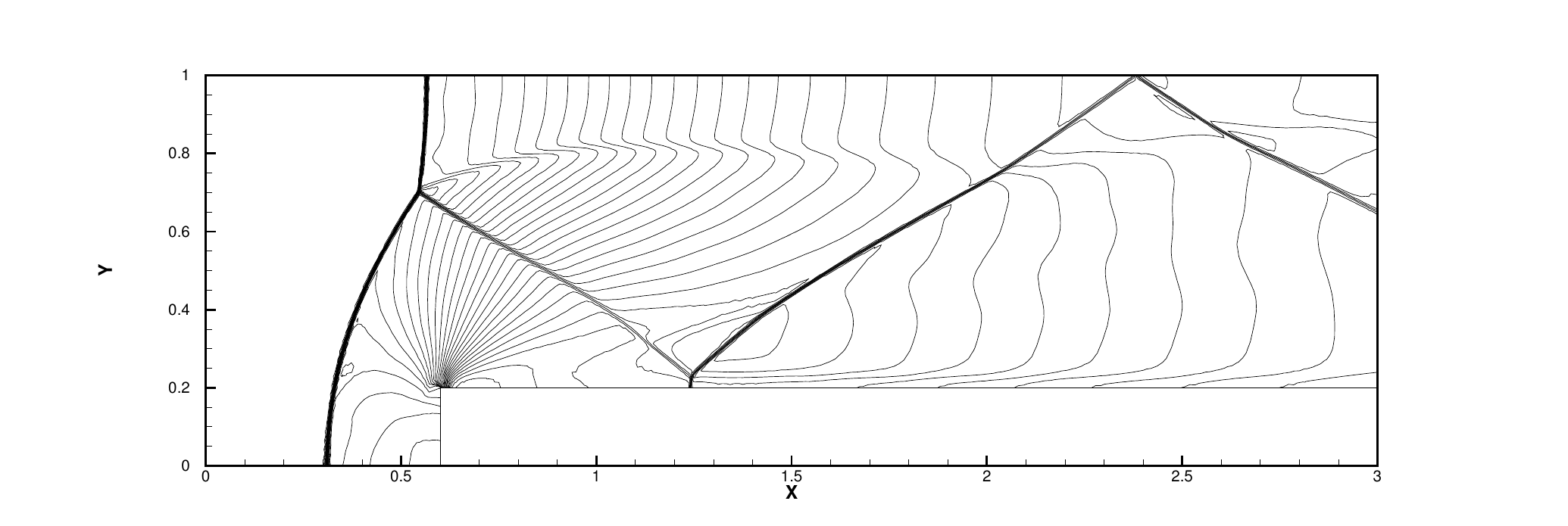}}\quad
    \subfigure[Moving mesh]
    {\includegraphics[width=8cm]{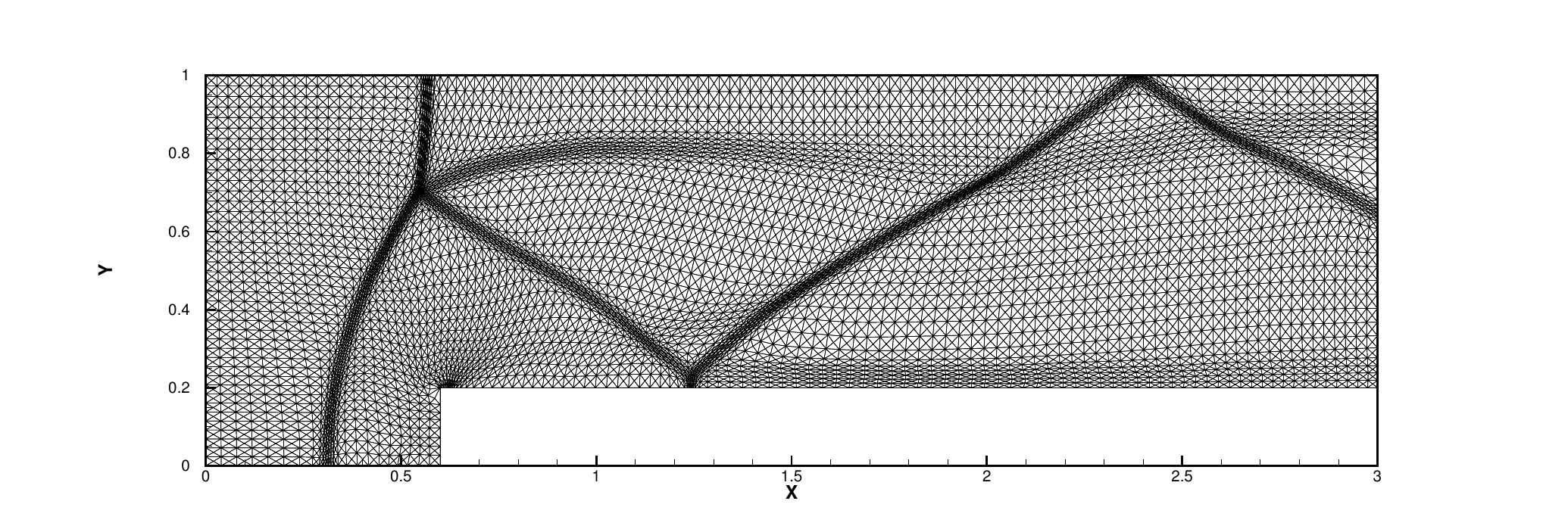}}
    }

    \caption{$P^1$ elements are used. 30 equally spaced density contours from 0.32 to 6.15 are used in the contour plots.}
    \label{fig:edge15}
    \end{center}
    \end{figure}

\begin{figure}[hbtp]
  \begin{center}
  \mbox{\subfigure[Fixed Mesh, $N = 120\times 40 \times 4$]
  {\includegraphics[width=8cm]{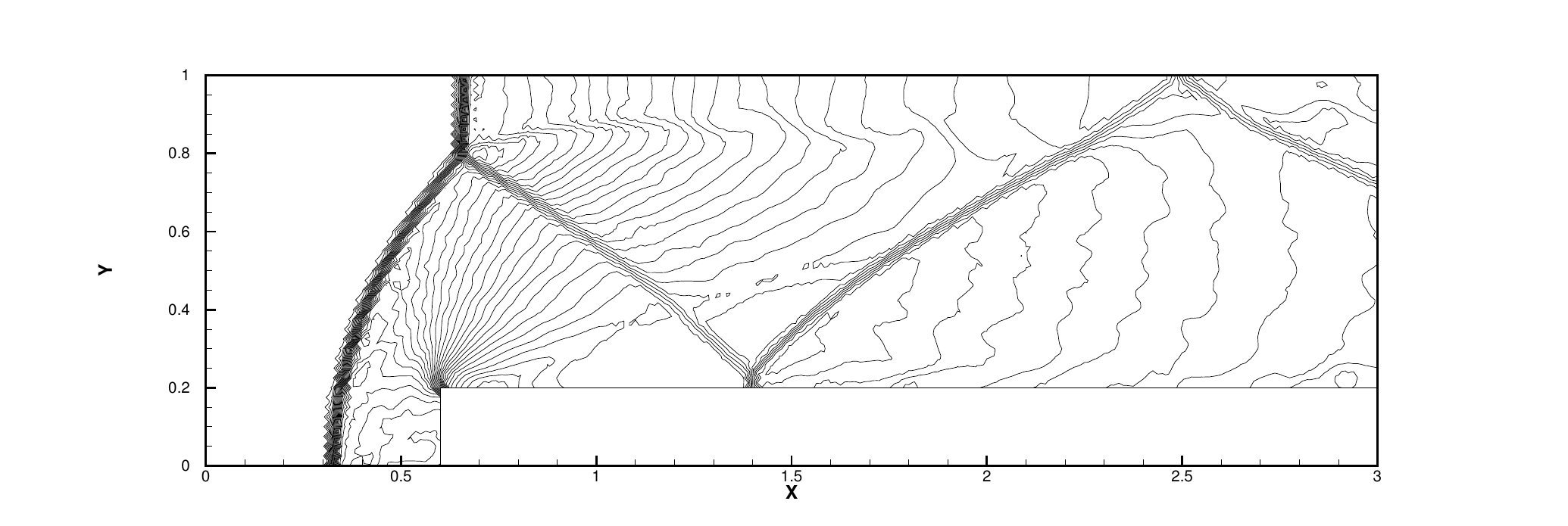}}\quad
    \subfigure[Fixed Mesh, $N = 240\times 80 \times 4$]
    {\includegraphics[width=8cm]{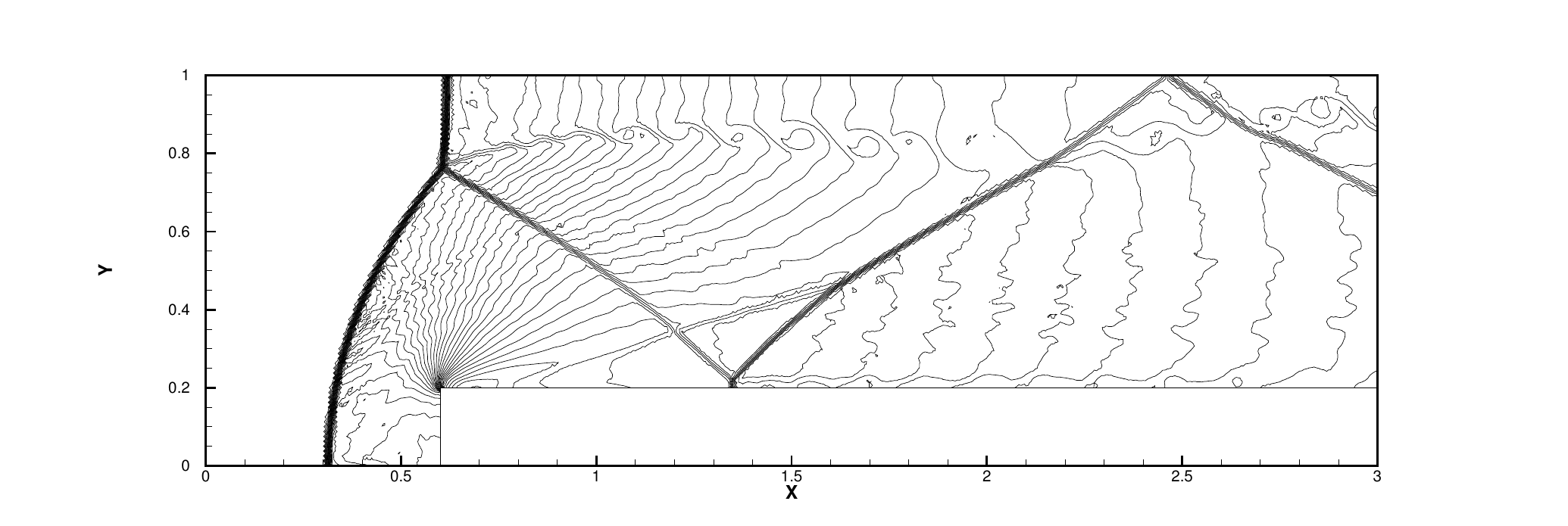}}
    }

  \mbox{\subfigure[Moving Mesh, $N = 120\times 40 \times 4$]
  {\includegraphics[width=8cm]{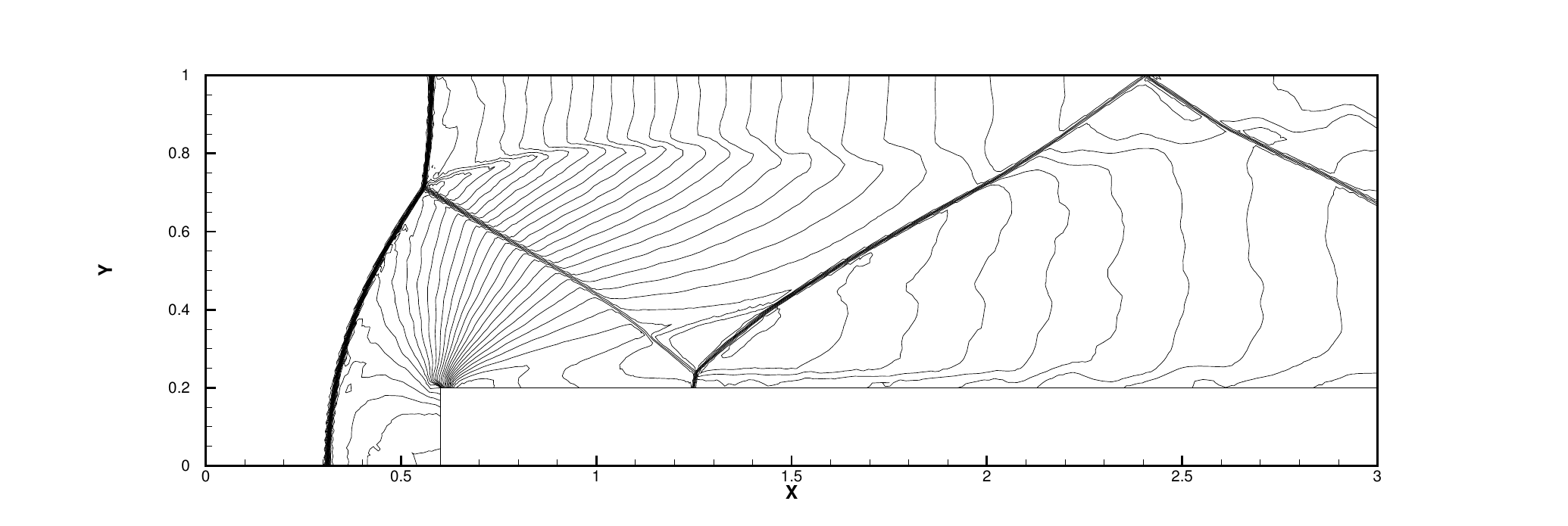}}\quad
    \subfigure[Moving mesh]
    {\includegraphics[width=8cm]{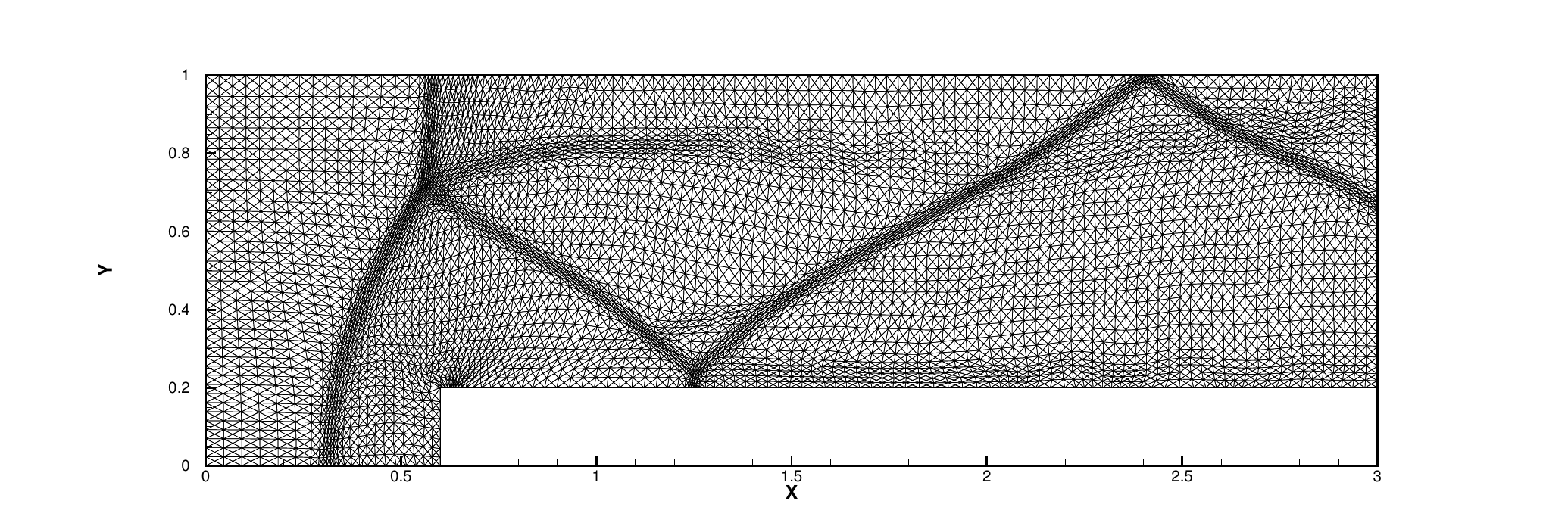}}
    }

    \caption{$P^2$ elements are used. 30 equally spaced density contours from 0.32 to 6.15 are used in the contour plots.}
    \label{fig:edge16}
    \end{center}
    \end{figure}

\begin{exam}{\em
\label{exam4.10}
The last example is the forward step problem \cite{EH26}. We solve the Euler equations (\ref{2d}) in a computational domain of $(0,3)\times (0,1)$. The problem is set up as follows: the wind tunnel is $1$ unit wide and $3$ units long. The step is $0.2$ units high and is located $0.6$ units from the left-hand end of the tunnel. The problem is initialized by a right-going Mach $3$ flow, namely,
\begin{equation}
(\rho,\mu,\nu,p)=(1.4,3,0,1).
\notag
\end{equation}
 Reflective boundary conditions are employed along the wall of the tunnel and inflow and outflow boundary conditions are applied at the entrance and exit, respectively. The final time is $T=4$.
 
The density contours from $0.32$ to $6.15$ are plotted in Figs. \ref{fig:edge15} and \ref{fig:edge16}. One can clearly see that the numerical solutions with the moving mesh method have better resolution than those with a uniform mesh of the same number of points.

}\end{exam}
 
\section{Conclusions}
\label{sec4}
\setcounter{equation}{0}
\setcounter{figure}{0}
\setcounter{table}{0}

In the previous sections we have presented a moving mesh DG method for the numerical solution
of hyperbolic conservation laws. The mesh is moved using the MMPDE moving 
mesh strategy where the nodal mesh velocities are defined as the gradient system of
an energy function associated with mesh equidistribution and alignment.
Moreover, hyperbolic conservation laws are discretized on a moving mesh in the quasi-Lagrangian fashion
with which the mesh movement is treated continuously and thus no interpolation is needed for physical variables
from the old mesh to the new one. Furthermore, the mesh movement introduces two extra convection terms
in the finite element formulation of conservation laws and their discretization can be incorporated
into the DG discretization naturally.

The numerical results for a selection of one- and two-dimensional examples have been presented.
They show that the moving mesh DG method achieves the theoretically predicted order of convergence for problems
with smooth solutions and is able to capture shocks and concentrate mesh points in non-smooth regions.
Moreover, it is shown that the numerical solution with a moving mesh is generally more accurate than that with a uniform
mesh of the same number of points and often comparable with the solution obtained with a much finer uniform mesh.
Furthermore, numerical results for problems with smooth and discontinuous solutions in one and two dimensions
have shown that the accuracy of the method is not sensitive to the smoothness of the mesh, which is in contrast
with the situation for the moving mesh finite difference WENO method \cite{EH08}.

We recall that the Hessian of a physical variable has been used in this work to guide the mesh adaptation
(cf. (\ref{mer}).
This is based on linear interpolation error \cite{EH30} and has been known to work for many problems.
Nevertheless, it would be advantageous to define the metric tensor based on some a posteriori error estimate.
Investigations of using residual-based metric tensors for the moving mesh DG method presented
in this work have been underway.

\section*{Acknowledgements}
The work was supported in part by China NSAF grant U1630247 and NSFC grant 11571290.
D. M.  Luo gratefully acknowledges the financial support from  China Scholarship Council under the grant 201506310112 to visit the Department of Mathematics, University of Kansas from September 2015 to September 2017. He also thanks Professor Jun Zhu at Nanjing University of Aeronautics and Astronautics for his help with the HWENO limiters.

\end{document}